\documentclass[12pt]{article}
\usepackage{amsmath}
\usepackage{amssymb}
\usepackage{amscd}
\usepackage{epsfig}
\usepackage{amsthm}
\usepackage{mysects}

\def\no{\noindent}
\numberwithin{equation}{subsection}

\theoremstyle{plain}
\newtheorem{theorem}[equation]{Theorem}

\newtheorem{proposition}[equation]{Proposition}
\newtheorem{lemma}[equation]{Lemma}
\newtheorem{corollary}[equation]{Corollary}

\theoremstyle{remark}
\newtheorem{remark}[equation]{Remark}

\theoremstyle{definition}
\newtheorem{definition}[equation]{Definition}

\newtheorem{example}[equation]{Example}

\newcommand{\lra}{\longrightarrow}
\newcommand{\ra}{\rightarrow}
\newcommand{\restr}{\mbox{\Large \(|\)\normalsize}}
\newcommand{\N}{\mathbb N}
\newcommand{\R}{\mathbb R}

\newcommand{\E}{\mathbb E}
\newcommand{\Z}{\mathbb Z}

\newcommand{\tits}{\partial_{T}}

\newcommand{\acts}{\curvearrowright}
\newcommand{\dcup}{\amalg}
\newcommand{\<}{\langle}
\renewcommand{\>}{\rangle}

\def\cangle{\widetilde\angle}

\newcommand{\al}{\alpha}
\newcommand{\be}{\beta}
\def\de{\delta}
\def\De{\Delta}

\def\eps{\epsilon}
\def\ga{\gamma}

\def\la{\lambda}

\def\lra{\longrightarrow}

\def\ol{\overline}

\def\ra{\rightarrow}

\def\si{\sigma}
\def\Si{\Sigma}

\def\th{\theta}

\def\geo{\partial_{\infty}}
\def\tits{\partial_{T}}

\def\tangle{\angle_{T}}
\def\defeq{:=}
\def\eqdef{=:}

\newcommand{\Tm}{{\cal T}}

\newcommand{\St}{{\cal S}}

\newcommand{\itin}{{\cal I}}
\newcommand{\hd}{d_{\cal H}}

\begin{document}

\title{The geodesic flow of a nonpositively curved graph manifold}
\author{Christopher B. Croke\thanks{Supported by NSF grants DMS-95-05175 DMS-96-26-232 and DMS 99-71749.}\\
Bruce Kleiner\thanks{Supported by a Sloan Foundation
Fellowship, and NSF grants
DMS-95-05175, DMS-96-26911, DMS-9022140.}}
\date{\today}
\maketitle

\begin{abstract}
We consider discrete cocompact isometric actions $G\stackrel{\rho}{\acts}X$
where $X$ is a locally compact Hadamard space\footnote{Following \cite{ballmann}
we will refer to $CAT(0)$ spaces (complete, simply connected length spaces
with nonpositive curvature in the sense of Alexandrov) as Hadamard spaces.}, 
and $G$ belongs
to a class of groups (``admissible groups'') which includes fundamental groups
of $3$-dimensional graph manifolds.
We identify invariants (``geometric data'') of the action $\rho$ which
determine, and are determined by, the equivariant homeomorphism type of the
action $G\stackrel{\geo\rho}{\acts}\geo X$ of $G$ on the
ideal boundary of $X$.   
Moreover, if $G\stackrel{\rho_i}{\acts}X_i$
are two actions with the same geometric data
 and $\Phi:X_1\ra X_2$ is a $G$-equivariant
quasi-isometry, then for every geodesic ray $\ga_1:[0,\infty)\ra X_1$,
there is a geodesic ray $\ga_2:[0,\infty)\ra X_2$ (unique up to equivalence)
so that $\lim_{t\ra\infty}\frac{1}{t}d_{X_2}(\Phi\circ\ga_1(t),
\ga_2([0,\infty)))=0$.  This work
was inspired by (and answers) a question of Gromov in
\cite[p. 136]{asyinv}.
\end{abstract}

\setcounter{section}{1}
\setcounter{subsection}{0}

\subsection{Introduction}

As a consequence of the Morse lemma 
on quasi-geodesics,
geodesic flows are especially simple
and well understood in the Gromov hyperbolic \mbox{case\hspace{2pt}:}

a. If $\phi:M_1\ra M_2$ is a homotopy equivalence
between closed negatively curved manifolds, then there
is an orbit equivalence $\hat\phi:SM_1\ra SM_2$ between
the unit sphere bundles, which covers $\phi$ up to homotopy \cite{3remarks}. 

b.  If $G$ is a hyperbolic group,  $G\stackrel{\rho_i}{\acts}X_i$
is a discrete, cocompact, isometric action on a Hadamard space
$X_i$ for $i=1,\,2$, and $\Phi:X_1\ra X_2$ is a $G$-equivariant
quasi-isometry, then $\Phi$ maps each geodesic $\ga_1\subset X_1$
to a subset at uniformly bounded Hausdorff distance from
a geodesic $\ga_2\subset X_2$.  Moreover, $\Phi$ induces
an equivariant homeomorphism $\geo\Phi:\geo X_1\ra\geo X_2$
between ideal boundaries, \cite{hypgps}.

c.  When $G\stackrel{\rho}{\acts}X$ is 
a discrete, cocompact action of a hyperbolic group
on a Hadamard space $X$, then
the induced action $G\stackrel{\geo\rho}{\acts}\geo X$
of $G$ on the boundary of $X$ is a finitely presented 
dynamical system, \cite{hypgps,codepa}. 

\medskip
\no
Naturally one may ask if  properties b and c hold
without the assumption of Gromov hyperbolicity. 
It turns out that they do not:  one can readily
produce examples of  pairs of discrete, cocompact,
isometric actions $G\acts X_1$, $G\acts X_2$ where $G$-equivariant 
quasi-isometries $X_1\ra X_2$ do not induce boundary 
homeomorphisms\footnote{Let $M_1$ and $M_2$
be closed surfaces with nonpositive curvature, and let $N_1$ and $N_2$
be the Riemannian products $N_i\defeq M_i\times S^1$.  Suppose $f_0:M_1\ra M_2$
is a homotopy equivalence, $f\defeq f_0\times id_{S^1}:N_1\ra N_2$
is the corresponding map between the $N_i$'s,
and $\hat f:\tilde N_1\ra \tilde N_2$ is a lift of $f$ to a map 
between the universal covers.
Then it turns out that $\hat f$ extends continuously to up to the ideal boundary
$\geo N_1$ if and only if $f_0$ is homotopic to a homothety.}
(this was observed independently by Ruane \cite{ruane}).
In  \cite{asyinv}
Gromov asked whether two actions $G\acts X_i$ 
induce $G$-equivariantly homeomorphic boundary actions
$G\acts \geo X_i$.  The answer to this is also no: S. Buyalo \cite{buy} and the authors
independently found pairs of actions which induce  inequivalent boundary
actions\footnote{Boundaries can even fail to be (non-equivariantly)
homeomorphic: \cite{leebcx} describes a pair of homeomorphic
nonpositively curved $2$-complexes whose universal covers have
nonhomeomorphic boundary (see also \cite{jwilson}).}.   Finally, we remark that the boundary 
action $G\acts\geo X$ is finitely presented if and only if $G$ is hyperbolic\footnote{
The action $G\acts\geo X$ is expansive if and only if $G$
is hyperbolic.}.

In this paper we examine actions $G\acts X$ where $G$
belongs to a class of groups which generalize 
fundamental groups of $3$-dimensional graph manifolds.
We develop a kind of ``coding'' for geodesic rays in $X$,
which allows us to understand the boundary action
$G\acts\geo X$ and the Tits metric on $\geo X$.
Before stating our main result in complete generality,
we first formulate it for
nonpositively curved $3$-dimensional graph manifolds.  

By the theorem of 
\cite{schroeder},  when $M$ is a $3$-dimensional graph manifold with
a nonpositively curved Riemannian metric, then $M$ has the
following structure.  There is a collection $M_1,\ldots,M_k$
of compact nonpositively curved $3$-manifolds with nonempty totally geodesic
boundary (the geometric Seifert components of $M$), and Seifert fibrations 
$M_i\stackrel{p_i}{\ra}N_i$ where
the metric on $M_i$ has local product structure compatible
with the fibration $p_i$, and the $N_i$ are nonpositively 
curved orbifolds;  $M$ is obtained from the
disjoint union $\dcup_i \,M_i$ by gluing boundary components isometrically
in pairs via gluing isometries which are incompatible with the boundary
fiberings.  
In what follows we  will only consider graph manifolds whose Seifert
fibered components have orientable fiber.  
Note that for each $1\leq i\leq k$,
the universal cover of $M_i$ is isometric to a Riemannian product $\tilde N_i\times\R$;
the action of $\pi_1(M_i)$  on $\tilde M_i$ 
preserves this product structure and so there is an induced action
of $\pi_1(M_i)$ on the $\R$ 
factor by translations. Hence we get a homomorphism
$\tau_i:\pi_1(M_i)\ra\R$ for each $i$.   We may also define a class function
$MLS_i:\pi_1(M_i)\ra \R_+$ by taking the minimum of the
displacement function for the induced action $\pi_1(M_i)\acts\tilde N_i$,
i.e. $MLS_i(g)=\inf\{d_{\tilde N_i}(gx,x)\mid x\in\tilde N_i\}$; this corresponds to
the marked length spectrum of the nonpositively curved
orbifold $N_i$.

Now suppose  $M$ and $M'$ are graph manifolds as above,
and $f:M\ra M'$ is a homotopy equivalence.  Embedded incompressible tori in Haken
manifolds are determined up to isotopy by their fundamental groups
up to conjugacy \cite{laudenbach},
so we may  assume after
isotoping $f$ that it is a homeomorphism which  induces  
homeomorphisms $ f_i:M_i\ra M_i'$ from the Seifert
components  of $M$ to the
Seifert components  of $M'$ (and hence isomorphisms on the corresponding fundamental groups).  
We may then use the maps $ f_i$ to compare
the invariants $\tau_i,\,\tau_i'$ and $MLS_i,\,MLS_i'$.

\begin{theorem}
\label{graphmfldcase}
The following are equivalent:

1. The functions $MLS_i$ and $\tau_i$
are preserved up to scale by $f_i$: for $i=1,\ldots,k$ there are
constants $a_i$ and $b_i$ so that 
$MLS_i=a_if^*(MLS_i')$ and $\tau_i=b_if^*(\tau_i')$.

2. Any lift $\tilde f:\tilde M\ra \tilde M'$ of $f$ 
extends continuously to a 
map $\bar f:M\cup\geo M\ra M'\cup\geo M'$
between the standard compactifications.

3. If $\tilde f:\tilde M\ra \tilde M'$ is any lift of $f$,
then $\tilde f$ maps geodesic rays to geodesic rays, 
up to uniform sublinear error:
there is a function $\th:\R_+\ra\R_+$
with $\lim_{r\ra\infty}\th(r)=0$ so that if $\ga:[0,\infty)\ra \tilde M$
is a unit speed geodesic ray, then there is a ray $\ga':[0,\infty)\ra \tilde M'$
where $d(\tilde f\circ\ga(t),\ga'([0,\infty)))<(1+t)\th(t)$.

4. If we identify $\pi_1(M)$ with $\pi_1(M')$ via $f$, then
the induced boundary actions $\pi_1(M)\acts \geo\tilde M$
and $\pi_1(M)\acts \geo\tilde M'$ are equivariantly homeomorphic
(by a unique equivariant homeomorphism).

\medskip
\no
If 1 holds and in addition the constants $a_i$ and $b_i$ are independent of
$i$, then the unique equivariant homeomorphism $\geo\tilde M\ra
\geo\tilde M'$ in 4 is an isometry with respect to Tits metrics.
\end{theorem}
In general (see Lemma \ref{constants}) the structure of $\pi_1(M)$ forces the $a_i$'s and
$b_i$'s in condition 1 to satisfy $\#\{a_1,\ldots,a_k,b_1,\ldots,b_k\}\leq 2$,
and except in special circumstances they all coincide.  
The condition $MLS_i=a_if^*(MLS_i')$
means that the homotopy equivalence $N_i\ra N_i'$
induced by $f_i$ preserves the marked length spectrum
of the nonpositively curved orbifolds
up to the scale factor $a_i$.
Although closed nonpositively curved surfaces with the
same marked length spectrum are isometric by
\cite{croke,otal,croke-fathi-feldman},
compact nonpositively curved surfaces with geodesic boundary
can have the same
marked length spectrum without being isometric:

\begin{example}
\label{bikini}
Let $N$ be a pair of pants with a 
 (constant curvature $-1$) hyperbolic
metric where the boundary
components are geodesics with length $L$,
and let $\{c_1,c_2\}\subset N$ be
the fixed point set of the order $3$ isometry of 
$N$.  If
$L$ is sufficiently large (so that $N$ looks like a bikini)
then a closed geodesic in $N$ cannot pass near 
 $\{c_1,c_2\}$.  This means that one can change the 
metric near $\{c_1,c_2\}$ without disturbing the 
marked length spectrum of $N$.  Note that one can modify this 
example slightly so that the metric is flat in a neighborhood of the boundary geodesics.
\end{example}
\noindent
Suppose $M$ is a nonpositively curved graph manifold 
with a Seifert
component $M_i$ isometric to $N\times S^1$, where $N$
is as in the example.  One can
change the metric on the $N$ factor as in example
\ref{bikini} to get a Riemannian 
manifold $M'$  so that the conditions of Theorem
\ref{graphmfldcase} hold (with $f=id$), but $M'$ is not isometric to $M$.

In section \ref{examplesection} we give an example to show  that the uniform
sublinear divergence estimate in condition 3 
cannot be improved to a bounded distance estimate
as in the Gromov hyperbolic case.

We now sketch some of the main points in the proof
of Theorem \ref{graphmfldcase}.  

First consider
a single nonpositively curved graph manifold $M$
with geometric Seifert components $M_1,\ldots,M_k$.
The universal cover $\tilde M_i$ is isometric to $\tilde N_i\times \R$ --
a nonpositively curved $3$-manifold with a countable 
collection of totally geodesic boundary components
isometric to $\E^2$.
The universal cover $\tilde M$ of $M$ is tiled by a countable
collection of copies of the  universal covers $\tilde M_i$
for $i=1,\ldots k$; we call these subsets {\em vertex spaces}. 
We refer to boundary components of vertex spaces as
{\em edge spaces}.  
Two vertex spaces are either disjoint, or intersect along
an edge space. 
Let $T$ be the incidence graph
for the collection of vertex spaces: $T$ is the graph which
has one vertex for each vertex space, and an edge joining
two vertices whenever the corresponding vertex spaces intersect.
$T$ is isomorphic to the Bass-Serre tree of the graph of 
groups associated with the decomposition $\dcup_i M_i\ra M$ (see section \ref{basicbassserre}).
If $v\in V\defeq Vertex(T)$ (resp. $e\in E\defeq Edge(T)$)
we will use the notation $\tilde M_v$ (resp. $\tilde M_e$)
for the vertex space (resp. edge space) associated with $v$
(resp. $e$); and we let $\tilde M_v\simeq \tilde N_v\times\R$
be the Riemannian product decomposition of $\tilde M_v$.  It
is not difficult to check (Lemma \ref{fixedinbdy}) 
that if  $G_v\defeq Stabilizer(\tilde M_v)\subset G\equiv \pi_1(M)$, 
then the center $Z(G_v)$ of $G_v$ is isomorphic to
$\Z$, and the fixed point set of $Z(G_v)$ in $\geo\tilde M$
is just $\geo \tilde M_v$; similarly, if $\tilde M_e$
is an edge space then the fixed point set of 
$\Z^2\simeq G_e\defeq Stabilizer(\tilde M_e)\subset G$ in $\geo \tilde M$
is $\geo \tilde M_e$. 

Let $p\in \tilde M$ be an interior point of  a
vertex space,  pick $\xi\in\geo M$, and let
$\ol{p\xi}$ denote the geodesic ray starting at
$p$ which is asymptotic to $\xi$.  The
ray $\ol{p\xi}$  encounters a (possibly finite) sequence 
of vertex and edge spaces called the {\em itinerary}
of $\ol{p\xi}$.  The convexity of vertex and
edge spaces forces the itinerary 
$v_0,e_1,v_1,e_2,\ldots$ of $\ol{p\xi}$
to be the sequence
of successive vertices and edges of a geodesic segment
or ray in $T$.  In order to understand the rays with itinerary
$v_0,e_1,v_1,e_2,\ldots$, we construct a piecewise flat complex $\Tm$ --
a {\em template} -- in $\tilde M$ as follows.  First let
$\ga_i\subset \tilde M_{v_i}$ be a shortest geodesic from 
$\tilde M_{e_{i}}$ to $\tilde M_{e_{i+1}}$ for $i>0$,
and let $\ga_0$
be a shortest path from $p\in\tilde M_{v_0}$ to $\tilde M_{e_1}$.
For $i\geq 0$ define $\St_i\subset \tilde M_{v_i}$ to be the flat
strip which is the
union of the geodesics in $\tilde M_{v_i}$
which are parallel to the $\R$-factor of $\tilde M_{v_i}$ and
which pass through $\ga_i$.
We define $\Tm$ to be 
the union of the edge spaces $\{\tilde M_{e_i}\}$
with the strips $\{\St_i\}$;    then
$\Tm$ is a Hadamard space with respect to the induced path
metric.  A key technical step in the proof of  Theorem \ref{graphmfldcase}
is Theorem \ref{shadthm}, which shows that for any geodesic ray
 $\ol{p\zeta}$ in the Hadamard space
$\Tm$, there is a  unique geodesic ray $\ol{p\zeta'}$
in $\tilde M$ with the property that for all 
 $x\in \ol{p\zeta}$, 
$$d_{\tilde M}(x,\ol{p\zeta'})\leq
\th(d_\Tm(x,p))(1+d_\Tm(x,p))$$
for some function $\th:[0,\infty)\ra[0,\infty)$ with $\lim_{r\ra\infty}\th(r)=0$
which is independent of the choice of itinerary.  
Using Theorem \ref{shadthm} one finds that the
set  of boundary points $\zeta\in\geo\tilde M$
for which  the ray $\ol{p\zeta}$ has a given infinite itinerary $v_0,e_1,\ldots$ 
is homeomorphic to the set boundary points $\zeta'\in\Tm$ so that the 
$\Tm$-ray $\ol{p\zeta'}$ passes through $\tilde M_{e_i}$ for every $i$.
One sees (Proposition \ref{interval}) that the latter is either
a single point or is homeomorphic to a closed interval, depending on the 
geometry of $\Tm$ (which depends, in turn, on the choice of itinerary and the geometry of $M$).

We now consider a second nonpositively curved $3$-manifold $M'$,
and use primes to denote the vertex spaces, edge spaces, etc for $M'$.   
Let $f:M\ra M'$ be a homeomorphism as in Theorem \ref{graphmfldcase},
and identify the deck groups $G\defeq \pi_1(M)\simeq\pi_1(M')$ via
a lift $\tilde f$ of $f$. Then $\tilde f$ maps vertex (resp. edge)
spaces of $M$ homeomorphically to vertex (resp. edge)
spaces of $M'$, so we may use $\tilde f$ to identify the 
incidence tree $T'$ with $T$.   Suppose $\phi:\geo\tilde M\ra\geo\tilde M'$ is a 
$G$-equivariant homeomorphism.  Using the remarks about fixed point
sets made above, it follows
that $\phi(\geo\tilde M_v)=\geo\tilde M'_v$ and 
$\phi(\geo\tilde M_e)=\geo\tilde M'_e$ for every $v\in V$
and every $e\in E$.   Also, if $p\in \tilde M_{v_0}$, 
and $v_0,e_1,v_1,e_2,\ldots$ is an infinite itinerary, then
$\phi(S)=S'$ where $S\subset\geo \tilde M$ and $S'\subset\geo\tilde M'$ 
are the subsets corresponding to the itinerary $v_0,e_1,v_1,e_2,\ldots$ 
(Corollary \ref{preslabelling});
in particular, either $S$ and $S'$ are both points or they are both intervals.
By considering all possible infinite itineraries and exploiting
this correlation, we are able to see (section \ref{recovering})
that the invariants
$MLS_i,\,MLS_i'$ and $\tau_i,\,\tau_i'$ must agree as in condition 1 of
Theorem \ref{graphmfldcase}.  Conversely, if condition 1 holds
and $p\in \tilde M$, one shows (section \ref{samedata}) 
that for each itinerary the 
corresponding templates in $\tilde M$ and  $\tilde M'$
have sufficiently similar geometry that their geodesics are ``similar'';
and this implies that $\tilde f$ extends to the compactifications
as in 2 of Theorem \ref{graphmfldcase}.

\medskip
Our main result generalizes Theorem \ref{graphmfldcase}
and applies to {\em admissible groups}, a class of 
(fundamental groups of) graphs of groups, see section
\ref{admissible} for the precise definition. When 
 an admissible group $G$ acts discretely and cocompactly
on a Hadamard space $X$
then we associate geometric data to each vertex group
$G_v\subset G$ consisting of a class function $MLS_v:G_v\ra \R^+$
and a homomorphism $\tau_v:G_v\ra \R$ (see section \ref{sectionvertexandedge}).

\begin{theorem}
\label{main}
Let $G\acts X$ be a discrete, cocompact, isometric action of 
an admissible group on  a Hadamard space $X$.  
Then for every vertex $v$, $MLS_v$ and $\tau_v$
are determined up to scale factors $a_v$ and $b_v$
 by the topological 
conjugacy class of the boundary action $G\acts\geo X$,
and vice-versa.   If $G\acts X'$ is another such
action, then the following are equivalent:

1.  $G\acts X$ and $G\acts X'$ have the same geometric 
data up to scale.

2. $G$-equivariant quasi-isometries 
$X\ra X'$ extend canonically to the compactifications
$X\cup\geo X\ra X'\cup\geo X'$.

3.  The boundary actions $G\acts\geo X$ and $G\acts \geo X'$
are $G$-equivariantly homeomorphic (by a unique\footnote{It follows
from the methods of \cite[III.3]{ballmann} that if $G\acts X$ is a 
cocompact isometric action on a Hadamard space, $g\in G$
is an axial isometry,  and $\ga_g\subset$ is an axis for $g$  which does not
bound a flat half plane, then the orbit of $\geo\ga\subset \geo X$ under
the action  $G\acts\geo X$ is dense in $\geo X$.   Hence the set
of  points  in $\geo X$ which are the unique attracting fixed point of some
element of $G$  is dense in $\geo X$.   Any $G$-equivariant
homeomorphism $\geo X\ra \geo X$ must fix this dense set pointwise, and
must therefore be the identity.  The uniqueness  statement follows immediately
from this.} $G$-equivariant homeomorphism).  

4.  If $f:X\ra X'$ is a $G$-equivariant quasi-isometry, then
there is a function $\th:\R_+\ra\R_+$ with $\lim_{r\ra\infty}\th(r)=0$
so that for every unit speed geodesic ray $\ga:[0,\infty)\ra X$ 
there is a ray $\ga':[0,\infty)\ra X'$ with 
$d(f\circ\ga(t),\ga'(t))<(1+t)\th(t)$.

  Furthermore, if there is a
single scale factor $s$ so that $MLS_v'=sMLS_v$ and $\tau_v'=s\tau_v$
for every vertex $v$, then the unique  $G$-equivariant homeomorphism
$\geo X\ra\geo X'$ is an isometry with respect to the Tits metrics.
\end{theorem}

The authors proved  theorem \ref{main} while attempting to digest
the negative answer to Gromov's question about boundary actions.
A key factor in our example was the (unanticipated) presence of intervals
in the Tits boundary.   After  the examples and their
properties had been announced, similar structure was found in 
other manifolds, \cite{humsch}.   The paper \cite{buysch} also contains
some discussion of the Tits boundary of universal covers of nonpositively
curved graph manifolds.

\medskip
\no
{\bf Open questions.}
The results in this paper raise a number of questions.  First of
all, for each group $G$ one may ask for a generalization of Theorem
\ref{main}, where the geometric data  $MLS_v$ and $\tau_v$  are
replaced with suitable substitutes.   Our methods actually yield more
information about the behavior of geodesics than is stated in 
Theorem \ref{main} alone.   We are able to give a good
description of all the geodesic rays in the Hadamard space $X$ in terms
of concrete geometric information; it seems likely that other classes
of groups are amenable to a similar treatment. 
 The  fundamental groups of the real-analytic manifolds considered
in \cite{humsch} are natural candidates for this, as they
have structure similar to graph manifold groups.    Here are two other questions:

1. What determines the (non-equivariant) homeomorphism
type of $\geo X$, when $X$ is  a Hadamard space with 
an action $G\acts X$ by an admissible group $G$?

2.  Does part 4 of Theorem \ref{main} have an analog
where rays are replaced by complete geodesics?  This seems
within reach.

\bigskip

\tableofcontents

\subsection{Preliminaries}

\subsubsection{Coarse geometry}

Let $X$ and $X'$ be metric spaces, and let $\Phi:X\ra X'$ be a map.

\begin{definition}
1. $\Phi$ is {\em $(L,A)$-Lipschitz} if for all $x_1,x_2\in X$,
$$d(\Phi(x_1),\Phi(x_2))\leq Ld(x_1,x_2)+A.$$
$\Phi$ is {\em coarse Lipschitz} if it is $(L,A)$-Lipschitz for some
$L,\,A>0$.

2.  $\Phi$ is an {\em $(L,A)$-quasi-isometric embedding} if it is
$(L,A)$-Lipschitz and for all $x_1,x_2\in X$,
$$d(\Phi(x_1),\Phi(x_2))\geq L^{-1}d(x_1,x_2)-A.$$
The constants $(L,A)$ will often be suppressed.  
A {\em quasi-geodesic} (respectively {\em segment/ray}) 
is a quasi-isometric embedding $\Phi:\R\ra X$ 
(respectively $\Phi:[a,b]\ra X$, $\Phi:[0,\infty)\ra X$).
We sometimes refer to the image of a quasi-geodesic 
as a quasi-geodesic.

3.  $\Phi$ is an {\em $(L,A)$-quasi-isometry} if it is an $(L,A)$-quasi-isometric
embedding and for all $x'\in X'$, $d(x',\Phi(X))<A$.

4. $\Phi$ is a {\em $D$-Hausdorff approximation} if it is a $(1,D)$-quasi-isometry.
\end{definition}

We will use the following well-known lemma:

\begin{lemma}
If $G\acts X$ is a discrete, cocompact, isometric
action of a group $G$ on a length space $X$, then there
is a $G$-equivariant quasi-isometry $\Phi:Cayley(G)\ra X$, where $Cayley(G)$ is any Cayley graph of $G$.
\end{lemma}

\subsubsection{Hadamard spaces}
\label{basichadamard}

We refer the reader to \cite{ballmann} for the material recalled here.

\medskip
\no
{\bf Geodesics and the boundary.}
Let $X$ be a locally compact Hadamard space.  If $p,\,q\in X$  then
$\ol{pq}\subset X$ denotes the segment from $p$ to $q$.  If $p,x,y\in X$
and  $p\not\in\{x,y\}$, then $\cangle_p(x,y)$ (respectively $\angle_p(x,y)$)
denotes the comparison angle (respectively angle) of the triangle $\De pxy$
at $p$.
We will use $\geo X$ 
to denote the set of asymptote classes of geodesic rays in $X$, with
the cone topology.  If $p\in X$ and $\xi\in\geo X$, then 
$\ol{p\xi}$ denotes the ray leaving $p$ in the asymptote class
of $\xi$.  $\bar X\defeq X\cup\geo X$ denotes the usual compactification:
a sequence $x_i\in\bar X$ converges if and only if for any basepoint $p\in X$
the sequence of geodesic segments/rays $\ol{px_i}$ converges in the compact open topology.
We denote the Tits angle between $\xi_1,\xi_2\in\geo X$ by
$\tangle(\xi_1,\xi_2)$, and $\tits X$ denotes the underlying set of $\geo X$
equipped with the Tits angle metric (which usually induces a topology different from 
the one defined above).  The metric space $\tits X$ is 
a $CAT(1)$ space with respect to this metric.  
When $\xi_1,\xi_2\in \tits X$ and
$\tangle(\xi_1,\xi_2)<\pi$ 
then there is a segment between
$\xi_1$ and $\xi_2$ in $\tits X$, which we denote by $\ol{\xi_1\xi_2}\subset\tits X$. 
This segment is the limit set in $\bar X$ of any sequence of segments
$\ol{x_1^kx_2^k}$ where $x_i^k$ tends to infinity along a ray asymptotic
to $\xi_i$.   
We will not use the Tits path metric.
We recall that $\geo$ and $\tits$ behave nicely with respect to products: 
$\geo (X_1\times X_2)=\geo X_1\circ\geo X_2$ and 
$\tits (X_1\times X_2)=\tits X_1\circ\tits X_2$ where 
in the first case $\circ$ represents the topological 
join and in the second the $\frac \pi 2$-metric join.  
We will use this in the case where $X_2=\R$.

We will let $N_R(C)$ be the closed metric tubular neighborhood of radius $R$ of a subset $C\subset X$.  
A closed convex subset $C\subset X$ is also a locally compact Hadamard space 
as is $N_R(C)$, since it is also convex.

Standard comparison arguments show the following.

\begin{lemma}
\label{prebdyconvex}
Let  $X$ be a locally compact Hadamard space, and let $C\subset X$
be a closed convex subset.  Then for any $R>0$, $\xi\in \geo N_RC$, 
and  $z\in C$ we have $\ol{z\xi_\infty}\subset C$.  In particular, $\geo C=\geo N_RC$.
\end{lemma}

One consequence is:

\begin{lemma}
\label{bdyconvex}
Let  $X$ be a locally compact Hadamard space, and let $C\subset X$
be a closed convex subset.  If $p\in X$, $\xi_i\in\geo X$, and $\xi_i\ra\xi_\infty$,
$\ol{p\xi_i}\cap C\neq\emptyset$ for all $i$, then either 
$\ol{p\xi_\infty}\cap C\neq\emptyset$ or $\xi_\infty\in\geo C$.  \end{lemma}
\proof
Pick $x_i\in\ol{p\xi_i}\cap C$.  If $\liminf d(x_i,p)<\infty$ then a subsequence 
of $x_i$ converges to $x_\infty\in C\cap\ol{p\xi_\infty}$.  On the other hand, 
if $\liminf d(x_i,p)=\infty$ then for some subsequence $\ol{px_i}\ra\ol{p\xi_\infty}$.  
By the convexity of $N_R(C)$,  $\ol{p\xi_\infty}
\subset N_R(C)$ for $R=d(p,C)$ and hence 
Lemma \ref{prebdyconvex} yields the result. 
\qed

\begin{lemma}
\label{sublinearbending}
Let $\Phi:X\ra X'$ be a quasi-isometric embedding, and assume
there is a point $x\in X$ and a function $\th:\R_+\ra\R_+$
with $\lim_{r\ra\infty}\th(r)=0$ so that 
 for every $y\in X$, $z\in\ol{xy}$, we have 
\begin{equation}
\label{bendingest}
d_{X'}(\Phi(z),\ol{\Phi(x)\Phi(y)})\leq
(1+d_X(z,x))\th(d_X(z,x)).
\end{equation}
Then there is a unique extension $\bar\Phi:\bar X\ra\bar X'$
of $\Phi$ which is continuous at $\geo X$, and $\geo\Phi\defeq
\bar\Phi\restr_{\geo X}$ is a topological embedding.
\end{lemma}
\proof
Let $\Phi$ be an $(L,A)-$quasi isometric embedding, and $\xi\in\geo X$, and $y_k\in X$ be such that 
$\ol{xy_k}$ converges to $\ol{x\xi}$.  By the convergence we can choose
$R_k\ra\infty$ so that for all $k\geq n$ we have
$d(y_k,x)\geq R_n$ and $y_{kn}\defeq \ol{xy_k}\cap S(x,R_n)
\subset N_1(\ol{x\xi})$.  
Note that the point on $\xi$ closest to $y_{kn}$ lies in $\xi([R_n-1,R_n+1])$, 
hence by the triangle inequality $d(y_{kn},y_{ln})\leq 4$ for $k,l\geq n$.  
Using (\ref{bendingest}),
for every $k\geq n$ choose $y_{kn}'\in\ol{\Phi(x)\Phi(y_k)}$
with $d_{X'}(y_{kn}',\Phi(y_{kn}))\leq (1+R_n)\th(R_n)$.
Then for every $k,l\geq n$ we have $d_{X'}(\Phi(y_{kn}),\Phi(y_{ln}))
\leq 4L+A$, and so $d_{X'}(y_{kn}',y_{ln}')\leq 2(1+R_n)\th(R_n)
+4L+A$. This, along with the fact that $d(\Phi(x),y'_{kn})\geq L^{-1}R_n-A-(1+R_n)\th(R_n)$, 
forces $\cangle_{\Phi(x)}(y_{kn}',y_{ln}')$
to zero as $n\ra\infty$.
This in turn forces $\ol{\Phi(x)\Phi(y_k)}$ to converge to a ray $\ol{\Phi(x)\Phi(\xi)}$ 
since for each $R>0$ we have for large enough $k$ that the sequence $\{\ol{\Phi(x)\Phi(y_k)}\cap S(\Phi(x),R)\}$ is Cauchy and hence converges.
This proves that $\Phi$ has a unique extension
$\bar\Phi:\bar X\ra\bar X'$ which is continuous at $\geo X$.
The map $\geo \Phi\defeq \bar\Phi\restr_{\geo X}$ clearly
has the property that for all $\xi\in\geo X$ and all
$y\in x\xi$, 
\begin{equation}
\label{stillsublinear}
d(\Phi(y),\ol{\Phi(x)\geo\Phi(\xi)})\leq
(1+d(x,y))\th(d(x,y)).
\end{equation}
When $\xi_1,\xi_2\in\geo X$ are distinct, the rays
$\ol{x\xi_i}$ diverge linearly, and hence $\Phi(\ol{x\xi_1})$ and $\Phi(\ol{x\xi_2})$ 
diverge linearly since $\Phi$ is a quasi isometric embedding.
Now if $\geo \Phi(\xi_1)=\geo \Phi(\xi_2)$
then (\ref{stillsublinear}) would
imply that $\Phi(\ol{x\xi_1})$ and $\Phi(\ol{x\xi_2})$ would each diverge 
sublinearly from $\ol{\Phi(x)\geo \Phi(\xi_1)}$ and hence diverge sublinearly 
from each other.  Thus we conclude that $\geo \Phi(\xi_1)\neq\geo \Phi(\xi_2)$.
\qed

\medskip
\no
\subsubsection{Groups acting on Hadamard spaces.}
\label{groupsonx} 

Let $X$ be
a Hadamard space.  We denote the displacement function
of an isometry $g:X\ra X$ by $d_g$, and the infimum
of $d_g$ by $\de_g$.  When $g$ is
axial, we let $Minset(g)$ denote  the convex subset 
where $d_g$ attains its minimum.  We recall that
$Minset(g)$ splits as a metric product $C\times\R$
where $C$ is convex and $g$ acts trivially on the $C$ factor and by
translation on the $\R$ factor.

Let $G\acts X$ be a discrete, cocompact, isometric action
of a group $G$ on a Hadamard space $X$.  If $H\subset G$
is a subgroup isomorphic to $\Z^k$, we let 
$Minset(H)\defeq \cap_{h\in H}Minset(h)$.  We recall that
 $Minset(H)=\cap_{h\in S} Minset(h)$ for any 
generating set $S\subset H$, and that $Minset(H)$
splits isometrically as a metric product $C\times\E^k$
so that $H$ acts trivially on the $C$ factor, and as a
translation lattice on the $\E^k$ factor.  The centralizer
$Z(H,G)$ of $H$ in $G$ preserves $\de_h$ for every 
$h\in H$, and hence also $Minset(H)$.  If $S\subset H$
is a finite generating set, then the function 
$\sum_{h\in S}\de_h:X\ra\R$ descends to a proper function
on $X/Z(H,G)$; in particular, $Minset(H)/Z(H,G)$ is compact.

\subsubsection{Gromov hyperbolic groups and spaces}
\label{basichyperbolic}
For background on the material in this section see 
\cite{hypgps}, \cite{delaharpe}, and \cite{codepa}.
Some standard facts that we will use: 
If a Gromov hyperbolic group $G$ acts cocompactly on a 
Hadamard space $X$ then (since $X$ is then quasi-isometric to 
$Cayley(G)$ and Gromov hyperbolicity is a quasi isometry invariant) 
$X$ is Gromov hyperbolic (i.e. $\delta$-hyperbolic for some $\delta$).  
Further $\geo X$ is homeomorphic to $\geo G$, and all infinite order
elements $g\in G$  are axial.
The Tits metric $\tits X$ is the discrete metric with any two distinct points having distance $\pi$.  

In this section we make use of the Morse lemma for quasi-geodesic segments (see \cite{hypgps,codepa}):

\begin{lemma}(Morse Lemma)
\label{qgmorse}
Given $\delta>0$, $L>0$ and $A\geq 0$ there is a constant $C=C(\delta,L,A)$ 
such that if $\ga_1$ and $\ga_2$ are $(L,A)$-quasi-geodesic segments 
with the same endpoints sitting in a
$\de$-hyperbolic  space, then their Hausdorff distance satisfies $\hd(\ga_1,\ga_2)<C$. 
\end{lemma}

Two geodesics $\ga_1$ and $\ga_2$ in a Hadamard space $X$ are 
 {\em parallel} if they stay a bounded distance apart.  The parallel set 
$P(\ga)\subset X$ of a geodesic $\ga$ is the union of all geodesics parallel to $\ga$.  
By the flat strip theorem, $P(\ga)$ is a convex subset of $X$, and is
isometric to  $C_\ga\times\R$ where 
$C_\ga\subset X$ is convex.  A bounded convex set $C$ always contains a unique 
circumcenter: the center of the smallest metric ball containing $C$.

\begin{lemma}
\label{geodesicunion}
Let $X$ be a $\de$-hyperbolic Hadamard space.  Then

1. If $\ga\subset X$ is a geodesic and $P(\ga)\simeq C_\ga\times\R$
is its parallel set, then $Diam(C_\ga)<\de$.  In particular
$P(\ga)$ contains a canonical geodesic $z\times\R\subset C_\ga\times\R$
where $z\in C_\ga$ is the circumcenter of $C_\ga$.

2.  If $\ga_1,\,\ga_2\subset X$ are geodesics, $x_i\in \ga_i$, then 
$\ga_1\cup\ol{x_1x_2}\cup \ga_2$ is $2\de$-quasi-convex\footnote{A
subset $Z\subset X$ is $C$-quasi-convex if for all $x,\,y\in Z$
we  have $\ol{xy}\subset N_C(Z)$.}.

3. 
Suppose $\ga_1,\,\ga_2\subset X$ are geodesics with
$\geo\ga_1\cap\geo\ga_2=\emptyset$, and let $\eta$
a minimal geodesic segment between $\ga_1$ and $\ga_2$.
Then any geodesic segment running from $\ga_1$ to $\ga_2$
will pass within distance $D=D(\ga_1,\ga_2)$ of both
endpoints of $\eta$; when $d(\ga_1,\ga_2)>4\de$ then we may
take $D=2\de$. 
 \end{lemma}
\proof
1 follows from the fact that a $\de$-hyperbolic Euclidean 
strip has width at most $\de$.  2 and 3 follow from repeated
application of the $\de$-thinness property of geodesic triangles.
\qed

\bigskip

In the following lemma is a slight variation on results from
\cite{hypgps}.  It shows that discrete isometric actions on 
Gromov hyperbolic spaces behave like free group
actions  on trees.  

\begin{lemma}
\label{almostfree}
Let $X$ be a $\de$-hyperbolic Hadamard space, and let $\star\in X$.
Suppose $(g_i)_{i\in\Z}$ is a periodic sequence of axial isometries
of $X$ with period $k$ (i.e. $g_{i+k}=g_i$ for all $i$;
and in particular $g_{0}=g_{k}$ 
and $g_{-1}=g_{k-1}$), and let the attracting (respectively repelling)
fixed point of $g_i$ be $\xi_i^+\in\geo X$ (respectively $\xi_i^-\in\geo X$).  If
for every $i$ we have $\xi_i^-\neq\xi_{i+1}^+$, 
then there are constants $L$, $A$, $N$, and $D$ with the following property.

1. If $(m_i)_{i\in\Z}$ is a sequence with $m_i>N$, then the broken
geodesic with vertices
\begin{equation}
\label{brokengeodesic}
\ldots,v_{-2}=
(g_{-1}^{-m_{-1}}g_{-2}^{-m_{-2}})(\star),\,v_{-1}= g_{-1}^{-m_{-1}}(\star),\,v_0=\star,\,
v_1= g_0^{m_0}(\star),\,v_2= (g_0^{m_0}g_1^{m_1})(\star),\,\ldots
\end{equation}
is an $(L,A)$ quasi-geodesic, and 
\begin{equation}
\label{expecteddist}
|d(v_i,v_{i+l})-\sum_{j=i}^{i+l-1}m_j\de_{g_j}|<lD
\end{equation}

2. If $(m_i)_{i\in\Z}$ is a sequence with $m_i>N$ and period $k$,
then $g\defeq g_0^{m_0}\ldots g_{k-1}^{m_{k-1}}$ is an axial isometry with
an axis $\ga$ within Hausdorff distance $D$ of the $g$-invariant broken geodesic
with vertices (\ref{brokengeodesic}), and the minimal displacement of $g$
satisfies
\begin{equation}
\label{dispest}
|\de_g-(m_0\de_{g_0}+\ldots+m_{k-1}\de_{g_{k-1}})|<D.
\end{equation}
Furthermore, as $m_0\ra\infty$ (respectively $m_{k-1}\ra\infty$),
the attracting (respectively repelling) fixed point of $g$
tends to $\xi_0^+$ (respectively $\xi_{k-1}^-$).
\end{lemma}
\proof
Since $\xi_i^-\neq\xi_{i+1}^+$ for all $i\in\Z$, there are constants
$L_1$, $A_1$, and $N_1$, so that when $m_\pm>N_1$ then for any $i$ the 
broken geodesic with vertices $g_i^{-m_-}(\star)$, $\star$, $g_{i+1}^{m_+}\star$
is an $(L_1,A_1)$ quasi-geodesic segment.    
Let  $(m_i)_{i\in\Z}$ be a sequence, and let $\eta:\R\ra X$ be
the broken geodesic  with vertices (\ref{brokengeodesic}).
By the local
characterization of quasi-geodesics given in \cite[Chapitre 3]{codepa},
we get constants $L=L(L_1,A_1,\de)$, $A=A(L_1,A_1,\de)$,
and $N_2=N_2(L_1,A_1,\de,\{g_i\}_{i\in\Z})\geq N_1$  so that $\eta$ is an
$(L,A)$ quasi-geodesic provided $m_i\geq N_2$ for all $i$.   

We now assume that $m_i\geq N_2$ for all $i$.
By the Morse lemma there is a  $D_1=D_1(L,A,\de)$ so that
there is a geodesic at Hausdorff distance at most $D_1$ from $\eta(\R)$,
and any geodesic $\ga\subset X$ with $\geo\ga=\geo\eta$
has Hausdorff distance at most $D_1$ from $\eta$.   Fix such a geodesic
$\ga\subset X$, and for each $i\in \Z$ let $w_i\in\ga$ be the point
in $\ga$ nearest $v_i$.   By the triangle inequality we have
\begin{equation}
\label{disterror1}
|d(v_i,v_{i+1})-d(w_i,w_{i+1})|\leq 2D_1.
\end{equation}
Choose $c_1=c_1(\star,\{g_i\})$ so that the distance from $\star$
to the nearest axis of $g_i$ is less than $c_1$; then for all $i\in \Z$
\begin{equation}
\label{dispest0}
|d(v_i,v_{i+1})-m_i\de_{g_i}|=|d(\star,g_i^{m_i}(\star))-m_i\de_{g_i}|<2c_1.
\end{equation}
Since for each $i$, the broken segment with vertices $v_{i-1},\,v_i,\,v_{i+1}$
is an $(L_1,A_1)$ quasi-geodesic, the Morse Lemma gives
\begin{equation}
\label{triangleqn}
d(v_{i-1},v_{i+1})\geq d(v_{i-1},v_i)+d(v_i,v_{i+1})-2D_1.
\end{equation}
This gives
\begin{equation}
\label{forw's}
d(w_{i-1},w_{i+1})\geq d(w_{i-1},w_i)+d(w_i,w_{i+1})-8D_1.
\end{equation}
Therefore there is an $N=N(\star,\{g_i\})\geq N_2$ so that if $m_i\geq N$
then $w_i$ lies between $w_{i-1}$ and $w_{i+1}$ for all $i$.
So when $m_i\geq N$ we have
$$
|d(v_i,v_{i+l})-\sum_{j=i}^{j+l-1}m_j\de_{g_j}|
\leq 2D_1+|d(w_i,w_{i+l})-\sum_{j=i}^{j+l-1}m_j\de_{g_j}|
$$
\begin{equation}
\label{internalest}
\leq 2D_1+\sum_{j=i}^{j+l-1}|d(w_j,w_{j+1})-m_j\de_{g_j}|
\leq 2D_1(l+1)+2lc_1.
\end{equation}
We now set $D\defeq (2k+4)D_1+2kc_1$, and note that we have proved
1.
When  the sequence $(m_i)$ has period $k$, $m_i\geq N$ for all $i$,
and $g\defeq
g_0^{m_0}\ldots g_k^{m_{k-1}}$, then we may take $\ga$ to be
an axis for $g$.  We have $\de_g=d(w_0,w_k)$, and (\ref{dispest})
follows from (\ref{internalest}).  The last assertion follows immediately from
the fact that as $m_0\ra\infty$ and $m_{k-1}\ra\infty$, the segments
$\ol{\star g_0^{m_0}(\star)}$ and $\ol{\star g_{k-1}^{-m_{k-1}}(\star)}$ 
converge to the rays $\ol{\star\xi_0^+}$ and $\ol{\star\xi_{k-1}^-}$
respectively.\qed

\begin{lemma}
\label{hyplemmas}
Let $G\stackrel{\rho}{\acts} X$ be a discrete, cocompact isometric
action of a hyperbolic group $G$ on a Hadamard space $X$.  There is a constant
$D=D(\de,\rho)$ so that for every $x_1,x_2\in X$ there is a $g\in G$
and an axis $\ga$ for $g$ with $d(x_i,\ga)<D$ for $i=1,2$, and
$d(g(x_1),x_2)<D$.
\end{lemma}
\proof If $G$ is elementary, then either $G$ is finite 
(in which case the result holds trivially) or there is a hyperbolic
element $g\in G$ with an axis $\ga$ so that $X=N_R(\ga)$ for 
some  $R$; this implies 2 in this case.  
So we may assume that $G$ is nonelementary,
and hence $G$ does not fix any $\xi\in\geo X$.  

Pick $\star\in X$ and a finite generating set $\Si\subset  G$.  
Fix $\si_0\in \Si$, let $\Si'=\{\si_0\}\cup \{\si\si_0|\si\in \Si\}$, 
and let $C(\Si')=\min\{d(\star,\si'(\star))|\si'\in \Si\}$.  

We note that by the cocompactness of the action it is sufficient to 
prove the theorem when $x_1$ is $\star$; for then (with a larger $D$) 
if $g_1(x_i)$ is near $\star$ (within the diameter of the 
fundamental domain)  and $g$ is the solution for $\star$ and 
$g_1(x_2)$ then $g_1^{-1}gg_1$ works for $x_1$ and $x_2$ 
(since $g_1^{-1}(\ga_g)$ is an axis for $g_1^{-1}gg_1$).

\no
{\em Claim.  There are constants $L_1$, $A_1$ such that for all
$g\in G$, there is a $\si'\in\Si'$ so that the broken geodesic with 
vertices $(g\si')^{-1})(\star),\,\star,\,(g\si')(\star)$ is an $(L_1,A_1)$-quasi-geodesic.}

\no
{\em Proof of claim.}  If not, there is a sequence $g_k\in G$ with $d(g_k(\star),\star)\ra\infty$,
so that for every $\si'\in\Si'$ the broken geodesic with vertices
$(g_k\si')^{-1})(\star),\,\star,\,(g_k\si')(\star)$ is not a $(k,k)$-quasi-geodesic.
This clearly implies that after passing to a subsequence (which works for all $\si'$), the segments 
$\ol{\star[(g_k\si')(\star)]}$, $\ol{\star[(g_k\si')^{-1}(\star)]}$ converge to some
ray $\ol{\star\xi(\si')}$.  But since $d(g_k(\star),(g_k\si')(\star))\leq C(\Si') $ 
we see $(g_k)(\star)$ converges to $\xi(\si')$, so $\xi(\si')=\xi$  is independent 
of $\si'$. The fact that $(g_k\si')^{-1}(\star)$ converges to $\xi$ tells us 
(by applying $\si'$ to the sequence) that $g_k^{-1}(\star)$ converges 
to $\si'(\xi)=\xi'$ which is again independent of the choice of $\si'$.  
Thus $\si'(\xi)=\xi'$ for all $\si'\in \Si$, and hence for every $\si\in\Si$, 
$\si(\xi')=(\si\si_0)\si_0^{-1}(\xi')=\xi'$, which is a contradiction.
\qed

\no
{\em Proof of Lemma \ref{hyplemmas} continued.}  By the claim and an
application of \cite[Chapitre 3]{codepa} as in the proof of Lemma \ref{almostfree},
we see that there are constants $L$, $A$, and $D_1$ so that if $g\in G$
and $d(g(\star),\star)>D_1$, then there is a $\si'\in\Si'$ so that the 
broken geodesic with $i^{th}$ vertex $(g\si')^i(\star)$ is an $(L,A)$-quasi-geodesic.
By the Morse Lemma and Lemma \ref{geodesicunion}, $(g\si')$ has
an axis at distance $<D_2(L,A,\de)$ from $\star$ and $(g\si')(\star)$.
Since $d((g\si')(\star),g(\star))=d(\si'(\star),\star)\leq C(\Si')$, the lemma
clearly follows.
\qed

\begin{lemma}
\label{preservinggeodesics}
Consider two discrete, cocompact, isometric actions $G\stackrel{\rho}{\acts} X$, $G\stackrel{\rho'}{\acts} X'$ 
where $X$ and $X'$ are $\de$-hyperbolic Hadamard spaces.  Assume that 
the minimum displacement  of any $g\in G$ in $X$ is the same as the minimal
displacement in $X'$.  Then any $G$-equivariant $(L,A)$-quasi-isometry
$\Phi:X\ra X'$ maps unit speed geodesics $\ga$ to within $D=D(L,A,\delta,\rho,\rho')$ of
a unit speed geodesic $\ga'$, that is $d(\Phi(\ga(t)),\ga'(t))\leq D$.
\end{lemma}

\proof
By the Morse lemma on quasi-geodesics, it suffices to show that
$\Phi$ is a $D_1=D_1(L,A,\delta,\rho,\rho')$ Hausdorff approximation.  
Let $D_5=\max\{D_4(\delta,\rho),D_4(\delta,\rho')\}$ where the 
$D_4$'s come from  Lemma \ref{hyplemmas}.  For 
$x_1,x_2\in X$ (resp. $\Phi(x_1),\Phi(x_2)\in X'$) let $g\in G$ 
(resp $g'\in G$) be the elements guaranteed by  
Lemma \ref{hyplemmas}.  Since $x_1$ is $D_5$ close to 
an axis of $g$ we know that $2D_5+\delta_g\geq d(x_1,g(x_1))\geq \delta_g$ 
and hence $3D_5+\delta_g\geq d(x_1,x_2)\geq \delta_g-D_5$.  
Now $$d(\Phi(x_1),\Phi(x_2))\geq d(\Phi(x_1),\Phi(g(x_1)))-
d(\Phi(g(x_1)),\Phi(x_2))\geq $$$$\geq \delta_g-LD_5-A\geq d(x_1,x_2)-(L+3)D_5-A$$
and similarly
$$d(x_1,x_2)\geq d(x_1,g'(x_1))-d(g'(x_1),x_2)\geq $$
$$\geq \delta_{g'}-LD_5-A\geq d(\Phi(x_1),\Phi(x_2))-(L+3)D_5-A.$$
Hence we can take $D_1=(L+3)D_5+A$.
\qed

\begin{lemma}
\label{discreteaxis}
Let $G\stackrel{\rho}{\acts} X$ be a discrete, cocompact 
isometric action on a $\de$-Hyperbolic Hadamard space $X$.

1. If $\{\ga_g|g\in S\subset G\}$ is a collection of distinct 
axis of distinct elements $g\in G$ such that $\{\de_g\}$ 
is bounded then $\{\ga_g\}$ forms a discrete set of geodesics.

2. If two axial elements $g_1,g_2\in G$ have a common fixed point 
in $\geo X$ then they have a common axis (i.e. they have both fixed 
points in common).
\end{lemma}
\proof
Assume some sequence $\ga_{g_i}$ converges to a geodesic $\ga$ and 
let $\star \in \ga$ then the boundedness of $\{\de_{g_i}\}$ says 
that there is a $C$ such that $d(\star,g_i(\star))<C$ but this 
cannot be true for infinitely many distinct $g_i$.  

To see the second statement we can assume that the attracting and repelling fixed points satisfy  $\xi^+_1=\xi^+_2$ and $\xi^-_1\not=\xi^-_2$ (the other cases are similar).  In this case $\{g_1^{-k}g_2g_1^{k}\}$ are distinct since they have distinct fixed point sets $\{\xi_1^+,g_1^{-k}(\xi_2^-)\}$ in $\geo X$ while the axes converge (after taking a subsequence) to an axis of $g_k$. But again there is a $C$ such that $d(\star,g_1^{-k}g_2g_1^{k}(\star))<C$ giving the desired contradiction.

\qed

\subsubsection{Graphs of groups and their Bass-Serre trees}
\label{basicbassserre}

For the remainder of the paper, all group actions on 
simplicial trees will be assumed to be simplicial actions
which do not invert edges, and geodesic segments/rays
in simplicial trees will be unions of edges. 

References for the material in this section are \cite{serre,scottwall,dundicks}. 

\begin{definition}
A {\em graph of groups} is a connected graph ${\cal G}$ together with
a group $G_\si$ labeling each $\si\in Vertex({\cal G})\cup Edge({\cal G})$,
and a monomorphism $G_e\ra G_v$ for each pair $(e,v)$
consisting of an oriented edge $e$
entering a vertex $v$.   An isomorphism of two graphs
of groups is an isomorphism of labeled graphs 
which is compatible with edge monomorphisms.
\end{definition}
Let $G\stackrel{\rho}{\acts} T$ be an action of a group $G$ on a simplicial
tree $T$.  We can define an associated graph of groups 
 ${\cal G}$ as follows.   We let the graph underlying ${\cal G}$
be $G/T$.  For each $\si\in Vertex({\cal G})\cup Edge({\cal G})$
we may label $\si$ with the stabilizer of a lift $\hat\si\subset T$
of $\si$.  An for each pair $(e,v)$, where $e\subset T/G$ is an oriented edge 
with terminus $v\in T/G$,  we can define an edge monomorphism
$G_e\ra G_v$ by composing the inclusion $G_e\defeq G_{\hat e}\ra G_{g{\hat v}}$
($g{\hat v}\in T$ is the terminus of $\hat e$) with the isomorphism 
$G_{g{\hat v}}\ra G_{\hat v}=G_v$ induced by conjugation by $g^{-1}$.
We refer to this as the {\em graph of groups associated with the action
$G\acts T$}.  
\begin{lemma}
If ${\cal G}$ is a graph of groups, then there is a group $G$, a simplicial
tree $T$, and an action $G\stackrel{\rho}{\acts}T$ so that:

1. If $\bar {\cal G}$ is the graph of groups associated with the action
$G\stackrel{\rho}{\acts}T$, then $\bar{\cal G}\simeq{\cal G}$.

2. If $G'\stackrel{\rho'}{\acts}T'$ is another action on a simplicial tree
satisfying 1, then there is an isomorphism $G'\simeq G$ so that
the actions $\rho$ and $\rho'$ become simplicially isomorphic.
\end{lemma}
\no
The (isomorphism class of the) group $G$ is the {\em fundamental
group of ${\cal G}$}, and the tree $T$ (or really the action 
$G\stackrel{\rho}{\acts}T$) is called the Bass-Serre tree of ${\cal G}$.
We note that if $v$ is a vertex of $T$ and $G_v$ is its stabilizer,
then the $G_v$-orbits of $Link(v)$  correspond bijectively
to the elements of $Link(\bar v)$ where $\bar v\in {\cal G}$ is the 
corresponding vertex of $G/T\simeq {\cal  G}$, and the stabilizer of
$\xi\in Link(v)$ is just $G_e$ where $e$ is the edge associated with
$\xi$.

Let ${\cal G}$ be a graph of groups and let $\bar e$ be an edge of ${\cal G}$
with endpoints $\bar v_1$ and $\bar v_2$. 
We let ${\cal G}'$ be the graph of groups determined by $\bar e$.
The fundamental
group of ${\cal G}'$ is
a free product with amalgamation if $\bar e$ is embedded in ${\cal G}$
and an $HNN$ extension if $\bar e$ is a loop.   Choose a lift $e=\ol{v_1v_2}\subset T$
of $\bar e$ to the Bass-Serre tree $T$.  We may identify $G_{\bar v_i}$
with $G_{v_i}$ and $G_{\bar e}$ with $G_e$ in a fashion compatible
with the edge inclusions $G_{\bar e}\ra G_{\bar v_i}$, $G_e\ra G_{v_i}$.
When $\bar e$ is a loop we may choose $t\in G=\pi_1({\cal G})$ so that
$t(v_1)=v_2$ and the composition $G_e\ra G_{v_1}\stackrel{t(\cdot)t^{-1}}{\lra}
G_{v_2}$ agrees with the edge monomorphism $G_e\ra G_{v_2}$.
Set $G'\defeq \< G_{v_1},G_{v_2}\>$ when $\bar e$ is embedded and
set $G'\defeq \< G_{v_1},t\>$ when $\bar e$ is a loop.  Then the orbit
$T'\defeq G'(e)\subset T$ is a $G'$-invariant subtree of $T$, and the action
$G'\acts T'$ is the Bass-Serre action for ${\cal G}'$.     When $\bar e$ is 
embedded we choose subsets  $\Si_i\subset G_{v_i}$ which intersect
each right coset of $G_e$ exactly once;  then any $g\in G'$ can be written 
uniquely in the form 
\begin{equation}
\label{amalgform}
s_1\ldots s_kr
\end{equation}
where $r\in G_e$, $s_i\not\in G_e$,  and the $s_i$'s belong alternately to $\Si_1$
and $\Si_2$.  
The combinatorial distance  from the edge $g(e)$ to $e$
is $k$ and $e_i=s_1\ldots s_i(e)$ is the sequence of edges along the path from $e$ to $g(e)$.   When $\bar e$ is a loop, we choose a cross-section 
$\Si_1\subset G_{v_1}$ (respectively $\Si_{-1}$) of the right cosets of $G_e$
(respectively $t^{-1}G_e t$). 
Then any $g\in G'$ can be written uniquely in the
form 
\begin{equation}
\label{hnnform}
s_1t^{\eps_1}s_2t^{\eps_2}\ldots s_kt^{\eps_k}r
\end{equation}
where for $i=1,\ldots,k$, $\eps_i=\pm 1$, $s_i\in \Si_{\eps_i}$, $r\in G_{v_1}$, and
if $s_i\in G_e$ then $\eps_{i-1}=-\eps_i$.

\subsection{Graphs of groups and the structure of Hadamard spaces on which
they act}

\subsubsection{Admissible groups and actions}

\label{admissible}

 \begin{definition} 
\label{admissibledef}

A graph of groups ${\cal G}$ is {\em admissible} if

1. ${\cal G}$ is a finite graph with at least one edge.

2. Each vertex group $\bar G_v$ has center $Z(\bar G_v)\simeq \Z$,
$\bar H_v\defeq \bar G_v/Z(\bar G_v)$ is a nonelementary hyperbolic
group, and every edge subgroup $\bar G_e$ is isomorphic to
$\Z^2$.

3. Let $e_1,\,e_2$ be distinct directed edges entering a vertex $v$,
and for $i=1,2$ let $K_i\subset \bar G_v$ be the image of the
edge homomorphism $\bar G_{e_i}\ra \bar G_v$.  Then for
every $g\in \bar G_v$, $gK_1g^{-1}$ is not commensurable 
with $K_2$, and for every $g\in \bar G_v-K_i$, $gK_ig^{-1}$ 
is not commensurable with $K_i$.

4. For every edge group $\bar G_e$, if $\al_i:\bar G_e\ra \bar G_{v_i}$
are the edge monomorphisms, then the subgroup generated by $\al_1^{-1}(Z(\bar G_{v_1}))$ and $
\al_2^{-1}(Z(\bar G_{v_2}))$ has finite index in $G_e\simeq \Z^2$.

\no
A group $G$ is {\em admissible} if it is the fundamental group of an
admissible graph of groups.
\end{definition}
Let $G$ be the fundamental group of an admissible graph of groups
${\cal G}$, and let $G\acts T$ be the action of $G$ on the 
associated  Bass-Serre tree.   We let $V\defeq Vertex(T)$ and $E\defeq Edge(T)$
denote the vertex and edge sets of $T$, and when $\si\in V\cup E$
we let $G_\si\subset G$ denote corresponding stabilizer.  Properties
1-4 of definition \ref{admissibledef} \mbox{imply\hspace{2pt}:}

\begin{lemma}
\label{admissiblebassserre}
1. $T$ is an unbounded tree with infinite valence at each vertex,
and $G$ acts on $T$ with quotient ${\cal G}$.

2. Each vertex group $ G_v$ has center $Z(G_v)\simeq \Z$,
$ H_v\defeq  G_v/Z( G_v)$ is a nonelementary hyperbolic
group, and every edge subgroup $ G_e$ is isomorphic to
$\Z^2$.

3.  If $e_1,\,e_2$ are distinct edges emanating from $v\in V$,
then $G_{e_1}$ is not commensurable with $G_{e_2}$.  In particular,
$Z(G_v)\subset G_{e_i}$ since $g\in Z(G_v)$ implies $G_{e_i}=gG_{e_i}g^{-1}
=G_{ge_i}$ forcing $ge_i=e_i$, i.e. $g\in G_{e_i}$. 

4.  If $e\in E$ has endpoints $v_1,\,v_2\in V$,
then $(Z(G_{v_1})\cup Z(G_{v_2}))\subset G_e$ generates a finite
index subgroup of $G_e$.
\end{lemma}

Most of the time we will work with the action $G\acts T$ and ignore
the graph of groups that produced it.

Examples of admissible groups:

1. (Graph manifolds) Let $M$ be a $3$-dimensional nonpositively curved graph
manifold as in Theorem \ref{graphmfldcase}, and let
$M_i$, $i=1,\ldots k$ be the geometric Seifert components of $M$.
Let ${\cal G}$ be the graph of groups which has one  
vertex labeled with $\pi_1(M_i)$ for each $i$, and an
edge labeled by $\Z^2$ for each pair of totally geodesic boundary
tori in the disjoint union $\cup_i M_i$ which are glued
to form $M$.  The edge monomorphisms come from the two 
different embeddings of a gluing torus into Seifert components.

2. (Torus complexes)  Let $T_0,\,T_1,\,T_2$ be flat two-dimensional tori.  For
$i=1,2$, we choose primitive closed geodesics
$a_i\subset T_0$ and $b_i\subset T_i$ with 
$length(a_i)=length(b_i)$, and we glue $T_i$ to $T_0$ by identifying
$a_i$ with $b_i$ isometrically.  We assume that $a_1$
and $a_2$ lie in distinct free homotopy classes, and intersect
at an angle $\al\in (0,\frac{\pi}{2}]$. Let ${\cal G}$
be the graph of groups associated with the decomposition
$(T_0\cup T_1)\dcup (T_0\cup T_1)\ra \cup T_i$.  Note that $T_0\cup T_1$ is 
homeomorphic to $S^1 \times (S^1\wedge  S^1)$, 
so $\pi_1(T_0\cup T_i)=\Z\times F_2$ where $F_2$ is the free group on two generators.

\bigskip

\begin{lemma}
\label{basicgroup}
1.  If $e_1,\,e_2\in E$ are distinct edges incident to $v\in V$,
then $Z(G_v)\simeq\Z$ is a finite index subgroup of 
$G_{e_1}\cap G_{e_2}$.  In particular $G_{e_1}\cap G_{e_2}\simeq \Z$.

2.  If $v_1,v_2\in V$ are the endpoints of an edge $e\in E$,
then $Z(G_{v_1})\cap Z(G_{v_2})=\{id\}$.
\end{lemma}
\proof
Note that $G_{e_1}\cap G_{e_2}$ has infinite index in each
$G_{e_i}$, for otherwise the $G_{e_i}\simeq\Z^2$ would be commensurable,
contradicting 3 of Lemma \ref{admissiblebassserre}.  Thus $G_{e_1}\cap G_{e_2}\simeq \Z$.  Also by 3 of Lemma 
 \ref{admissiblebassserre} we  have $Z(G_v)\subset
G_{e_1}\cap G_{e_2}$, so both are rank $1$ free abelian groups and 1 follows.

2 follows immediately from 4 of Lemma \ref{admissiblebassserre}, since
the $Z(G_{v_i})$ are rank $1$ subgroups of $G_e\simeq \Z^2$ which
generate a finite index subgroup of $G_e$. 
\qed

\begin{lemma}
\label{vertexintersections}
If $v_1,\,v_2\in V$, then

1.  If $d(v_1,v_2)>2$ then $G_{v_1}\cap G_{v_2}=\{ id\}$.

2.  If $d(v_1,v_2)=2$ and $v\in V$ is the vertex between them,
then $Z(G_v)$ is a finite index subgroup of $G_{v_1}\cap G_{v_2}$.

3. If $d(v_1,v_2)=1$, then $G_{v_1}\cap G_{v_2}=G_e$ where
$e=\ol{v_1v_2}$.
\end{lemma}
\proof

We will prove the assertions in reverse order.  Part 3 is immediate since
the action $G\acts T$ does not invert edges, see section \ref{basicbassserre}. 
To prove 2 we let $e_i$ be the edge between $v_i$ and $v$.  Clearly 
$G_{v_1}\cap G_{v_2}=G_{e_1}\cap G_{e_2}$, which by Lemma \ref{basicgroup}
contains  $Z(G_v)$ as a subgroup of finite index.  To prove 1, let $e_1,e,e_2\in E$
be three consecutive edges of the segment $\ol{v_1v_2}$, and let 
$w_1,\,w_2\in V$ be the endpoints of $e$.  Then by 1 of Lemma \ref{basicgroup}, 
$Z(G_{w_i})$ ($\simeq \Z$) has finite index in $G_{e_i}\cap G_e$ ($\simeq \Z$), so we have 
$(G_{e_1}\cap G_e)\cap (G_e\cap G_{e_2})=\emptyset$ since otherwise $Z(G_{w_1})\cup Z(G_{w_2})$ would generate a  cyclic group in $G_e$ ($\simeq \Z^2$) contradicting
4 of Lemma \ref{admissiblebassserre}.  Since
$$G_{v_1}\cap G_{v_2}\subset (G_{e_1}\cap G_e)\cap (G_e\cap G_{e_2})$$
1 follows.
\qed

\begin{lemma}
\label{actionproperties}
For $v\in V$, the fixed point set of $Z(G_v)$ is the closed star $\ol{Star(v)}$.  
For $e\in E$, the fixed point set of $G_e$ is $e$.
In particular, for any $\sigma\in V\cup E$,  $G_\sigma$ leaves no point in 
$\geo T$ fixed.  Further, if $\ol{Star(v)}$ is invariant under $Z(G_\sigma)$ 
then $\sigma \in \ol{Star(v)}$.

\end{lemma}
\proof
Fix a vertex $v$.  Since $Z(G_v)\subset G_e$ for all $e$ in the star 
it is clear that the closed star is in the fixed point set of $Z(G_v)$.  
On the other hand  1 of Lemma \ref{vertexintersections} says we 
need only consider vertices $v_1$ such that $d(v,v_1)=2$. 
Let $w$ be the vertex between $v$ and $v_1$. Now 2 of 
Lemma \ref{vertexintersections} says that $Z(G_w)$ is a 
finite index subgroup of $G_{v}\cap G_{v_1}$, while 2 of 
Lemma \ref{basicgroup} says that $Z(G_v)\cap Z(G_w)=\emptyset$.  
Thus $Z(G_v)\cap G_{v_1}=\emptyset$ and the first statement follows.

For an edge $e=\ol{v_1v_2}$ the first part (since $G_{{v_i}}\subset G_e$) says that the fixed point set of $G_e$ is contained in $\ol{Star(v_1)}\cap\ol{Star(v_2)}=e$.  So the second statement follows.  

The last two statements follow from the first two.
\qed

\begin{lemma}
\label{properdisp}
If $v\in V$ then the centralizer of $Z(G_v)$ in
$G$ is just $G_v$.  If $e\in E$ the centralizer of $G_e$
in $G$ is  $G_e$.
\end{lemma}
\proof
By Lemma \ref{actionproperties}, the fixed point set of $Z(G_v)$ in $T$ is just
the closed star of $v$ in $T$.  Hence any $g\in G$
which commutes with $Z(G_v)$ must take the star of $v$ to itself and
hence fix $v$.

If $e\in E$ and $e=\ol{v_1v_2}$, then a finite index subgroup $G_e$
is generated by $Z(G_{v_1})\cup Z(G_{v_2})$.
So the centralizer of $G_e$ in $G$ is a subgroup of the intersection
of the centralizers of $Z(G_{v_1})$ and $Z(G_{v_2})$, i.e. $G_{v_1}\cap
G_{v_2}$ which is
$G_e$ itself.
\qed

\begin{lemma}(Uniqueness of decomposition)
\label{uniquedecomp}
Let $G\acts T$ and $G'\acts T'$ be the Bass-Serre actions associated with two
admissible graphs of groups, and suppose $G\stackrel{\phi}{\ra} G'$
is an isomorphism.  Then after identifying $G$ with $G'$ via $\phi$,
the trees $T$ and $T'$ become $G$-equivariantly isomorphic.
\end{lemma}

\proof
We will use primes to denote the vertex and edge set of $T'$.
Pick $v\in V$.   

\medskip
\no
{\em Claim:  $G_v$ fixes a unique vertex in $T'$.}  Let $g\in Z(G_v)$ be  
a generator.   If $Fix(g,T')$ is empty, then $g$ translates a unique
geodesic $\ga\subset T'$, and since $g\in Z(G_v)$ the whole vertex
group $G_v$ must preserve $\ga$, and act on it by translations.
The signed translation distance yields a homomorphism $G_v\ra\Z$
with nontrivial kernel.   But then $Ker(G_v\ra\Z)$ fixes $\ga$ pointwise,
which contradicts 1 of Lemma \ref{vertexintersections}.   Consequently
$Fix(g,T')$ is nonempty, and by 1 of Lemma \ref{vertexintersections}
this is a subcomplex of $T'$ with diameter at most $2$.   So $G_v$
must fix the center of $Fix(g,T')$.   It can fix nothing more, since
no edge stabilizer can contain the nonabelian $G_v$.  Thus we have proved the
claim.

Now consider the $G$-equivariant map $f:V\ra V'$ which assigns to each $v\in V$ the 
unique vertex in $T'$ fixed by $G_v$; and define a map $f':V'\ra V$
by reversing the roles of $T$ and $T'$.   For all $v\in V$, 
$G_v$ fixes $f'\circ f(v)$, so we must have $f'\circ f(v)=v$; similar reasoning
applies to $f\circ f'$, and we see that $f$ and $f'$ are inverses. 
The maps $f$ and $f'$ are adjacency preserving since two vertices
are adjacent iff their stabilizers intersect in a subgroup isomorphic to 
$\Z^2$.  It is now straightforward to see that $f$ defines a $G$-equivariant isomorphism $T\ra T'$.
\qed

Lemma \ref{uniquedecomp} justifies use of the phrase ``$G\stackrel{\rho}{\acts}T$ is
the Bass-Serre tree of the admissible group $G$.''

\subsubsection{Vertex spaces, edge spaces, and geometric data for admissible actions}
\label{sectionvertexandedge}

\begin{definition}
\label{admissibleactiondef}
We say that $G\stackrel{\rho}{\acts} X$ is an {\em admissible action} if $G$ is an admissible group,
$X$ is a Hadamard space, and the action is discrete, cocompact, and isometric.
\end{definition}
For the remainder of this section $G\acts X$ will be a fixed admissible 
action.  In particular, all constants depend on $G\acts X$ (i.e. the group $G$, the Riemannian manifold $X$, and the action) in addition to other explicitly mentioned quantities.
By Lemma \ref{uniquedecomp} there is an essentially unique 
admissible graph of groups associated with $G$, and we will let $G\acts T$
be the corresponding Bass-Serre tree.    

We refer the reader to section \ref{groupsonx}
for properties of Minsets that we use here.  For each $v\in V$ we let $Y_v\defeq
Minset(Z(G_v))\defeq\cap_{g\in Z(G_v)}Minset(g)$ (this will be the $Minset$ of a generator), and
for every $e\in E$ we let $Y_e\defeq Minset(G_e)\defeq\cap_{g\in G_e}
Minset(g)$.

$Minset(\alpha)$ is a convex subset of $X$, invariant under the
centralizer of $\alpha$, which is a metric product of $\R$ with a
Hadamard space.  If $\alpha$ belongs to a group of isometries that acts
cocompactly on $X$ then the centralizer of $\alpha$ acts cocompactly on
$Minset(\alpha)$ (see section \ref{groupsonx}).  Thus $Y_v$ is the product of $\R$
with a Hadamard space $\bar{Y_v}$.  $Z(G_v)$ acts by translation on the
$\R$ factor and the induced action of $H_v$ on $\bar{Y_v}$ is 
discrete and cocompact.
$Y_e$ is the product of $\R^2$ with a compact Hadamard space $\bar{Y_e}$,
and $G_e = \Z^2$ acts by translations on the $\R^2$ factor (section \ref{groupsonx}).

Note that the assignments $v\mapsto Y_v$ and $e\mapsto Y_e$
are $G$-equivariant with respect to the natural $G$ actions.
  The minimal displacement of a generator of
$Z(G_v)$ is the same as that of a generator of
$Z(G_{g(v)})=gZ(G_v)g^{-1}$. By the finiteness of $\cal{G}$ there is a
number $C$ such that for all $v\in V$ the minimal displacement of a
generator of $Z(G_v)$ is less than C.

\begin {definition}
\label{geometricdatadef}
Let $G\stackrel{\rho}{\acts}X$ be an admissible action, and let
$T$ be the Bass-Serre tree for $G$.  For each $v\in V$ we
choose a generator $\zeta_v\in Z(G_v)$ in a $G$-equivariant way.
We have an isometric splitting $Y_v\simeq \bar Y_v\times\R$,
which is preserved by $G_v$.  The choice of generator
$\zeta_v$ defines an orientation of the $\R$ factor of $Y_v$.
We have a map $MLS_v:G_v\ra \R_+$ which assigns to each
$g\in G_v$ the minimum displacement of the induced isometry
$\bar Y_v\ra\bar Y_v$.  $MLS_v$ descends to $G_v/Z(G_v)\simeq
H_v$ since $Z(G_v)$ acts trivially on $\bar Y_v$.   We define a 
homomorphism $\tau_v:G_v\ra\R$
by sending $g\in G_v$ to the signed distance that
$g$ translates the $\R$ factor of $Y_v\simeq \bar Y_v\times\R$.  
The collections of functions $MLS_v$ and $\tau_v$ constitute 
the {\em geometric data} of the action.   Both $MLS_v$ and $\tau_v$
descend to functions of the vertex groups of the graph of groups
${\cal G}$ defining $G$; we will sometimes find it more convenient to
think of the geometric data in this way.
\end{definition}
We remark that it follows from the discreteness of the action
$H_v\acts \bar Y_v$ that  $g\in G_v$ then $MLS_v(g)=0$ iff $g$ projects
to an element of finite order in $H_v$.
\begin{lemma}
\label{yslocallyfinite}
The collections $\{ Y_v\}_{v\in V}$ and $\{Y_e\}_{e\in E}$
are locally finite.  More precisely, for every $R$ there
is an $N$ so that if $x\in X$ then there are at most
$N$ elements $\si\in V\cup E$ so that $Y_\si\cap B(x,R)\neq\emptyset$.
\end{lemma}
\proof
Suppose $v\in V$ and $p\in Y_v$.  Then $p$ has  displacement $<C$ under
the generators of $Z(G_v)$.   Therefore if $p\in B(x,R)$,
then $x$ has displacement $<2R+C$ under the generators of $Z(G_v)$.
But there are only finitely many $g\in G$ with $d(g\cdot x,x)<2R+C$;
since $Z(G_{v_1})\cap Z(G_{v_2})=\{e\}$ when $v_1\neq v_2$,
the local finiteness of $\{ Y_v\}_{v\in V}$ follows.  Similar
reasoning proves the local finiteness
of $\{ Y_e\}_{e\in E}$.  The fact that $N$ can be chosen independent of
$x$ follows from the cocompactness of the $G$ action.
\qed

\bigskip

The lemma implies that for any $D$, the collection of
$D$-tubular neighborhoods of the $Y_\cdot$'s is locally
finite.  We also have the following consequences:

\begin{lemma}
\label{finitepairs}
For every $D$ there are only finitely many pairs $(\si_1,\,\si_2)\in (V\cup E)\times (V\cup E)$  
-- modulo the diagonal action of $G$ --
with $N_D(Y_{\si_1})\cap N_D(Y_{\si_2})\neq
\emptyset$.  
\end{lemma}
\proof
By the finiteness of $[V\cup E]/G$ we need only show that for fixed $\sigma$ there 
are only finitely many $\sigma_2$ modulo $G_\sigma$ such that $N_D(Y_{\si})\cap N_D(Y_{\si_2})\neq
\emptyset$.
This follows from Lemma \ref{yslocallyfinite},
since $G_\si$ acts cocompactly on $N_D(Y_\si)$ and hence for some 
$g\in G_\si$ $N_D(Y_{\si})\cap (N_D(Y_{g(\si_2)})$ intersects a fixed ball.
\qed

\begin{lemma}
\label{actscocompactly}
For $\si_1,\,\si_2\in V\cup E$, $G_{\si_1}\cap G_{\si_2}$ acts cocompactly on the
intersection $N_D(Y_{\si_1})\cap N_D(Y_{\si_2})$.  Thus, in particular, 
the diameter of $N_D(Y_{\si_1})\cap
N_D(Y_{\si_2})/[G_{\si_1}\cap G_{\si_2}]$ is uniformly bounded by a function
of $D$.
\end{lemma}
\proof
This follows from the local finiteness of the family
$\{Y_\si\}_{\si\in V\cup E}$ and the discreteness
of the cocompact action $G\acts X$.  Pick $D>0$
and $\si_1,\,\si_2\in V\cup E$.  If 
$x_k\in N_D(Y_{\si_1})\cap N_D(Y_{\si_2})$, 
we may choose a sequence $g_k\in G_{\si_1}$ such that
$g_k(x_k)\ra x_\infty$ for some $x_\infty\in N_D(Y_{\si_1})$.
Then $g_k(\si_2)$ lies in a finite subset of $V\cup E$ (since
$g_k(N_D(Y_{\si_2}))$ intersects some ball $B(x_\infty,R)$ for all $k$)
so after passing to a subsequence if necessary we may assume
that $g_k\si_2$ is constant.  Then $g_1^{-1}g_k\in G_{\si_1}\cap G_{\si_2}$
and $(g_1^{-1}g_k)(x_k)\ra g_1^{-1}(x_\infty)$.  Thus $G_{v_1}\cap G_{v_2}$ 
acts cocompactly.  The second statement now follows from Lemma \ref{finitepairs}.

\qed

\begin{lemma}
\label{separatingys}
For every $D$ there is a $D'$ (depending only on D) such that if $\si\in
V\cup E$ separates $\si_1\in V\cup E$ from $\si_2\in V\cup E$,
then $N_D(Y_{\si_1})\cap N_D(Y_{\si_2})\subset N_{D'}(Y_\si)$.
In particular, if $T_1$ and $T_2$ are the closures of distinct connected components
of $T-\si$  then
$$[\cup_{\hat\si\subset T_1}N_D(Y_{\hat\si})]\cap [\cup_{\hat\si\subset T_2}N_D(Y_{\hat\si})]\subset
N_{D'}(Y_\si).$$
\end{lemma}
\proof
Pick $D>0$.  Suppose $(\si_1,\si_2,\si)$ is a triple with 
$\si_i,\,\si\in V\cup E$, $\si$ separates $\si_1$
from $\si_2$ in $T$, and $N_D(Y_{\si_1})\cap N_D(Y_{\si_2})\neq\emptyset$.
Then $G_{\si_1}\cap G_{\si_2}\subset G_\si$ and $G_{\si_1}\cap
G_{\si_2}$ acts cocompactly on $N_D(Y_{\si_1})\cap N_D(Y_{\si_2})$ by 
Lemma \ref{actscocompactly};  hence $d(Y_\si,\cdot)$ is bounded
on $N_D(Y_{\si_1})\cap N_D(Y_{\si_2})$.   By Lemma \ref{finitepairs}
there are only
finitely many such triples $(\si_1,\si_2,\si)$ modulo $G$, so the lemma follows.
\qed

\begin{definition}
\label{vertexspaces}
Since $G$ acts cocompactly on $X$
 we can now fix a $D$ so that $\cup_{v\in V}N_D(Y_v)
=\cup_{e\in E}N_D(Y_e)=X$.  We
define $X_v\defeq N_D(Y_v)$ for all $v\in V$.   Let  $D'$ denote the
constant in the previous lemma, $D''=\max (D,D')$,
and set $X_e\defeq N_{D''}(Y_e)$ for all
$e\in E$.  We will refer to the $X_v$'s and $X_e$'s as {\em vertex spaces}
and {\em edge spaces} respectively.  
\end{definition}

We note that Lemma \ref{prebdyconvex} 
implies that for any $\sigma\in V\cup E$ we have $\geo X_\sigma = \geo Y_\sigma$.
We summarize the properties of vertex
and edge spaces:

\begin{lemma}
\label{propsofvespaces}
There is a constant $C_1$ with the following property.

1. $\cup_{v\in V}X_v=\cup_{e\in E}X_e=X$.

2. If $\hat{e}\in E$ and $T_1$ and $T_2$ are the distinct connected components
of $T-Int(\hat{e})$, then $[\cup_{v\in T_1}X_v]-X_{\hat{e}}$ and 
$[\cup_{v\in T_2}X_v]-X_{\hat{e}}$ are disjoint
closed and open subsets of $X-X_{\hat{e}}$; and
$\cup_{e\in T_1}X_e-N_{C_1}(X_{\hat e})$ and $\cup_{e\in T_2}X_e-N_{C_1}(X_{\hat e})$
are disjoint closed and open subsets of $X-N_{C_1}(X_{\hat e})$.

3. If $\si_1,\,\si_2\in V\cup E$ and $X_{\si_1}\cap X_{\si_2}\neq\emptyset$
then $d_T(\si_1,\si_2)<C_1$.
\end{lemma}
\proof
1 and 2 follow from the definition of vertex/edge spaces and 
Lemma \ref{separatingys}.  By Lemma \ref{finitepairs}
we can choose $C_1$ so that 3 holds.
\qed

\begin{corollary}
\label{titscaptits}
For any $v\in V$, $\tits X_v$ is isometric to the metric suspension of
an uncountable discrete space, and for every $e\in E$, $\tits X_e$
is isometric to a standard circle.  Pick $v_1,\,v_2\in V$.  

1. If $d(v_1,v_2)>2$, then $\tits X_{v_1}\cap\tits X_{v_2}=\emptyset$.

2. If $d(v_1,v_2)=2$ and $v$ is the vertex in between $v_1$ and $v_2$,
then $\tits X_{v_1}\cap\tits X_{v_2}=\tits \ga$ where $\ga\subset
Y_v$ is a geodesic of the form $\{p\}\times\R\subset \bar Y_v\times \R=Y_v$;
i.e. $\tits X_{v_1}\cap\tits X_{v_2}$ is the pair of suspension points of
$\tits X_v$.

3. If $d(v_1,v_2)=1$, then $\tits X_{v_1}\cap\tits X_{v_2}=\tits X_e\simeq S^1$,
where $e\defeq\ol{v_1v_2}$.
\end{corollary}
\proof
Since $Y_v\simeq\bar Y_v\times\R$, we have $\tits X_v=\tits Y_v=\Si(\tits\bar Y_v)$,
and since $\bar Y_v$ admits a discrete cocompact action by the non-elementary
hyperbolic group $H_v\defeq G_v/Z(G_v)$, $\tits Y_v$ is a discrete set with
the cardinality of $\R$.   For all $e\in E$, $Y_e\simeq\bar Y_e\times\R^2$ where $\bar Y_e$
is compact, so $\tits X_e=\tits Y_e\simeq\tits\R^2$, and the latter is the standard
circle.

Pick $v_1,\,v_2\in V$, and choose $R$ large enough that 
$Z\defeq N_R(X_{v_1})\cap N_R(X_{v_2})\neq\emptyset$.  Then
$\tits X_{v_1}\cap \tits X_{v_2}=\tits Z$.  The lemma now follows
from Lemmas \ref{vertexintersections}  and \ref{actscocompactly}.
\qed

\begin{lemma}
\label{findingz's}
There is  a constant $C_2$ with the following property.
Suppose $v,\,v'\in V$,  $e_1,\ldots,e_n\in E$ are the consecutive edges
of the segment $\ol{vv'}\subset T$, $x\in X_v$ and 
$y\in X_{v'}$.  Then for $1\leq i\leq n$ we can find
points $z_i\in\ol{xy}$ such that 

1. $d(z_i,X_{e_i})<C_2$

2. For all $1\leq i\leq j\leq n$
we have $d(z_i,x)\leq d(z_j,x)$.  

3.  For every $p\in X$ we have $\#\{z_i\in B(p,1)\}<C_2$.
\end{lemma}

\proof
Pick $v,\,v'\in V$, $x\in X_v$, and $y\in X_{v'}$.
Suppose $\hat e\in E$ and let $T_1$ and $T_2$ be the two
connected components of $T-Int(\hat e)$.  If
$\ol{xy}\cap N_{C_1}(X_{\hat e})=\emptyset$ (hence in particular $\ol{xy}\cap X_{\hat e}=\emptyset$) then by the first part of 2 of Lemma \ref{propsofvespaces},
$\ol{xy}$ is contained in one of the two disjoint open
sets  $(\cup_{e\subset T_i}X_e)-N_{C_1}(X_{\hat e})$
for $i=1$ or $i=2$ . 
It follows that $\ol{xy}\cap N_{C_1}(X_{e_i})$ is nonempty
for every $1\leq i\leq n$.  Let $w_i\in \ol{xy}$
be the point in $\ol{xy}\cap N_{C_1}(X_{e_i})$ closest to $x$.
Let $z_1=w_1$, and let $z_i$ be the element of
$\{ w_i,\ldots,w_n\}$ closest to $x$.  So we have
either $z_i=w_i$, or $z_i=w_{i'}$ for some
$i'>i$.  In the latter case $\ol{xz_i}\subset \cup_{e\subset T'}X_e$ (hence in particular $z_i=w_{i'}
\in \cup_{e\subset T'}X_e$) where $T'\subset T$
is the component of $T-Int(e_i)$ containing $e_1$,
so by Lemma \ref{separatingys} we have 
$z_i\in N_{D'}(X_{e_i})$ where $D'$ depends only 
on $C_1$.    If $p\in X$
and $z_i\in B(p,1)$, then $X_{e_i}\cap B(p,1+D')\neq
\emptyset$ and thus $Y_{e_i}\cap B(p,1+D'+D'')\neq
\emptyset$ so by Lemma \ref{yslocallyfinite}
we have 
$\#\{z_i\in B(p,1)\}<N$ where $N$ depends only on $D'$.
Setting $C_2\defeq \max\{D',N\}$, the lemma follows.
\qed

\subsubsection{Itineraries}
\label{itinsection}
Our next objective is to associate an itinerary to any
ray $\ol{p\xi}\subset X$ which is not contained
in a finite tubular neighborhood of a single
vertex space; the itinerary of $\ol{p\xi}$ is a 
ray in $T$ which (roughly speaking) records the 
sequence of vertex spaces visited by $\ol{p\xi}$.

Let $\rho:X\ra V\subset T$ be a $G$-equivariant coarse Lipschitz map
from
the Hadamard space $X$ to the vertex set of the tree $T$ with
the property that for every $x\in X$ we  have $x\in X_{\rho(x)}$.
Such a $\rho$ may be constructed as follows.  Let $\Si\subset X$
be a set theoretic cross-section for the free action $G\acts X$;
define $\rho_0:\Si\ra T$ so that $\si\in X_{\rho_0(\si)}$ for every
$\si\in\Si$, and then extend $\rho_0$ to an equivariant
map $X\ra T$.  Let $L$ be such that $N_{2D}(Y_{v_1})\cap
N_{2D}(Y_{v_2})\neq \emptyset$ implies $d(v_1,v_2)<L$, which exists by
Corollary \ref{finitepairs}.  In particular if $d(x,y)< 2D$ then
$d(\rho(x),\rho(y))\leq L$.  In general, by dividing $\ol{xy}$ into less than 
$\frac {d(x,y)} {2D}+1$ segments of length less than $2D$ and adding 
the previous estimates we see that $\rho$ will be coarse 
Lipschitz; i.e. $d(\rho (x),\rho (y))\leq \frac L{2D}\ d(x,y)+L$.

\begin{lemma}
If $\ga:[0,\infty)\ra X$ is a geodesic, then
$\rho\circ\ga:[0,\infty)\ra T$ has the bounded backtracking
property\footnote{A
map $c:[0,\infty)\to T$ has the {\em bounded backtracking
property} if for every $r\in (0,\infty)$ there is an $r'\in (0,\infty)$
such
that if $t_1<t_2$, and $d(c(t_1),c(t_2))>r'$, then $d(c(t),c(t_1))>r$
for every $t>t_2$.}.
\end{lemma}
\proof
Let $e\in E$ be an edge in $T$, and let $T_1,\,T_2\subset T$ be the
connected components of $T-Int(e)$.  Suppose $\rho\circ\ga(t_1)\in T_1$,
and $\rho\circ\ga(t_2)\in T_2-N_{C_1}(e)$, where $t_2>t_1$.
By Lemma \ref{separatingys} we have
$$[\cup_{v\in T_1}X_v]\cap [\cup_{v\in T_2}X_v]\subset
X_e.$$
Therefore there is a $t_3\in [t_1,t_2)$ such that
$\ga(t_3)\in X_e$.  Since $d(\rho\circ\ga(t_2),e)\geq {C_1}$,
the choice of ${C_1}$ and lemma \ref{propsofvespaces} implies that $\ga(t_2)\not\in X_e$.
Hence the convexity of $X_e$ gives
$\ga([t_2,\infty))\subset X- X_e$, which forces $\rho\circ
\ga([t_2,\infty))\subset T_2$.
This property clearly implies uniformly bounded backtracking.
\qed

\begin{lemma}
\label{itindichotomy}
If $\ga:[0,\infty)\ra X$ is a geodesic ray, then
one of the following holds:

1. $\rho\circ\ga:[0,\infty)\ra T$ is unbounded, and
$\rho\circ\ga([0,\infty))$
lies in a uniform tubular neighborhood of a unique geodesic ray, $\tau$,
in $T$ starting at $\rho(\ga(0))$.  The geodesic $\ga$ intersects
$X_e$ for all but finitely many edges $e$ of $\tau$.

2.  $\rho\circ\ga:[0,\infty)\ra T$ is bounded, and $\ga$ eventually
lies in $N_{D'}(Y_v)$ (where $D'$ comes from Definition \ref{vertexspaces}) for some $v\in V$.  In this case
there is a subcomplex $T_\ga\subset T$ defined by the property that
for each simplex $\si$ in $T$, $\si\in T_\ga$ if and only if $\ga$ is asymptotic to $X_\si$ .  
The possibilities for $T_\ga$
are: a single vertex $v\in V$, a single edge $e\in E$ along with its vertices,
or the closed star $\ol{Star(v)}$ for some $v\in V$.   

\end{lemma}

\proof
Pick $v\in V$.  By the convexity of $N_{D'}(Y_v)$, either
$\ga$ is eventually contained in $N_{D'}(Y_v)$, or
$\ga$ is eventually contained in $X-N_{D'}(Y_v)$.
In the latter case $\rho\circ\ga$ eventually remains
in a unique component of $T-v$, by Lemma
\ref{separatingys}.

If for every $v\in V$ the ray $\ga$ eventually lies
in $X-N_{D'}(Y_v)$, then clearly $\rho\circ\ga$
is unbounded and hence it must lie within uniform distance
of a ray in $T$ by the bounded backtracking property.
So we may assume that $\ga$ is eventually contained in
$N_{D'}(Y_v)$ for some $v\in V$.  We note that if $e\in T_\ga$ 
then any vertex $v'$ of $e$ must also be in $T_\ga$ since 
$Y_e\subset Y_{v'}$. Also if $v,v'\in V$ and $d(v,v')> 2$ then part 1 of 
Lemma \ref{vertexintersections} along with Lemma \ref{actscocompactly} 
says that for any $K$, $N_KX_v\cap N_KX_{v'}$ is compact so 
$\ga\not\in N_KX_v\cap N_KX_{v'}$ and hence at most one of $v$ and $v'$ 
can be in $T_\ga$.  

If there are vertices $v_1$ and $v_2$ in $T_\ga$ 
with $d(v_1,v_2)=2$ and $v$ is the vertex between them then 
part 1 of Lemma \ref{vertexintersections} along with Lemma 
\ref{actscocompactly} says that $Z(G_v)$ acts cocompactly on 
$N_KX_{v_1}\cap N_KX_{v_2}$ which contains $\ga$ for some $K$, and hence
 there is a $K'$ such that for all $t>0$ there is a $g_t\in Z(G_v)$ 
such that $d(\ga(t),g_t(\ga(0))<K'$ and hence $\ga$ stays a 
distance at most $K'+d(\ga(0),Y_v)$ from a geodesic in the $\R$ 
direction of $Y_v=\bar Y_v\times \R$ (since $Z(G_v)$ translates 
the $\R$ direction).
Thus for every $e$ with $v$ as a vertex we have $\ga$ is 
asymptotic to a geodesic in $Y_e$ and hence $e\in T_\ga$.  
Thus $\ol{Star(v)}\subset T_\ga$. But since vertices in $T_\ga$ are at 
most distance 2 apart we see that $\ol{Star(v)}= T_\ga$.

The only cases left for $T_\ga$ are the two mentioned and the case 
of two vertices a unit distance apart.  But in the final case a similar 
argument shows that if $e$ is the edge between them then $\ga$ 
stays a bounded distance from $Y_e$ and hence $e$ must also be in $T_\ga$.
\qed

\begin{definition}
\label{itindef}
Let $\ga$ be a geodesic ray in
$X$.  If case 1 of Lemma \ref{itindichotomy} applies
then we will say that $\gamma$  has {\em itinerary}
 $\tau$,  and otherwise we say that the itinerary of $\ga$ 
 is the subtree $T_\ga\subset T$ described in case 2
 of the lemma. In either case we denote the itinerary of 
  $\ga$ by $\itin(\ga)$.
\end{definition}
 
One immediate consequence of the proof of Lemma \ref{itindichotomy} is 
\begin{corollary} 
\label{staritin}
If $\itin (\ga)=\ol{Star(v)}$ then $\ga$ is asymptotic to either the positive 
or the negative $\R$ direction in the decomposition $Y_v=\bar Y_v\times \R$.
\end{corollary}

\begin{lemma}
\label{asympitin} 
 If  $\ol{p_1\xi},\,\ol{p_2\xi}\subset X$ are asymptotic
 geodesic rays, then either both $\itin(\ol{p_1\xi})$ and $\itin(\ol{p_2\xi})$ 
are finite subtrees, in which case they agree, or both $\itin(\ol{p_1\xi})$ 
and $\itin(\ol{p_2\xi})$ are rays, in which case    $\geo\itin(\ol{p_1\xi})=\geo\itin(\ol{p_2\xi})$.
   In other words, $\itin(\ol{p_1\xi})$ is a ray in $T$
   if and only if $\itin(\ol{p_2\xi})$ is a ray in $T$ asymptotic to 
   $\itin(\ol{p_1\xi})$. 
 \end{lemma}

 \proof
 $\ol{p_1\xi}$
 lies in a tubular neighborhood of some $Y_v$ if and only if 
 $\ol{p_2\xi}$ does, thus the case where $\itin(\ol{p_1\xi})$ (or $\itin(\ol{p_2\xi})$) 
is finite follows.  Thus $\ol{p_1\xi}$ has itinerary a ray $\tau_1$
 if and only if $\ol{p_2\xi}$ has itinerary $\tau_2$ for some ray
 $\tau_2\subset T$.  But the sets $\rho(\ol{p_1\xi})$
 and $\rho(\ol{p_2\xi})$ are at finite Hausdorff distance from 
 one another since $\rho$ is coarse Lipschitz; hence the
 $\tau_i$ are asymptotic.
 \qed

 By the lemma we have a well-defined $G$-equivariant
 map  from $\geo X$ to the union
 $$\geo T\cup (\mbox{finite subsets of $T$})$$
  which assigns to each
 $\xi\in \geo X$ either 
 $\geo \itin(\ol{p\xi}), \,p\in X$ if $\itin(\ol{p\xi})$
 is a ray or $\itin(\ol{p\xi})$ otherwise; we will also
 denote this map by $\itin$.
 If $\eta\in\geo T$, we use $\geo^\eta X$ to denote
 the corresponding subset: $\geo^\eta X\defeq
 \itin^{-1}(\eta)\subset\geo X$.
 We will say that $\geo^\eta X$ is
{\em trivial} if $\geo^\eta X$ is a point or {\em nontrivial} otherwise; 
(in the latter case we will see that $\geo^\eta X$ is homeomorphic 
to a closed interval and is in fact an interval in the Tits metric.).

In particular $\geo X=(\cup_{v\in V}\geo X_v)\cup (\cup_{\eta\in \geo T}\geo^\eta X)$, 
where $\cup_{v\in V}\geo X_v$ is disjoint from $\cup_{\eta\in \geo T}\geo^\eta X$.

The cone topology and Tits metric on $\geo X_v=\geo Y_v=\geo (\bar Y_v\times \R)$ 
is described in sections \ref{basichadamard} and \ref{basichyperbolic}.  
We see that in the cone topology $\geo X_v$ is just the suspension 
$\Sigma (\geo H_v)$ and is independent of the metric on $X$.  
The Tits metric is just the metric suspension
 of the discrete metric.

\bigskip

We now study the dynamics of the action of $G$
on $\geo X$.

\begin{lemma}
\label{fixedinbdy}
1.  For every $v\in V$, the fixed point set of
$Z(G_v)$ in $\geo X$ is $\geo X_v$; this set is
homeomorphic to the suspension of  $\geo H_v$
where $H_v$ is the nonelementary hyperbolic group $G_v/Z(G_v)$.

2. For every $e\in E$, $Fix(G_e,\geo X)=\geo X_e$ which is homeomorphic to a circle.
\end{lemma}
\proof Let $\xi\in \geo X$ be fixed by $Z(G_v)$ and $p\in Y_v$.  
If  $\itin(\ol{p\xi})$ is a ray then by Lemma \ref{asympitin} $\geo \itin(\ol{p\xi})$ 
is fixed by $Z(G_v)$.  But this can not happen since by Lemma 
\ref{actionproperties} $Z(G_v)$ leaves no point in $\geo T$ fixed.  
So $\itin(\ol{p\xi})$ is a finite subtree which by Lemma \ref{asympitin} 
is invariant under $Z(G_v)$.   Thus by Lemma \ref{itindichotomy} and 
Lemma \ref{actionproperties} $v\in \itin(\ol{p\xi})$.  Thus $\ol{p\xi}\subset Y_v$ 
and hence $\xi\in \geo X_v=\geo Y_v$.  On the other hand geodesics rays 
$\ol{p\xi}$ in $Y_v$ are translated by a fixed amount by elements 
$g\in G_v$, so $g(\ol{p\xi})$ is asymptotic to $\ol{p\xi}$ and so $\xi$ 
is fixed by $g$.  The rest of part 1 follows from sections \ref{basichadamard} and \ref{basichyperbolic} as above. 

Let $\xi$ be fixed by $G_e$ and $p\in Y_e$.  Again since 
$G_e$ leaves no point in $\geo T$ fixed, by Lemma \ref{actionproperties} 
$\itin(\ol{p\xi})$ is a finite subtree that is invariant under $G_e$ and 
hence by Lemma \ref{itindichotomy} and Lemma \ref{actionproperties} 
must contain $e$.  Hence we see $\ol{p\xi}\subset Y_e$ and 
$\xi \in \geo Y_e=\geo X_e$.  Again, since $G_e$ acts by translations 
on $Y_e$, we see that if $\xi\in \geo X_e=\geo Y_e$ then it is left fixed by 
$G_e$. Since $Y_e=\R^2\times \ol{Y_e}$ where $\ol{Y_e}$ is compact,
$\geo Y_e=\geo \R^2$ is homeomorphic to a circle.
\qed

\bigskip

\subsection{Templates and the behavior of their geodesics}

\subsubsection{Templates}
\label{subsectemplates}

In this section we study ``Templates''.  These are piecewise Euclidean
Hadamard spaces (which can be embedded in $\R^3$) which approximate
certain subspaces 
of the spaces we are studying, and carry much of the
information about the spaces at infinity.

A {\em template} is a Hadamard space $\Tm$ obtained from a 
disjoint collection of Euclidean planes $\{W\}_{W\in Wall_\Tm}$
(called {\em walls}) and directed Euclidean strips\footnote{A 
direction for a strip $\St$ is an orientation for its
$\R$-factor $\St\simeq \R\times I$.} 
$\{\St\}_{\St\in Strip_\Tm}$ by isometric gluing\footnote{In 
general one may also have to complete the resulting quotient
space to get a Hadamard space.} subject to the following 
conditions:

1.  The boundary geodesics of each strip $\St\in Strip_\Tm$,
which we will refer to as {\em singular geodesics}, are
glued isometrically to distinct walls in $Wall_\Tm$.

2. Each wall $W\in Wall_\Tm$ is glued to at most two strips,
and the gluing lines are not parallel.

3. $\Tm$ is connected.

\noindent
One can think of $\Tm$ as sitting in $\R^3$ so that its
walls are parallel planes and the strips
meet the walls orthogonally.
Two walls $W_1,\,W_2\in Wall_\Tm$ are {\em adjacent}
if there is a strip $\St\in Strip_\Tm$ with 
$S\cap W_i\neq\emptyset$.  The incidence graph $Graph(\Tm)$ of $\Tm$
-- the graph with vertex set $Wall_\Tm$ and one edge for each pair of
incident walls --
is a graph isomorphic to a connected subcomplex of $\R$ with 
the usual triangulation (where the vertices are the integers).  
A wall is an {\em interior wall}
if it is incident to two strips, and a strip is an {\em interior
strip} if it is incident to two interior walls; $Wall^o_\Tm$ and 
$Strip^o_\Tm$ denote the interior walls and strips respectively.
For every interior wall $W\in Wall^o_\Tm$ we have a distinguished
point $o_W\defeq W\cap \St_1\cap \St_2$, where $\St_i\in Strip_\Tm$,
$i=1,\,2$,
are the strips incident to $W$.  Let $Strip^+_\Tm$ be the 
collection of oriented interior strips; an orientation of 
a strip $\St\in Strip_\Tm$
combines with the direction of $\St$ to give an orientation 
of the interval factor of $\St\simeq \R\times I$, and also an
ordering of the two incident walls.  We can define
a function $\eps:Strip^+_\Tm\ra\R$ as follows: if $W_-,\,W_+$
are incident to $\St^+\in Strip^+_\Tm$ and $W_-<W_+$
with respect to the ordering defined by $S^+$, then 
$\eps(\St^+)\in\R$ is defined to be the signed distance that
$o_{W_+}$ lies ``above'' $o_{W_-}$ in the strip $\St^+$.  We also
have a strip width function $l:Strip_\Tm\ra (0,\infty)$
and an angle function $\al:Wall^o_\Tm\ra (0,\pi)$
which give the angle between the oriented lines
$W\cap \St_i$ where $\St_1,\,\St_2$ are incident to $W$.

We will sometimes enumerate the consecutive walls and strips
of $\Tm$ so that $Wall_\Tm=\{W_i\}_{a<i<b}$ and
$Strip_\Tm=\{\St_i\}_{a<i<b-1}$ where $a\in\{-\infty,0\}$
and $b\in\N\cup\{\infty\}$.  We then define  $L_i^-\defeq
W_i\cap\St_{i-1}$ for $a+1<i<b$ and $L_i^+\defeq W_i\cap\St_i$
for $a<i<b-1$.

An {\em equivalence} between two templates $\Tm_1$ and $\Tm_2$
is an isometry $\psi:\Tm_1\ra\Tm_2$ which respects strip
directions.  Two templates are equivalent if and only if there is
an incidence preserving bijection $Wall_{\Tm_1}\cup
Strip_{\Tm_1}\ra Wall_{\Tm_2}\cup Strip_{\Tm_2}$ which
respects the functions $l,\eps$, and $\al$.
We will call a template {\em uniform} if there is a
$\frac \pi 2 \geq \beta>0$ so that the angle function
$\al:Wall^o_\Tm\ra(0,\pi)$ satisfies  $\pi-\beta\geq \alpha\geq \beta$,
and if the strip widths are bounded away from zero.  We are mostly
interested in uniform templates.  A
template $\Tm$ is {\em full} if $Graph(\Tm)\simeq\R$, {\em half}
if $Graph(\Tm)\simeq \R_+$,  and {\em finite}
if $|Wall_\Tm|<\infty$.

If $W\in Wall^o_\Tm$  and $\St_1,\,\St_2$ are the incident strips,
then the oriented lines $W\cap \St_1$ and $W\cap \St_2$ divide the plane $W$ into four
sectors which we call {\em Quarter Planes} and which we label as $Q_I$,
$Q_{II}$, $Q_{III}$, and $Q_{IV}$ as usual.  If we are given a choice
$Q_W$ of quarter planes in $W$ for each $W\in Wall_\Tm$  then
there is an isometric immersion ${\cal
D}:[\cup_{W\in Wall_\Tm}Q_W]\cup[\cup_{\St\in Strip_\Tm}\St]\ra \R^2$ (the development) which takes any
geodesic ray $\ga\subset \Tm$ such that $\ga\cap W\subset Q_W$, to a Euclidean ray (see Figure \ref{scale} in section \ref{selfsimsec} for an example of the developement of a special kind of template).

When $\Tm$ is a half template we will be primarily 
interested in geodesic rays $\gamma\subset\Tm$ that
start at a given base point and intersect all but finitely
many walls of $\Tm$.  From
the separation properties of walls it
is clear that such a ray intersects the walls $W\in Wall_\Tm$
in order.  We let $\geo^\infty \Tm$ (resp. $\tits^\infty \Tm$)
denote the corresponding subset of $\geo\Tm$ (resp. $\tits \Tm$).
In section \ref{pointorinterval} we will show that $\tits^\infty
\Tm$ is isometric to either a point, in which case $\Tm$ is called trivial,
 or an interval of length $<\pi$.

\begin{remark} One can show directly that any two half templates 
such that corresponding angles agree, and both corresponding strip widths and
displacements differ by a bounded amount will have $\geo^\infty$'s with the 
same Tits length.  We will only need a weaker version (that will follow from 
Theorem \ref{shadthm}) in this paper so we will not digress to prove it here.
\end{remark}

\subsubsection{Templates associated with itineraries in $T$}
\label{templatesection}

We now return to the setting of our paper: $G$ is an admissible
group with a discrete cocompact isometric action on a Hadamard
space $X$.  We now want to associate a template with each
geodesic segment/ray in $T$; these templates capture the asymptotic
geometry of geodesic segments/rays in $X$ which pass near the
corresponding edge spaces.

We first choose, in a $G$-equivariant way, a plane $F_e\subset Y_e$ for
each edge $e\in E$.  Then for every pair of adjacent edges $e_1,e_2$ we
choose, again equivariantly, a minimal geodesic from $F_{e_1}$
to $F_{e_2}$; by the convexity of $Y_v=\bar Y_v\times \R$, $v\defeq e_1\cap e_2$,
 this geodesic determines a Euclidean strip\footnote{$\St_{e_1,e_2}$ may have 
 width zero.} $\St_{e_1,e_2}\defeq \gamma_{e_1,e_2}
\times \R$ for some geodesic segment $\gamma_{e_1,e_2} \subset \bar Y_v$;
note that $\St_{e_1,e_2}\cap F_{e_i}$ is an axis of $Z(G_v)$.
Hence if $e,\,e_1,\,e_2\in E$, $e_i\cap e=v_i\in V$ are distinct
vertices, then  the angle between the geodesics
$\St_{e_1,e}\cap F_{e}$ and $\St_{e_2,e}\cap F_e$ 
is bounded away from zero (since only finitely many angles show
up).  We also note that Definition \ref{vertexspaces} tells us that $d_H(F_e,X_e)$ is bounded by $D''+Diam(\bar Y_e)$ which, 
since there are only finitely many $e$ up to the action of $G$, is uniformly bounded.

The significance of the strips $\St_{e_1,e_2}$ can be seen in the
next two lemmas.

\begin{lemma}
\label{wallstripwallqconvex}
There is a constant $C_3$ so that if $e_1=\ol{vv_1}$ and $e_2=\ol{vv_2}$
are adjacent
edges, then $X_{e_1}\cup\St_{e_1,e_2}\cup X_{e_2}$ and 
$X_{v_1}\cup\St_{e_1,e_2}\cup X_{v_2}$
are $C_3$-quasi-convex.
\end{lemma}
\proof
Since the Hausdorff distance $d_H(X_{\hat e},F_{\hat e})$
is uniformly bounded for $\hat e\in E$, 
 it suffices to show that there is a constant $C$ so that the unions
$F_{e_1}\cup\St_{e_1,e_2}\cup F_{e_2}$ are $C$-quasi-convex for 
all pairs of adjacent edges.  But if $e_1,\,e_2\in E$
are adjacent then $F_{e_1}\cup\St_{e_1,e_2}\cup F_{e_2}
\subset Y_v\simeq \bar Y_v\times\R$, and we are reduced
to showing that $\bar F_{e_1}\cup \ga_{e_1,e_2}\cup \bar F_{e_2}\subset\bar Y_v$
is uniformly quasi-convex, where $\bar F_{e_i}$ is 
the image of $F_{e_i}$ under the projection
$Y_v\simeq \bar Y_v\times\R\ra \bar Y_v$.  This follows
from Lemma \ref{geodesicunion}.  This gives the first statement.  

The interesting part of the second statement is when we consider $\ol{xy}$ when $x\in X_{v_1}$ and $y\in X_{v_2}$.  In this case Lemma \ref{findingz's} gives us $z_1\in N_{C_2}X_{e_1}$ and $z_2\in N_{C_2}X_{e_2}$ on $\ol{xy}$.  Now the first statement along with the convexity of $N_{C_2+D''}(X_{v_i})$ (note $z_i\in N_{C_2}X_{e_i}\subset N_{C_2+D''}X_{v_i}$) yields the second statement.

\qed

\begin{lemma}
\label{templateqconvexity}
There is a constant $C_4$ so that if $e_1,\ldots,e_n\in E$
is a geodesic edge path in $T$ with initial vertex $v_1$ and
terminal vertex $v_n$, then 
$$Z\defeq X_{e_1}\cup \St_{e_1,e_2}\cup X_{e_2}\cup\ldots
\cup X_{e_{n-1}}\cup \St_{e_{n-1},e_n}\cup X_{e_n}$$
and 
$$Z'\defeq X_{v_1}\cup \St_{e_1,e_2}\cup X_{e_2}\cup\ldots
\cup X_{e_{n-1}}\cup \St_{e_{n-1},e_n}\cup X_{v_n}$$
are $C_4$-quasi-convex.
\end{lemma}
\proof  
Pick $x,\,y\in Z$.  We may assume without loss of generality
that $x\in X_{e_1}\cup\St_{e_1,e_2}\cup X_{e_2}\subset X_v$ where
$v=e_1\cap e_2$, and $y\in X_{e_{n-1}}\cup\St_{e_{n-1},e_n}\cup X_{e_n}
\subset X_{v'}$ where $v'=e_{n-1}\cap e_n$.  Applying Lemma
\ref{findingz's} we get points
$z_i\in\ol{xy}\,\cap N_{C_2}(X_{e_i})$ for $2\leq i\leq n-1$,
with $d(z_i,x)\leq d(z_j,x)$ when $i\leq j$.  If $C_4\defeq C_2+C_3$, 
then by Lemma
\ref{wallstripwallqconvex} we have $\ol{xz_2}\subset
N_{C_4}(X_{e_1}\cup\St_{e_1,e_2}\cup X_{e_2})$,
$\ol{z_iz_{i+1}}\subset N_{C_4}(X_{e_i}\cup\St_{e_i,e_{i+1}}\cup
X_{e_{i+1}})$ for $i=2,\ldots n-1$, and $\ol{z_{n-1}y}\subset
N_{C_4}(X_{e_{n-1}}\cup\St_{e_{n-1},e_n}\cup X_{e_n})$. 

We omit the proof that $Z'$ is quasi-convex, as it is similar.
\qed

\medskip
Lemma \ref{templateqconvexity} suggests that we will understand the 
geodesic geometry of $X$ if the geometry of the sets $Z$
(as in the lemma) can be easily modeled.  To this end, we 
``approximate'' $Z$ with a template.

\begin{definition}
\label{ktemplate}
Suppose $\ga\subset T$ is a geodesic segment or ray.
Let $\Tm$ be a template with walls $\{W\}_{W\in Wall_\Tm}$
and strips $\{\St\}_{\St\in Strip_\Tm}$,
let $f:Wall_\Tm\ra E$ be an adjacency preserving bijection 
between the walls of $\Tm$ and the edges of $\ga$, and let $\phi:\Tm\ra X$ 
be a (not necessarily continuous) map. Then the triple $(\Tm, f,\phi)$ is a {\em $K$-template}
for $\ga$ if for all $W\in Wall_\Tm$ we have $\phi(W)\subset N_K(X_{f(W)})$
and $X_{f(W)}\subset N_K(\phi(W))$ and the following
conditions are met for every $\St\in Strip_\Tm$.

1. $Width(\St)\geq 1$.

2.  If $\St$ is incident to $W_1,\,W_2\in Wall_\Tm$ and
$x,\,y\in W_1\cup\St\cup W_2$ then we have
$|d_X(\phi(x),\phi(y))-d_\Tm(x,y)|<K$
and if $\ga_1:[0,1]\ra\ol{xy}$ and $\ga_2:[0,1]\ra\ol{\phi(x)\phi(y)}$
are constant speed parameterizations, then $d(\phi\circ\ga_1(t),\ga_2(t))<K$ for all $t\in [0,1]$.

3. If $W_1,\,W_2\in Wall_\Tm$ are adjacent to $\St$ then
the Hausdorff distance $d_H((X_{f(W_1)}\cup\St_{f(W_1),f(W_2)}\cup
X_{f(W_2)}),\phi(W_1\cup\St\cup W_2))<K$.
\end{definition}

Often when the value of $K$ is not relevant we will refer to the 
triple $(\Tm, f,\phi)$ as a template for $\ga$, by which we mean a $K$-template for some $K$.  

Let $\Tm'$ be another template such that the angles on corresponding walls agree with those of $\Tm$ and such that the other data (strip widths and displacements) differ from $\Tm$ by a bounded amount.  Then there is a natural (discontinuous) map $F:\Tm'\to \Tm$ which is an isometry on each wall and simply stretches the width of the strips.  It is easy to check that using $\phi'=\phi\circ F$ that we get a $K'$ template $(\Tm',f,\phi')$ for $\ga$ (see step 3 of the proof of Lemma \ref{itinscaleinvariant}).

For a suitably large $K$ we describe a construction that
gives a $K$-template for any  geodesic segment 
or geodesic ray $\ga\subset\Tm$. 
We will refer to these $K$-templates as {\em standard
$K$-templates}.  We begin with a disjoint collection of walls
$W_e$ and an isometry $\phi_e: W_e\ra F_e$ for each edge $e\subset\ga$.
For every pair $e,\,e'$ of adjacent edges of $\ga$, we let
$\hat\St_{e,e'}$ be a strip which is isometric to
$\St_{e,e'}\subset X$ if $Width(\St_{e,e'})\geq 1$,
and isometric to $\R\times[0,1]$ otherwise; we let
$\phi_{e,e'}:\hat\St_{e,e'}\ra\St_{e,e'}$ be an affine
map which respects product structure ($\phi_{e,e'}$ is
an isometry if $Width(\St_{e,e'})\geq 1$ and
compresses the interval otherwise).  We construct $\Tm$
by gluing the strips and walls so that the maps $\phi_e$ and $\phi_{e,e'}$
descend to continuous maps on the quotient.

The above construction yields

\begin{lemma}
\label{kbetaexist}
There is a constant $K=K(X)$ such that for every geodesic segment 
or ray, $\ga\subset T$, there is a $K$-template for $\ga$.  

There is a $\beta=\beta(X)>0$ such that for any $K$-template the 
angle function $\al:Wall^o_\Tm\ra(0,\pi)$ satisfies $0<\beta\leq\al\leq\pi-\beta<\pi$.  
\end{lemma}
\proof
We check that each condition of Definition
\ref{ktemplate} holds for the standard template
described above, for sufficiently large $K$.

First, since $\hd(F_e,X_e)$ is uniformly bounded
and $F_{f(W)}=\phi(W)$ for every $W\in Wall_\Tm$, we
have $\phi(W)\subset N_K(X_{f(W)})$
and $X_{f(W)}\subset N_K(\phi(W))$
for all $W\in Wall_\Tm$ for large enough $K$.

Conditions 1 and 3 follow immediately from the 
description of standard templates.

We now verify condition 2.  Pick adjacent walls
$W,W'\in Wall_\Tm$ and set $e\defeq f(W)$, $e'\defeq
f(W')$, and $v\defeq e\cap e'$.  Recall that
$F_e\cup\St_{e,e'}\cup F_{e'}\subset Y_v$
and $Y_v$ splits isometrically as 
$Y_v=\bar Y_v\times \R$ where $\bar Y_v$
is Gromov hyperbolic.  Furthermore, $\phi$
induces a map $W_e\cup\hat\St_{e,e'}\cup W_{e'}
\ra F_e\cup \St_{e,e'}\cup F_{e'}$ which is
compatible with the product structure.  Hence
condition 2 follows from part 3 of Lemma
\ref{geodesicunion} (and triangle inequalities)
when $d(F_e,F_{e'})\geq 4\de$
(where $\de$ is the maximum of the
hyperbolicity constants of the $\bar Y_v$'s);
modulo $G$ there are only finitely many 
cases when $d(F_e,F_{e'})<4\de$ (Lemma \ref{finitepairs}),
and each of these is also settled by part 3
of Lemma \ref{geodesicunion}.

For any $K$ template $(\Tm, f,\phi)$ for a $\ga$ containing an interior edge $e=\ol{v'v}$ we claim that the wall $W$ with $f(W)=e$ will have the 
same angle, up to taking supplements (i.e. $\alpha$ might be 
replaced by $\pi-\alpha$), as the angle $\alpha$ between the $\R$ 
factors of $Y_{v'}=\bar Y_{v'}\times \R$ and $Y_{v}=\bar Y_{v}\times \R$ 
in $Y_e = Y_v\cap Y_{v'}$.  The fact that these angles are positive 
and the finiteness of edges modulo $G$ will 
yield the result.  We note that $\alpha$ is the Tits angle between the $\R$ factors.

To see this we first note that Property 2 of Definition \ref{ktemplate} says that the angle for $W$, 
i.e. the angle between the gluing lines $L'$ and $L$, is the same as 
the comparison angle $\lim_{t\to\infty}\tilde \angle_{\phi(o)}\phi(L'(t)),\phi(L(t))$.  Also
Property 2 of Definition \ref{ktemplate} says that there are geodesic rays 
$\sigma=\lim_{i\to \infty} \ol{p\phi(L(t_i))}$ where $t_i\to \infty$ and $p\in Y_e$.    
Since $\phi(W)\subset N_KX_e$ we see that $\phi(L)\subset N_KX_e$ and hence $e\in \itin(\sigma)$.  Now if we let $e'$ be the other edge incident to $v$ in $\ga$ then the intersection of the wall $f^{-1}(e')$ with the strip between $W$ and $f^{-1}(e')$, is a line parallel to L and hence, again by 2 of Definition \ref {ktemplate}, $\phi(L)$ stays in a uniform neighborhood of $X_{e'}$ so $e'\in \itin(\sigma)$.  but $e,e'\in \itin(\sigma)$ implies by Lemma 
\ref{itindichotomy} that $\itin(\sigma)=\ol{Star(v)}$ and hence by Corollary \ref{staritin} 
any such $\sigma$ is asymptotic to  the $\R$ factor of  $Y_{v}=\bar Y_{v}\times \R$. 
Since $p\in Y_e\subset Y_v$, $\sigma$ is a half line of such an $\R$ and we assume 
without loss of generality that it points in the positive direction.  A similar argument works for $\phi(L')$.
Thus $\angle(L',L)= \lim_{t\to\infty}\tilde \angle_{\phi(o)}\phi(L'(t)),\phi(L(t))\leq \alpha$. 
Also the same arguments applied to $-L'$ and $L$ yield $\angle(-L',L)=\leq \pi-\alpha$.  
Thus we get equality and the result. 

\qed

The next proposition and Proposition \ref{templatetoambientpath}
are technical results that compare
template geometry with ambient geometry.

\begin{proposition}
\label{ambienttotemplatepath}
Suppose $K>0$.  There is a constant $C_5$ depending only on $K$ and the geometry
of $X$ with the following property.  Suppose $\ga\subset T$ 
is a geodesic segment or geodesic ray in $T$ with $i^{th}$
edge $e_i$, and set
$Z\defeq[\cup_{e\subset\ga}X_{e}]\cup
[\cup_{e,e'\subset\ga}\St_{e,e'}]$. If
$(\Tm,f,\phi)$ is a $K$-template for $\ga$, and $x,\,y\in Z$,
then there is a continuous map $\al:\ol{xy}\ra\Tm$
so that

1. $d(\phi\circ\al,id\restr_{\ol{xy}})<C_5$

2. For all $p,\,q\in \ol{xy}$ we have
\begin{equation}
\label{phidistortion}
d_X(p,q)-kC_5\leq d_\Tm(\al(p),\al(q))\leq length(\al\restr_{\ol{pq}})
\leq d_X(p,q)+kC_5
\end{equation} 
where the segment $\ol{\al(p)\al(q)}\subset \Tm$ intersects at most
$k-1$ strips and walls in $\Tm$.  In particular there are constants
$(L,A)$ depending only on $K$ and $X$ so that $\phi$
is an $(L,A)$ quasi-isometric embedding for every $K$-template
$(\Tm,f,\phi)$.
\end{proposition}

\proof
By the standard properties of Hadamard spaces we may reduce to the case that $x\in  Y_{e_x}\cup \St_{e_x,e_x'} \subset X_{v_x}$ ($v_x=e_x\cap e_x'$) and $y\in \St_{e_y,e_y'}\cup Y_{e_y'}\subset X_{v_y}$ since the original $x$ and $y$ are within a bounded distance of such.
Let $W_i\defeq f^{-1}(e_i)\in Wall_\Tm$ for $e_i$ $1\leq i\leq n$ the edges between $v_x$ and $v_y$.
We may apply Lemma \ref{findingz's} to the pair $x,\,y$
obtaining points $z_i\in\ol{xy}$.  We let $z_0=x$ and $z_{n+1} =y$. After making a small
perturbation of the $z_i$'s if necessary,  we may assume
that they satisfy $d_X(z_i,x)<d(z_j,x)$ when $i<j$.
For $1\leq i\leq n$ pick $w_i\in\Tm$ with $w_i\in W_i\subset\Tm$
with $d(z_i,\phi(w_i))\leq 1+\inf\{d(z_i,\phi(w))\mid w\in W_i\}\leq C_2+1$.  By Definition \ref{ktemplate} part 3 we can also choose $w_0\in W_0\cup \St \cup W_1$ and such that $d_X(x,\phi(w_0))\leq K$ and similarly choose $w_{n+1}$.  Now define $\al$ by the condition that $\al(z_i)=w_i$,
and $\al$ is a constant speed geodesic on the segment
$\ol{z_iz_{i+1}}$.

{\em Proof of 1.} Apply Definition \ref{ktemplate} to see that the 
constant speed parameterization $[0,1]\ra \ol{\phi(w_i)\phi(w_{i+1})}$
is at uniformly bounded distance from  the composition
of the constant speed parameterization 
$[0,1]\ra \ol{w_iw_{i+1}}\subset\Tm$ with
$\phi:\Tm\ra X$.  Since $d(\phi(w_i),z_i)$ is uniformly
bounded, we know that the constant speed parameterizations
$[0,1]\ra \ol{\phi(w_i)\phi(w_{i+1})}$ and $[0,1]\ra\ol{z_iz_{i+1}}$
are also at uniformly bounded distance from one another, so
there is a constant $c_1$ depending on $K$ 
so that $d(\phi\circ\al,id\restr_{\ol{xy}})<c_1$.

{\em Proof of 2.}
Assume $p\in \ol{z_{j-1}z_j}-z_{j-1}$ and $q\in\ol{z_{j'}z_{j'+1}}-z_{j'+1}$
for $j\leq j'$.  By Definition \ref{ktemplate} we have, for $c_2=2c_1+K$
$$ |length(\al\restr_{\ol{pz_j}})-d_X(p,z_j)|<c_2$$
$$ |length(\al\restr_{\ol{z_iz_{i+1}}})-d_X(z_i,z_{i+1})|<\mbox{$c_2$ for every $i=1,\ldots n-1$}$$
$$ |length(\al\restr_{\ol{z_{j'}q}})-d_X(z_{j'},q)|<c_2.$$
Hence there is a $c_3=c_3(K)$ so that 
$$length(\al_{\ol{pq}})=length(\al\restr_{\ol{pz_j}})+\ldots
+length(\al\restr_{\ol{z_{j'}q}})$$
$$\leq d_X(p,z_j)+\ldots +d_X(z_{j'},q)+(j'-j+2)c_2$$
$$\leq d_X(p,q)+kc_2.$$
To prove the remaining inequality of (\ref{phidistortion}) we 
break up the $\Tm$-geodesic $\ol{\al(p)\al(q)}$ into
at most $k$ subsegments $\ol{u_iu_{i+1}}$ so that each subsegment
lies in $W_i\cup\St_i\cup W_{i+1}$ for some $i$.
Then by definition \ref{ktemplate} we have
$|d_X(\phi(u_i),\phi(u_{i+1}))-d_\Tm(u_i,u_{i+1})|<K$
so 
$$d_X(p,q)\leq 2c_1+d_X(\phi\circ\al(p),\phi\circ\al(q))$$
$$\leq 2c_1+\sum d_X(\phi(u_i),\phi(u_{i+1}))$$
$$\leq 2c_1+kc_2+d_\Tm(\al(p),\al(q))$$
$$\leq d_\Tm(\al(p),\al(q))+kc_3.$$

where $c_3\defeq 2c_1+c_2$.  

To see the quasi-isometry property of $\phi$, let $x',y'\in \Tm$, $x=\phi(x')$, $y=\phi(y')$, and let $\alpha$ be the map defined above where we choose $w_0=x$ and $w_1=y$ (i.e $\alpha(x)=x'$ and $\alpha(y)=y'$).  Now \ref{phidistortion} applied to $p=x$ and $q=y$
along with $k\leq d_{\Tm}(x',y')+1$ (since strips have width at least 1) and $k\leq const_1 d_X(p,q)+const_2$ (as in the proof of the coarse lipschitz property of $\rho$ - see section \ref{itinsection})
yields the quasi-isometry property of $\phi$.
This completes the proof of Proposition 
\ref{ambienttotemplatepath}.
\qed

\begin{proposition}
\label{templatetoambientpath}
Pick $K>0$.  There is a constant $C_6=C_6(K,X)$ so that the following holds.
If $(\Tm,f,\phi)$ is a $K$-template for $\ga\subset T$,
and $x,\,y\in\Tm$, then there is a continuous map 
$\al:\ol{xy}\ra X$ where

1. $d(\al,\phi\restr_{\ol{xy}})<C_6$.

2.  For all $p,\,q\in\ol{xy}$ we have

\begin{equation}
d_\Tm(p,q)-kC_6\leq d_X(\al(p),\al(q))\leq length(\al\restr_{\ol{pq}})
\leq d_\Tm(p,q)+kC_6
\end{equation}
where the segment $\ol{pq}\subset\Tm$ intersects at most 
$(k-1)$ strips and walls in $\Tm$.
\end{proposition}
\proof
This is similar to the proof of Proposition
\ref{ambienttotemplatepath}, so we omit it.  
\qed

\begin{corollary}
\label{phidistortion2}
If $(\Tm,f,\phi)$ is a $K$-template and $C_6=C_6(K)$ is the constant
from Proposition \ref{templatetoambientpath}, then for
any $x,\,y\in\Tm$ we have 
$$d_\Tm(x,y)-kC_6\leq d_X(\phi(x),\phi(y))\leq d_\Tm(x,y)+kC_6$$
where $\ol{xy}\subset\Tm$ meets at most $k-1$ strips and walls.
\end{corollary}

\subsection{Shadowing}
\label{shadowingsection}

In this section we show that geodesic segments in a $K$-template
are sublinearly shadowed by ambient geodesic segments, and vice-versa.

\begin{theorem}[Shadowing]
\label{shadthm}
There is a function $\th:\R_+\ra\R_+$ 
depending on $K$ and the geometry of $X$ (sometimes denoted $\theta_{(X,K)}$) with 
$\lim_{R\ra\infty}\th(R)=0$ so that
if $(\Tm,f,\phi)$ is a $K$-template for a geodesic 
segment/ray $\ga\subset T$, then the following
hold.

1. If $x,\,y\in \Tm$, $z\in\ol{xy}$ and  $R\defeq d(z,x)$, then $d(\phi(z),\ol{\phi(x)\phi(y)})\leq
(1+R)\th(R)$.

2. If $x,\,y\in \Tm$, $z\in\ol{\phi(x)\phi(y)}$ and  $R\defeq d(z,\phi(x))$ then
$d(z,\phi(\ol{xy}))\leq (1+R)\th(R)$.

3. Let $\bar\Tm\defeq \Tm\cup\geo \Tm$ and $\bar X\defeq X\cup\geo X$
be the usual compactifications.  Then there is a unique topological
embedding $\geo\phi:\geo\Tm\ra \geo X$ so that 
$$\bar\phi\defeq \phi\cup\geo\phi:\bar\Tm\ra\bar X$$
is continuous at every $\xi\in\geo\Tm\subset\bar\Tm$.

4. The image of $\geo\phi$ is
$$[\cup_{W\in Wall_\Tm}(\geo X_{f(W)})]\cup[\geo^{\geo\ga}X],$$
and when $\ga$ is a ray with $\geo\ga=\eta$ then
$\geo\phi(\geo^\infty\Tm)=\geo^\eta X$ (see section \ref{itinsection} for the
definition of $\geo^\eta X$).

5. $\geo\phi\restr_{\geo^\infty\Tm}:\geo^\infty\Tm\ra\geo^\eta X$ is an isometric embedding with respect to the
Tits metric.
\end{theorem}

The proof of the theorem  breaks up into two pieces.  We first show
in Proposition \ref{wallcluster}
that a geodesic segment (in a template or in $X$) running through a sequence of consecutive
walls has to be ``close'' to any point $p$ which lies close to sufficiently
many walls in the sequence.  We then show in Theorem \ref{almostsqrtthm}
that a segment 
in a template (resp. in $X$) which
doesn't meet too many walls (i.e. encounters at most $Const\,\log R$ walls in the
segment $\ol{px}$, $d(p,x)=R$) is well shadowed by a geodesic
segment in $X$ (resp. in the template).  These two arguments 
are combined in section \ref{actualshadowing} to prove
Theorem \ref{shadthm}.

\subsubsection{Paths in a template which are close to a cluster of walls}

We begin with a result about templates. It estimates the excess length of a path
$\eta$ which connects two walls $W,\,W'$ 
while remaining outside a ball which intersects $W,\,W'$,
and all walls between them.

\begin{proposition}
\label{clusterprop}
Let $\Tm$ be a template with angle function 
$\al:Wall^o_\Tm\ra(0,\pi)$ satisfying $0<\beta\leq\al\leq\pi-\beta<\pi$.
Then there are positive constants $N_1=N_1(\beta),$ ($N_1\approx \frac{Const}{\beta}$), 
$C_1=C_1(\beta)$, and $C_2=C_2(\beta)$ 
with the following property.  Let $W_{n_0},\ldots,W_{n_1}\in  Wall_\Tm$
be a sequence of consecutive walls,
and suppose $W_i\cap B(p,R)\neq\emptyset$ for some $p\in\Tm$,
$R>0$ and every $n_0\leq i\leq n_1$.  Then for any $R'\geq N_1R$ 
and any path
$c:[0,1]\ra \Tm -B(p,R')$ with $c(0)\in W_{n_0}$ and $c(1)\in W_{n_1}$
we have $$length(c)\geq d_\Tm(c(0),c(1))+C_1(n_1-n_0-C_2)R'.$$ 
\end{proposition}

\proof Let $\St_i$ be the strip incident to $W_i$ and
$W_{i+1}$ for $i=n_0,\ldots,n_1-1$, set 
$L_i^-\defeq \St_{i-1}\cap W_i$ for $i=n_0+1,\ldots n_1$,
and set $L_i^+\defeq\St_i\cap W_i$ for $i=n_0,\ldots,n_1-1$. 

We prove the proposition with the help of some lemmas.

\begin{lemma}
\label{3walls}
There is a constant $c_1\approx\frac{Const}{\beta}$
so that if $p\in\Tm$ and $d(p,W_j)<R$ for $j=i\pm 1$
then $d(p,o_i)<c_1R$.
\end{lemma}

\proof
By joining $\ol{x_{i-1}p}$ to $\ol{px_{i+1}}$ for appropriate choices of 
$x_{i-1}\in W_{i-1}$ and $x_{i+1}\in W_{i+1}$ we get a path $\ga:[0,1]\ra\Tm$ 
of length at most $2R$ joining $W_{i-1}$ to $W_{i+1}$.  Therefore
there is a segment $[a,b]\subset[0,1]$ with $\ga(a)\in L_i^-$
and $\ga(b)\in L_i^+$.  So 
$$d(p,o_i)\leq R+\min(d(\ga(a),o_i),d(\ga(b),o_i))$$
$$\leq R+\frac{d(\ga(a),\ga(b))}{2\sin(\frac{\beta}{2})}\leq c_1R.$$
\qed

We now define $N_1\defeq\max (2c_1,[\frac{\pi}{\beta}]+2)$.

Consider a path $\eta:[0,1]\ra\Tm-B(p,R')$ where
$R'>N_1R$.  The ball $B(p,R')\subset \Tm$ is convex, 
so clearly $\Tm-B(p,R')$ is complete and locally compact with
respect to the induced path metric.  Therefore we may
assume that $c$ is a constant speed minimizing path from
$c(0)$ to $c(1)$ in $\Tm-B(p,R')$.  Since $\ol{B(p,R')}$ 
is a convex subset of the Hadamard space $\Tm$, the nearest
point projection $\Tm\ra\ol{B(p,R')}$ is distance non
increasing; it follows that  the set $c^{-1}(\ol{B(p,R')})$ is either empty
or a closed subinterval $[a,b]\subset [0,1]$.
$c\restr_{[0,a]}$ and $c\restr_{[b,1]}$ are constant speed
geodesic segments in the Hadamard space $\Tm$, and since
$R'>c_1R$ these segments lie in $\Tm-\{o_i\}_{n_0<i<n_1}$
by Lemma \ref{3walls}.

\begin{lemma}
$$c([0,a])\subset \left [\cup_{i=n_0}^{n_0+N_1}W_i\right ]\cup
\left [\cup_{i=n_0}^{n_0+N_1-1}\St_i\right ]$$
and 
$$c([b,1])\subset \left [\cup_{i=n_1-N_1}^{n_1}W_i\right ]\cup
\left [\cup_{i=n_1-N_1}^{n_1-1}\St_i\right ].$$
\end{lemma}

\proof
We prove the first assertion; the proof of the second is similar.
If the lemma were false, we would have $c(t)\in \St_{n_0+N_1}$
for some $t\in[0,a]$.   Therefore $c([0,t])$ must cross every strip
$\St_i$ for $n_0\leq i<n_0+N_1$, and for every
$n_0<i\leq n_0+N_1$ it must enter $W_i$ through
$L_i^-$ and exit through $L_i^+$.  Since 
$c([0,a])$ is disjoint from $\{o_i\}_{n_0<i<n_1}$
there is a flat convex strip $Y\subset\Tm$ containing
$c([0,a])$ in its interior.  Using $Y$ we can
define co-orientations for the segments
$c([0,a])\cap \St_i$ and $c([0,a])\cap W_i$
for $n_0\leq i\leq n_0+N_1$.  If two of the origins $o_i\in W_i$ for $n_0<i\leq n_0+N_1$ lie
on opposite sides of the corresponding
segments $c([0,a])\cap W_i$ with respect to
the co-orientations then the geodesic between them (of length less than 
$2c_1R$ by lemma \ref{3walls}) will intersect $c([0,a])$ and hence 
$d(o_i,c([0,a]))<c_1R$ for some $n_0<i\leq n_0+N_1$, and thus $d(p,c([0,a])<2c_1R$. 
But this cannot happen since $d(p,c([0,a])\geq R'>2c_1R$.
Thus all the origins
$o_i\in W_i$ for $n_0<i\leq n_0+N_1$ lie
on the same side of the corresponding
segments $c([0,a])\cap W_i$ with respect to
the co-orientations.  It follows that
the angle between $c([0,a])$ and $L_i^-$
increases by at least $\beta$ each time $c([0,a])$
passes through a wall.  Hence $(N_1-1)\beta<\pi$,
contradicting the definition of $N_1$.
\qed

\medskip
\no
{\em Proof of Proposition \ref{clusterprop} concluded.}  Let 
$[a',b']\subset[a,b]\subset[0,1]$ be the inverse
image of 
$$[\cup_{i=n_0+N_1+1}^{n_1-N_1-1}W_i]\cup[\cup_{n_0+N_1}^{n_1-N_1-1}\St_i]$$
under $c$.  We know that $c([a',b'])$ remains in the sphere
$S(p,R')$ while it passes through all the walls 
$W_i$ for $n_0+N_1<i<n_1-N_1$.  So for every $n_0+N_1<i<n_1-N_1$,
$c([a',b'])$ joins $L_i^-$ to $L_i^+$ outside
$B(p,R')\supset B(o_i,R'-c_1R)$.   Hence
$length(c([a',b'])\geq \beta(R'-c_1R)(n_1-n_0-(2N_1+2))\geq \frac \beta 2R'(n_1-n_0-(2N_1+2))$
while $d_\Tm(c(a'),c(b'))\leq 2R'$
so $length(c[a',b'])\geq d_\Tm(c(a'),c(b'))+ \frac \beta 2(n_1-n_0-(2N_1+2)-\frac 4 \beta)R'$ and hence
$$length(c)\geq d_\Tm(c(0),c(1))+C_1(n_1-n_0-C_2)R'$$
where $C_1,\,C_2$ depend only on $\beta$.\qed

\begin{corollary}
\label{clustercorollary}
Let $\Tm$, $N_1$, $C_1$, $C_2$, $W_{n_0},\ldots,W_{n_1}$, $p$, $R$ be as in 
Proposition \ref{clusterprop}.  If $n_1-n_0>C_2$,
 then any geodesic segment from $W_{n_0}$
to $W_{n_1}$ must pass through $B(p,N_1R)$.
\end{corollary}
The result corresponding to Corollary \ref{clustercorollary} in the space $X$ is:

\begin{proposition}
\label{wallcluster}
There are constants $N_2=N_2(X),\,R_0=R_0(X)$ with the following
property.  If $n\geq N_2$, $e_1,\ldots e_n\in E$ are 
the consecutive edges of a  geodesic segment, $\ga$, in the tree $T$, $p\in X$,
$R\geq R_0$, and $X_{e_i}\cap B(p,R)\neq\emptyset$
for $1\leq i\leq n$; then for any $C\geq 0$, and
any segment $\ol{xy}\subset X$
with $\ol{xy}\cap N_C(X_{e_i})\neq \emptyset$ for $i=1$ and $i=n$,
we have $\ol{xy}\cap B(p,N_2R+2C)\neq\emptyset$.
\end{proposition}

\proof Let $e_1,\ldots,e_n$, $p$, $X_{e_i}$ be as in the
statement of the proposition.  If $\ol{xy}\cap N_C(X_{e_i})\neq\emptyset$
for $i=1$ and $i=n$, then we have $x_0\in X_{e_1}$ and 
$y_0\in X_{e_n}$ with $d(x_0,\ol{xy}),d(y_0,\ol{xy})\leq C$.
By convexity of the distance function $d_X$ it suffices
to show that $\ol{x_0y_0}\cap B(p,N_2R)\neq\emptyset$. 

Let $K$ and $\beta$ be as in Lemma \ref{kbetaexist} and $(\Tm,f,\phi)$ 
be a $K$ template for $\ga$ whose angles are bounded by $\beta$. 
Let $\alpha:\ol{x_0y_0}\to \Tm$ be the map guaranteed by Proposition 
\ref{ambienttotemplatepath}.  So $Length(\alpha)\leq d(x_0,y_0)+nC_5$.  
By part 3 of definition \ref{ktemplate}
there is a $p'\in \Tm$ be such that $d(p,\phi(p'))<R+K$ and hence (since $X_{e_i}\subset N_K\phi(W_i)$) $d(\phi(p'),\phi(W_i))<2R+2K$.  Since, by 
Proposition \ref{ambienttotemplatepath}, $\phi$ is an $(L,A)$-quasi 
isometric embedding we have $B(p',R_2)\cap W_i\neq \emptyset$ for 
$R_2= L(2R+2K)+A$.
We will choose $R_0$ and $N_2$ large enough so that for $n>N_2$ 
and $R\geq R_0$ we will have $nC_5<C_1(n-C_2)N_1R_2$ for the 
$N_1(\beta)$, $C_1(\beta)$ and $C_2(\beta)$ of Proposition \ref{clusterprop}.    
Thus Proposition \ref{clusterprop} forces $\alpha$ to intersect $B(p',N_1R_2)$.  
So we conclude (again using Proposition \ref{ambienttotemplatepath}) that 
$d(\ol{x_0y_0},p)<C_5+d(\phi(\alpha),\phi(p'))+R+K\leq C_5+LN_1R_2+A+R+K$.  
The proposition now follows by taking $N_2$ and $R_0$ large enough.
\qed

\subsubsection{Paths with small length distortion}

\begin{proposition}
\label{almostsqrtthm}
Pick $M>0$ and $\al\in(\frac{1}{2},1]$.  Then there is
a constant $C=C(M,\al)$ so that if $1\leq A\leq B$,
$\eta:[A,B]\ra X$ is a (not necessarily continuous) map
to a Hadamard space $X$, and for all $A\leq t_1\leq t_2\leq B$
we have
\begin{equation}
\label{logdistortion}
|d_X(\eta(t_1),\eta(t_2))-(t_2-t_1)|\leq M(1+log(\frac{t_2}{t_1}))
\end{equation}
then 
\begin{equation}
\label{almostsqrt}
d(\eta(t),z)\leq C(1+t^\al)
\end{equation}
where $z\in\ol{\eta(A)\eta(B)}$ is the point with
$d(z,\eta(A))=\frac {t-A}{B-A}d(\eta(A),\eta(B))$.  Similarly,
if $A\geq 1$ and $\eta:[A,\infty)\ra X$ satisfies
(\ref{logdistortion}) for all $A\leq t_1\leq t_2$,
then there is a unique unit speed geodesic ray
$\ga:[A,\infty)\ra X$ with $\ga(A)=\eta(A)$ such that
$$d_X(\eta(t),\ga(t))\leq C(1+t^\al)$$
for all $t\in[A,\infty)$.
\end{proposition}

\proof
First note that we may assume that $A=1$, since the
map $\eta_1:[1,B-A+1]\ra X$ given by
$\eta_1(t)\defeq \eta(t+A-1)$ will satisfy the hypotheses of
the proposition, and the conclusion of the proposition applied
to $\eta_1$ will imply (\ref{almostsqrt})
for $\eta$.

{\em Step 1: When $1\leq s_1\leq s_2\leq 2s_1\leq B$ and 
$s_1$ is sufficiently large then the comparison angle
$\cangle_{\eta(1)}(\eta(s_1),\eta(s_2))\leq Const\, s_1^{\al-1}$.}

We will make use of the following lemma that follows from standard comparisons.

\begin{lemma}
\label{excesslemma}
Let $X$ be a Hadamard space, $x,\,y,\,z\in X$.  Set
$L\defeq d(x,z)$, and the excess $E\defeq d(x,y)+d(y,z)-d(x,z)$.
Then 
\begin{equation}
d(y,\ol{xz})\leq \frac{\sqrt{2LE}}{2}\sqrt{1+\frac{E}{2L}}
\end{equation}
\begin{equation}
\label{simpleest}
\leq \sqrt{LE} \mbox{ if $E\leq 2L$.}
\end{equation}
\end{lemma}
\proof Triangle comparison.\qed

\no
{\em Proof of Proposition \ref{almostsqrtthm} continued.}
Take $1\leq s_1\leq s_2\leq B$, and consider the triple
$\eta(1),\,\eta(s_1),\,\eta(s_2)$.  The excess for the triple
is
$$\leq M(1+\log s_2)+M(1+\log s_1)+M(1+\log(\frac{s_2}{s_1}))= 3M(1+\log s_2).$$
Since $d(\eta(1),\eta(s_2))\geq (s_2-1)-M(1+\log s_2)$ when 
$s_2> c_1= c_1(M)$, then the excess is $\leq 2d(\eta(1),\eta(s_2))$.
Thus since $d(\eta(1),\eta(s_2))\leq (s_2-1)+M(1+\log s_2)$ applying (\ref{simpleest}) we get
$$d(\eta(s_1),\ol{\eta(1)\eta(s_2)})\leq \sqrt{[(s_2-1+M(1+\log s_2))][3M(1+\log s_2)]}$$
\begin{equation}
\label{firstexcess}
\leq
c_2(1+s_2^\al)
\end{equation}
where $c_2=c_2(M,\al)$.  Therefore the comparison angle
$\cangle_{\eta(1)}(\eta(s_1),\eta(s_2))$ satisfies
$$\sin(\cangle_{\eta(1)}(\eta(s_1),\eta(s_2)))\leq \frac{c_2(1+s_2^\al)}{d(\eta(1),\eta(s_1))}
\leq \frac{c_2(1+s_2^\al)}{[(s_1-1)-M(1+\log s_1)]}.$$
So there are constants $c_3=c_3(M,\al)$ and $c_4=c_4(M,\al)$
so that if $c_3\leq s_1\leq s_2\leq 2s_1\leq B$
then 
\begin{equation}
\sin(\cangle_{\eta(1)}(\eta(s_1),\eta(s_2)))\leq \frac{2c_2(1+s_2^\al)}{s_1}
\leq\frac{c_4}{2}s_2^{\al-1}
\end{equation}
and 
\begin{equation}
\label{sinest}
\cangle_{\eta(1)}(\eta(s_1),\eta(s_2))\leq c_4s_2^{\al-1}.
\end{equation}

{\em Step 2: Estimating $d(\eta(t),\ol{\eta(1)\eta(B)})$.}
Now pick $t_0\in [1,B]$ with $t_0\geq c_3$ .
Let $t_i=2^it_0$ for $i=0,\ldots,n$ where $n$ is the integer
part of $\frac{\log B}{\log 2}$, and $t_{n+1}=B$.
Then  $\frac{B}{t_n}<2$.  Applying the estimate (\ref{sinest}) 
with $u_i\defeq \eta(t_i)$ we have
for $i=1,\ldots,n+1$:
$$\cangle_{\eta(1)}(u_{i-1},u_i)\leq c_4t_i^{\al-1}$$
since $\frac{t_i}{t_{i-1}}\leq 2$ and $t_i\geq c_3$.
Either $n=0$, in which case we have 
$$d(\eta(t_0),\ol{\eta(1)\eta(B)})\leq c_2(1+t_1^\al)\mbox{\quad by (\ref{firstexcess})}$$
$$\leq c_2(1+2^\al t_0^\al)\leq c_5(1+t_0^\al)$$
where $c_5=c_5(M,\al)$.
Otherwise, if $t_0>c_6=c_6(M,\al)$,
$$d(\eta(1),u_0)\leq (t_0-1)+M(1+\log t_0)$$
$$\leq (t_i-1)-M(1+\log t_i)\leq d(\eta(1),u_i).$$
So for $0\leq i\leq n+1$ we may pick
$v_i\in\ol{\eta(1)u_i}$ with $d(\eta(1),v_i)=d(\eta(1),u_0)=:R_0$.
By triangle comparison we have
$$d(v_{i-1},v_i)\leq R_0\cangle_{\eta(1)}(v_{i-1},v_i)\leq R_0\cangle_{\eta(1)}(u_{i-1},u_i)$$
$$\leq R_0c_4t_i^{\al-1}$$
so
$$d(u_0,\ol{\eta(1)\eta(B)})=d(u_0,\ol{\eta(1)u_{n+1}})\leq d(u_0,v_{n+1})$$
$$\leq \sum_{i=1}^{n+1}d(v_{i-1},v_i)\leq R_0c_4\sum_{i=1}^{n+1}t_i^{\al-1}$$
$$\leq R_0c_7t_0^{\al-1}$$
for $c_7=c_7(M,\al)$ (independent of $n$).
We have
$R_0\leq t_0-1+M(1+\log t_0)\leq 2t_0$ when $t_0$ is sufficiently
large, so when $t_0\geq c_8=c_8(M,\al)$ 
$$d(u_0,\ol{\eta(1)\eta(B)})\leq c_9t_0^\al$$
where $c_9=c_9(M,\al)$.

{\em Step 3:  Estimating $d(\eta(t),z)$ where $z\in\ol{\eta(1)\eta(B)}$ satisfies
$d(z,\eta(1))=\frac {t-1}{B-1}d(\eta(1),\eta(B))$.}  Let $x\in\ol{\eta(1)\eta(B)}$ be
the point nearest $\eta(t)$.  Then for $t\geq c_{10}$
we have 
$$d(\eta(t),z)\leq d(\eta(t),x)+d(x,z)$$
$$= d(\eta(t),x)+|d(x,\eta(1))-(t-1)\frac {d(\eta(1),\eta(B))}{B-1}|$$
$$\leq d(\eta(t),x)+[d(x,\eta(t))+\frac {(t-1)}{B-1}M(1+\log B)+M(1+\log t)]\leq 3c_9t^\al.$$

Thus the result holds for large $t$, but by choosing $C$ large enough we get the result for all $t$.

For the ray case we redo step 3 above when $d(\eta(1),\eta(B))>t$ for $z_1\in\ol{\eta(1)\eta(B)}$ the point which satisfies
$d(z_1,\eta(1))=t$.  Let $x\in\ol{\eta(1)\eta(B)}$ be
the point nearest $\eta(t)$.  Then for $t\geq c_{11}$
we have 
$$d(\eta(t),z_1)\leq d(\eta(t),x)+d(x,z_1)$$
$$= d(\eta(t),x)+|d(x,\eta(1))-(t-1)|$$
$$\leq d(\eta(t),x)+[d(x,\eta(t))+M(1+\log t)]\leq 3c_9t^\al.$$
 
Now choose $\ga$ as a limit of a subsequence of $\ol{\eta(1),\eta(B)}$ as $B$ goes to $\infty$ 
(the result will in fact imply that the sequence itself converges).  
Since none of the constants depended on $B$ the above estimate 
for $d(\eta(t),z_1)$ gives the result for rays. 
\qed

\subsubsection{The proof of Theorem \ref{shadthm}}
\label{actualshadowing}

Let $\beta=\beta(X)$ be the minimum angle between singular 
geodesics of a template for $X$ defined in Lemma \ref{kbetaexist}, 
and let $C_2(\beta),\,N_1(\beta)$
be the constants from Corollary \ref{clustercorollary}.
Let $N_2,\,R_0$ be the constants from Proposition \ref{wallcluster};
we will assume that $N_2\geq \max(C_2(\beta),N_1(\beta))$.

\begin{definition}
For every $\psi>0$ , we will say that $z\in\ol{xy}$
is a {\em $\psi$-cluster point} if the segment $\ol{B(z,\psi d(z,x))}\cap\ol{xy}$
intersects at least $N_2$ walls of $\Tm$, or if
$z\in\{x,y\}$.  We will let $P_\psi\subset\ol{xy}$
represent the set of $\psi$-cluster points.
\end{definition}

{\em Proof of 1.}
Let $(\Tm,f,\phi)$, $x$, $y$, be as in the statement of part
1 of the theorem.  We begin with two lemmas.

\begin{lemma}
\label{mustcrosswall}
There is a constant $c_1$ depending on $K$ such that if
$\ol{xy}$ intersects a wall $W\in Wall_\Tm$, then 
$\ol{\phi(x)\phi(y)}\cap N_{c_1}(X_{f(W)})\neq\emptyset.$
\end{lemma}

\proof 
If $\{x,y\}\cap W\neq\emptyset$ this is immediate 
since $\phi(W)\subset N_K(X_{f(W)})$ by definition
\ref{ktemplate}.  Otherwise $x$ and $y$ must lie
in distinct components of $\Tm-W$ because each component
is convex.  Say $x\in W_1\cup \St_1\cup W_1'$ and 
$y\in W_2\cup\St_2\cup W_2'$ where $W_i$ is adjacent to
$W_i'$.  Set $v\defeq f(W_1)\cap f(W_1')$ and
$v'\defeq f(W_2)\cap f(W_2')$.  Applying
Lemma \ref{findingz's} the lemma follows.
\qed

\begin{lemma}
\label{smalldistapplies}
Suppose $p_1,\,p_2\in\ol{xy}$, $d(p_1,x)\leq d(p_2,x)$,
and $\ol{p_1p_2}\cap P_{\psi}\subset \{p_1,p_2\}$.
Set $A\defeq 1+d(p_1,x)$, $B\defeq 1+d(p_2,x)$. Let
$\eta_0:[A,B]\ra \Tm$ be the unit speed parameterization
of $\ol{p_1p_2}$.  Then the
composition $\eta\defeq\phi\circ\eta_0:[A,B]\ra X$
satisfies the hypotheses of Proposition \ref{almostsqrtthm}
with $M=M(K,\psi)$.
\end{lemma}
\proof
Pick $A\leq t_1<t_2\leq B$.  By Corollary \ref{phidistortion2}
we have
\begin{equation}
\label{distortionestimate}
|d_X(\eta(t_1),\eta(t_2))-(t_2-t_1)|\leq kC_6
\end{equation}
where $\ol{\eta_0(t_1)\eta_0(t_2)}$ intersects at most
$(k-1)$ walls and strips of $\Tm$.  For each 
$z\in\ol{\eta_0(t_1)\eta_0(t_2)}$, the segment
$\ol{B(z,\psi d(z,x))}\cap\ol{xy}$ intersects at most
$N_2-1$ walls. We can cover
$\ol{\eta_0(t_1)\eta_0(t_2)}-B(x,1)$ with
at most $Const\,\log(\frac{t_2}{t_1})$ such segments,
and $\ol{B(x,1)}\cap \ol{\eta_0(t_1)\eta_0(t_2)}$
intersects at most $2$ walls (since strips have width at least 1), so the lemma follows from 
(\ref{distortionestimate}).
\qed

\medskip
{\em Step 1: Estimating $d(\phi(z),\ol{\phi(x)\phi(y)})$
when $z\in P_\psi$.}
Suppose $z\in P_\psi$.  Then $\ol{xy}\cap B(z,\psi R)$ intersects
at least $N_2$ consecutive walls $W_1,\ldots,W_n$ of $\Tm$.
By definition \ref{ktemplate} and the fact that $\phi:\Tm\ra X$
is a quasi-isometric embedding (Proposition \ref{ambienttotemplatepath})
there is a constant $c_2=c_2(K)$ so that
$d_X(\phi(z),X_{f(W_i)})\leq c_2(1+\psi R)$ for $i=1,\ldots n$.
By Lemma \ref{mustcrosswall} and Proposition \ref{wallcluster},
we conclude that $\ol{\phi(x)\phi(y)}\cap B(\phi(z),N_2c_2(1+\psi R)+2c_1)\neq
\emptyset$ provided $c_2(1+\psi R)\geq R_0$.  So there 
are positive constants $r_0=r_0(K,\psi), \,c_3=c_3(K,\beta)$ so that if $z\in P_\psi$
and $R\defeq d(z,x)\geq r_0$,
then 
\begin{equation}
\label{psiclusterclose}
d(\phi(z),\ol{\phi(x)\phi(y)})\leq c_3\psi R.
\end{equation}
Hence for any $z\in P_\psi$
\begin{equation}
\label{psiclusterclose1}
d(\phi(z),\ol{\phi(x)\phi(y)})\leq c_3\psi R+Lr_0+A=c_3\psi R+c_4
\end{equation}
where $(L,A)$ are the quasi-isometric embedding constants
of $\phi$, and $c_4\defeq Lr_0+A=c_4(K,\psi)$.

{\em Step 2: Estimating $d(\phi(z),\ol{\phi(x)\phi(y)})$ when
$z\not\in P_\psi$.} Pick $z\not\in P_\psi$.  There are points
$p_1,\,p_2\in \ol{xy}$ so that $z\in\ol{p_1p_2}$,
$d(p_1,x)\leq d(p_2,x)$, and either $\ol{p_1p_2}\cap P_\psi
=\{p_1,\,p_2\}$ or one or both of $p_1=x$ or $p_2=y$ hold.  Each step of the argument below holds in the special cases $p_1=x$ or $p_2=y$ (often for easier reasons) so we will ignore them.
By (\ref{psiclusterclose1}) we have
\begin{equation}
\label{pixsmall}
d(\phi(p_i),\ol{\phi(x)\phi(y)})\leq c_3\psi d(p_i,x)+c_4
\end{equation}
for $i=1,2$.  
Since $p_1,\,p_2$ satisfy the conditions of Lemma
\ref{smalldistapplies}, we may apply Proposition
\ref{almostsqrtthm} with $\alpha = \frac 3 4$ to get that for the 
$w\in \ol{\phi(p_1)\phi(p_2)}$ with 
$d(w,\phi(p_1))=\frac {d(z,p_1)}{d(p_1,p_2)}d(\phi(p_1),\phi(p_2))$ we have
\begin{equation}
\label{closetop1p2}
d(\phi(z),w)\leq C(1+[d(z,x)]^{\frac{3}{4}})\leq \psi d(z,x) +c_6.
\end{equation}
Where $c_6=c_6(K,\psi)$.

On the other hand for $\zeta=\frac {d(w,\phi(p_1))}{d(\phi(p_1),\phi(p_2))}=\frac {d(z,p_1)}{d(p_1,p_2)}$ 
the convexity of the distance function $d(*,\ol {\phi(x)\phi(y)})$ and \ref{pixsmall} says that 
\begin{equation}
\label{wclose}
d(w,\ol{\phi(x)\phi(y)})\leq c_3\psi ( (1-\zeta)d(p_1,x)+ \zeta d(p_2,x))+c_4
= c_3\psi d(z,x) +c_4
\end{equation}

Combining (\ref{wclose}) and (\ref{closetop1p2}) we get
\begin{equation}
\label{noncluster}
d(\phi(z),\ol{\phi(x)\phi(y)})\leq c_7\psi d(z,x)+c_8
\end{equation}
for $c_7=c_7(K)$ and $c_8=c_8(K,\psi)$.  Thus (\ref{psiclusterclose1}) and
(\ref{noncluster}) together imply that for every choice of $\psi>0$ 
there are constants $c_9(K,\beta)$ and $c_{10}(K,\psi)$ such that for all $z$
$$d(\phi(z),\ol{\phi(x)\phi(y)})\leq c_9\psi R+c_{10}.$$
The fact that $c_9$ does not depend on $\psi$ allows us to pick $\theta(R)$ 
which decays to 0 as $R\to \infty$ such that for each $R$ there is a 
$\psi=\psi(R)>0$ (for example choose $\psi(R)$ decaying to 0 such 
that $c_{10}(\psi)<R^{-\frac 1 2}$) so that 
$$c_9\psi R+c_{10}(\psi)\leq (1+R)\theta (R)$$
This implies part 1 of Theorem \ref{shadthm}.

\medskip
{\em Proof of 2.} We omit this, as it is similar to the proof of 1.

\medskip
{\em Proof of 3.} Part 3 follows directly from part 1 and 
Lemma \ref{sublinearbending}.

\medskip
{\em Proof of 4.} Part 2 along with reasoning similar to the proof of Lemma 
\ref{sublinearbending} shows that if $x,\,y_k\in\Tm$, $\xi\in \geo X$,
and $\ol{\phi(x)\phi(y_k)}\ra \ol{x\xi}$, then $\ol{xy_k}$ 
converges to some ray $\ol{x\xi'}\subset\Tm$.  Hence 
$\geo\phi$ is  a homeomorphism from $\geo \Tm$ onto the limit
set  of $\phi(\Tm)\subset X$, which we denote by $\geo\phi(\Tm)$.
It follows immediately from this that $\geo\phi$ maps 
$\geo\Tm-\geo^\infty\Tm=\cup_{W\in Wall_\Tm}\geo W$
homeomorphically onto $\cup_{W\in Wall_\Tm}\geo X_{f(W)}$.  
This implies 4 when $\ga\subset T$ is a geodesic segment,
so we now assume that $\ga$ is a geodesic ray, and we need only
show that $\geo \phi(\geo^\infty\Tm)=\geo^\eta X$ where
$\eta=\geo\ga\in\geo T$.  We will let $W_k\in  Wall_\Tm$ be the $k^{th}$ wall of $\Tm$.

We first show $\geo\phi(\geo^\infty\Tm)\subset\geo^\eta X$.
Suppose $\xi\in\geo^\infty\Tm$, $z_k\in\ol{x\xi}\cap W_k$, 
and $d(z_k,x)\ra\infty$. Thus $\ol{\phi(x)\phi(z_k)}\ra\ol{\phi(x)
\geo\phi(\xi)}$.  By lemma \ref{findingz's} part 1, for any $l$
we have $\ol{\phi(x)\phi(z_k)}\cap N_c(X_{f(W_l)})\neq
\emptyset$ for sufficiently large $k$.  So by Lemma
\ref{bdyconvex} either $\geo\phi(\xi)\in\geo N_c(X_{f(W_l)})
=\geo X_{f(W_l)}$ or $\ol{\phi(x)\geo\phi(\xi)}\cap N_c(X_{f(W_l)})
\neq\emptyset$.  The former case is impossible since
we already know that $(\geo\phi)^{-1}(\geo X_{f(W_l)})
=\geo W_l$.  So $\ol{\phi(x)\geo\phi(\xi)}\cap N_c(X_{f(W_l)})\neq
\emptyset$ for every $W\in Wall_\Tm$, forcing $\itin(\ol{\phi(x)\geo\phi(\xi)})=\ga$.  
So $\geo\phi(\geo^\infty\Tm)\subset\geo^\eta X$.

We now show $\geo^\eta X\subset\geo\phi(\geo^\infty\Tm)$.  
Suppose $x\in W_1\in Wall_\Tm$ and $\xi\in\geo^\eta X$.  
Then from
the definition of itineraries \ref{itindef} and Lemma \ref{itindichotomy} 
there are $z_k\in\ol{\phi(x)\xi}\cap N_C(X_{e_k})$ for 
$e_k\in T$ so that $d_T(e_k,f(W_1))\ra\infty$ and
$d_T(e_k,\ga)$ is uniformly bounded.  So for all but finitely
many $W$, $f(W)$ separates $e_k$ from $f(W_1)$ for sufficiently
large $k$, so by Lemma \ref{findingz's}
$$\ol{\phi(x)z_k}\cap N_c(X_{f(W)})\neq \emptyset$$
for sufficiently large $k$.  Hence $\xi$ belongs
to the limit set of $\cup_{W\in Wall_\Tm}X_{f(W)}$, which
is the same as $\geo(\phi(\Tm))=\geo\phi(\geo\Tm)$. 
Therefore $\xi\in(\geo\phi)(\geo\Tm-\cup_{W\in Wall_\Tm}
\geo W)=(\geo\phi)(\geo^\infty\Tm)$.

\medskip
{\em Proof of 5.}  Pick $x\in \Tm$ and $\xi\in\geo^\infty\Tm$.
Suppose there is a sequence $z_k\in\ol{x\xi}$ with
$\lim_{k\ra\infty}(z_k,x)=\infty$ so that $z_k$ is a $\psi_k$-cluster
point $\ol{x\xi}$ where $\psi_k\ra 0$.  Set $R_k\defeq d(z_k,x)$.
Applying Corollary \ref{clustercorollary} to each $z_k$,
we see that if $\xi'\in\geo^\infty\Tm$, then for sufficiently
large $k$ the intersection 
$\ol{x\xi'}\cap B(z_k,N_1\psi_kR_k)$ is nonempty, and this
 clearly forces $\ol{x\xi'}=\ol{x\xi}$.  In this case we have $\geo^\infty\Tm
 =\{\xi\}$, and 5 is immediate.  So we may assume that for
 every $\xi_1,\,\xi_2\in\geo^\infty\Tm$ there is a $\psi>0$
 such that $\ol{x\xi_1}$ and $\ol{x\xi_2}$ contain no
 $\psi$-cluster points.  Let $\bar\eta_i:[1,\infty)\ra \Tm$
 be the unit speed parameterization of $\ol{x\xi_i}$, and 
 let $\eta_i:[1,\infty)\ra X$ be the composition $\phi\circ\bar\eta_i
:[1,\infty)\ra X$.  Then by Lemma \ref{smalldistapplies} the $\eta_i$
satisfy the hypotheses of Proposition \ref{almostsqrtthm}
for a suitable $M$; so we have  unit speed geodesic rays
$\ga_i:[1,\infty)\ra X$ with $\ga_i(1)=\eta_i(1)$
and $d_X(\eta_i(t),\ga_i(t))\leq C(1+t^{\frac{3}{4}})$.
Now, for sufficiently large $k$, choose $t^i_k$ so that 
$\bar\eta_i(t_k^i)$ lies in the $k^{th}$ wall of 
$\Tm$.  By definition \ref{ktemplate} 
$$|d_X(\eta_1(t_k^1),\eta_2(t_k^2))-d_\Tm(\bar\eta_1(t_k^1),
\bar\eta_2(t^2_k))|<K.$$
Thus $$\tangle
(\ga_1,\ga_2)=\lim_{k\ra\infty}
\cangle_{\eta_1(1)}(\eta_1(t^1_k),\eta_2(t^2_k))=\lim_{k\ra\infty}
\cangle_x(\bar\eta_1(t^1_k),\bar\eta_2(t^2_k))=\tits(\xi_1,\xi_2).$$
This proves 5. 
\qed

\subsubsection{Applications of Theorem \ref{shadthm}}

Theorem \ref{shadthm} has a number of  corollaries:

\begin{corollary}
\label{geoetadetection}
Let $G$ be the fundamental group of an admissible graph of 
groups ${\cal G}$, and let $G\acts T$ be the Bass-Serre tree
of ${\cal G}$.     Let $\ga\subset T$ be a geodesic ray with 
$i^{th}$ edge $e_i\subset T$, and set $\eta\defeq\geo \ga\in\geo T$.
Then for any admissible action $G\acts X$, 
the subset $\geo^\eta X\subset\geo X$ defined in 
section \ref{itinsection} is precisely the set of limit points of the
sequence of subsets $\geo X_{e_k}\subset\geo X$.
\end{corollary}
\proof  Let $(\Tm,f,\phi)$ be  a template for $\ga\subset T$.  By parts
3 and 4 of Theorem \ref{shadthm}, $\geo\phi:\geo\Tm\ra\geo X$ is 
a homeomorphism onto $(\cup_k\,\geo X_{e_k})\cup(\geo^\eta X)$.
Therefore it suffices to show that $\geo^\infty\Tm\subset\geo\Tm$
is the set of limit points of the sequence $\geo f^{-1}(e_k)\subset\geo\Tm$.
To see this, observe that any geodesic segment $\ol{px}\subset\Tm$
which arrives at a wall $W\in Wall_\Tm$ via a strip $\St\in Strip_\Tm$
may be prolonged by a geodesic ray contained in $W$; this 
implies that every $\xi\in\geo^\infty\Tm$ is a limit of a sequence
$\xi_k\in\geo f^{-1}(e_k)$.  On the other hand, if $\xi_k\in\geo f^{-1}(e_k)$
converges to $\xi\in\geo\Tm$, then Lemma \ref{bdyconvex} implies
that either $\xi\in\geo^\infty\Tm$, or $\xi$ belongs to $\geo f^{-1}(e_k)$
for all sufficiently large $k$, which is absurd.
\qed

\medskip
\begin{remark}
Let $G\acts X$, $T$, $\ga$, and $\eta$ be as in
Corollary \ref{geoetadetection}.  When $X$ is a $3$-dimensional
Hadamard manifold, then $\geo X\simeq S^2$ and there
is an alternate characterization of $\geo^\eta X$
which uses little more than the definition of itineraries.
For each $i$, $\geo X_{e_i}\subset \geo X\simeq S^2$
determines two closed disks by the Jordan separation
theorem; let $D_i$ be the one which contains $\geo X_{e_j}$
for all $j\geq i$.  Then $\geo^\eta X=\cap_i D_i$.  
To see this note that if $F\subset X$ is a flat totally
geodesic plane and $p\in X-F$, then the two components
of $\geo X-\geo F$ are $\{\xi\in\geo X\mid \ol{p\xi}\cap F\neq
\emptyset\}$ and $\{\xi\in\geo X\mid \ol{p\xi}\cap F=\emptyset\}$.
\end{remark}

\begin{corollary}
\label{preslabelling}
Let $G\acts X$ and $G\acts X'$ be admissible actions,  let
$G\acts T$ be the Bass-Serre tree of $G$, and let $V,\,E$ be the
sets of vertices and edges of $T$, respectively.
If $\Phi:\geo X\ra \geo X'$
is any $G$-equivariant homeomorphism, then

1.  $\Phi$ maps $\geo X_\si$ homeomorphically to $\geo X_\si '$ for all $\si\in V\cup E$.

2. $\Phi$ maps $\geo^\eta X$ homeomorphically to $\geo^\eta X'$ for all $\eta\in\geo T$.
\end{corollary}
\proof
Part 1 follows from the characterization of $\geo X_\si$ as a fixed point set
which is stated  in Lemma \ref{fixedinbdy}.   Part 2 follows from part 1 and
Corollary \ref{geoetadetection}.
\qed

\begin{corollary}
\label{structuretits}
Let $G\acts X$ be an admissible action, and $G\acts T$ be the 
Bass-Serre action for $G$.  Then

1.  The union $\cup_{v\in V}\,\tits X_v\subset\tits X$ is a $CAT(1)$ space with respect
to the induced metric, and may be described metrically as follows.
First, $\tits X_v$ is a metric suspension of an uncountable discrete $CAT(1)$
space for each $v\in V$, and $\tits X_e$ is isometric to the standard circle
for every $e\in E$.
Take the disjoint union $\amalg_{v\in V}\,\tits X_v$, and for each
edge $e=\ol{v_1v_2}\in E$, glue $\tits X_{v_1}$ to $\tits X_{v_2}$
isometrically by identifying the copies of $\tits X_e\subset \tits X_{v_i}$;
the result is  isometric to $\cup_{v\in V}\,\tits X_v\subset \tits X$.

2.  The union $\cup_{v\in V}\,\tits X_v$ forms a connected component of $\tits X$.
The remaining components are contained in the subsets $\tits^\eta X$
for $\eta\in \geo T$.  We will show in Lemma \ref{interval}  that each $\tits^\eta X$
is either a point or isometric to an interval of length $<\pi$.
\end{corollary}
\proof  To prove 1, we first observe that if $e_1,\ldots,e_n$ 
is an edge path in $T$ with initial vertex $v_1$ and terminal
vertex $v_n$, $\xi_1\in\tits X_{v_1}$, $\xi_n\in\tits X_{v_n}$,
and $\tangle(\xi_1,\xi_n)<\pi$, then the Tits segment
$\ol{\xi_1\xi_n}\subset \tits X$  is contained in 
$\cup_{i=1}^n \,\tits X_{v_i}$.  To see this, pick $\xi
\in\ol{\xi_1\xi_n}$.  Recall that for a given base point $x$, $\ol{x\xi}$ may be obtained
as the limit of a sequence $\ol{xy_k}$ where $y_k$ lies
on a segment $\ol{x_1^kx_2^k}$ and $x_i^k\in\ol{x\xi_i}$ 
is a sequence tending to infinity.   The quasi-convexity property
of Lemma \ref{templateqconvexity} (applied in sucession to $\ol{x\xi_i}$, $\ol{x_1^kx_2^k}$, and then $\ol{xy_k}$) implies that $\ol{xy_k}\subset N_C(\cup_{i=1}^n
X_{v_i})$ for some $C$, and the convexity of the $N_C(X_{v_i})$'s
implies that $\xi\in\tits X_{v_i}$ for some $i\in\{1,\ldots n\}$.
Part 1 now follows from Corollary \ref{titscaptits}.

Before proving 2, we recall that open  balls of radius $\frac{\pi}{2}$ in $CAT(1)$ spaces
are geodesically convex, so two points in a $CAT(1)$ space belong
to the same connected component iff they can be joined by a unit
speed path.

Suppose $\eta\in\geo T$ and $\xi\in\tits^\eta X$.   Fix $v\in V$, $p\in X_v$,
and let $e_k\in E$ denote the $k^{th}$ edge of the ray $\ol{v\eta}\subset T$.
Part 1 of Lemma \ref{itindichotomy} implies that $\ol{p\xi}\cap X_{e_k}\neq\emptyset$
for all but finitely many $k$.  If $c:[0,L]\ra\tits X$ is a unit speed path starting
at $\xi$, then by Lemma \ref{bdyconvex}, we find that either $c([0,L])\subset\tits^\eta X$
or for all sufficiently large $k$ there is a $t_k\in [0,L]$ so that 
$\ol{pc(t_k)}\in\tits X_{e_k}$.   But part 1 shows that when $e,\,e'\in E$
and $d(e,e')\geq 2$, then  $d(\tits X_e,\tits X_{e'})$ is bounded away from $0$,
which gives a contradiction.   Therefore $c([0,L])\subset\tits^\eta X$ and
we have shown that the connected component of $\xi$ is contained in 
$\tits^\eta X$. 

Finally, we note that $\cup_{v\in V}\,\tits X_v$ is connected: if $v_1,\ldots,v_k$
are the consecutive vertices of a geodesic segment in $T$, then $\tits X_{v_i}\cap\tits X_{v_{i+1}}
=\tits X_{\ol{v_iv_{i+1}}}\neq\emptyset$.
\qed

\subsection{Geometric data and equivariant quasi-isometries}
\label{samedata}

Throughout this section, $G\stackrel{\rho}{\acts} X$ and 
$G\stackrel{\rho'}{\acts} X'$
will denote  admissible actions
of an admissible group $G$ on Hadamard spaces $X$ and $X'$, and
we let $MLS_v,\,\tau_v$ and $MLS_v',\,\tau_v'$ denote
their respective geometric data (see Definition \ref{geometricdatadef}).  The main result in this section
is Theorem \ref{equivdatasublinear}, which shows
 that a $G$-equivariant quasi-isometry
$\Phi:X\ra X'$  induces
an equivariant homeomorphism $\geo\Phi:\geo X\ra\geo X'$
provided $\rho$ and $\rho'$ have {\em equivalent} geometric data.

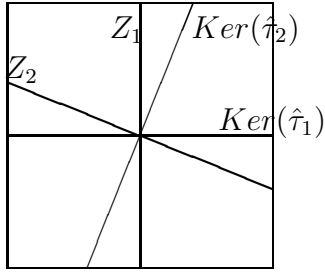
\begin{figure}
\begin{center}

\begin{picture}(100,100)
\put(0,0){\framebox(100,100)}
\thicklines

\put(50,0){\line(0,1){100}}
\put(0,70){\line(5,-2){100}}
\put(45,90){\makebox(0,0){$Z_1$}}
\put(5,75){\makebox(0,0){$Z_2$}}
\thinlines
\put(0,50){\line(1,0){100}}
\put(30,0){\line(2,5){40}}
\put(100,55){\makebox(0,0){$Ker(\hat \tau_1)$}}
\put(90,90){\makebox(0,0){$Ker(\hat \tau_2)$}}
\end{picture}

\bigskip

\caption{On $G_e\otimes\R\simeq \R^2$ we have two metrics $<\cdot,\cdot>$ and $<\cdot,\cdot>'$.  $Ker(\hat \tau_i)$ (which is the same as $Ker(\hat \tau'_i)$) for $i\in\{1,2\}$ is perpendicular to $Z_i$ with respect to both metrics.}
\label{twoscales}

\end{center}
\end{figure}

\begin{definition}
\label{geomdataequiv}
  We say
that $\rho$ and $\rho'$ have {\em equivalent} geometric data if 
there are functions $\la:V\ra\R_+$ and $\mu:V\ra\R_+$ 
so that for every $v\in V$
$$ MLS_v'=\la(v)MLS_v\mbox{\quad and \quad}\tau_v'=\mu(v)\tau_v.$$  
It follows from the $G$-invariance of the geometric data that
$\la$ and $\mu$ will be $G$-invariant.
\end{definition}
The structure of $G$ strongly restricts the possibilities for the 
functions $\la$ and $\mu$:

\begin{lemma}
\label{constants}
Suppose the geometric data for the actions
$G\acts X$ and  $G\acts X'$ are equivalent,  
and let $\la:V\ra\R_+$
 and $\mu:V\ra\R_+$ be as in Definition \ref{geomdataequiv}.
 Then either
 
 1. $\la\equiv\mu\equiv a$, for some $a>0$.
 
 or
 
 2. There are constants $a$ and $b$, $a\not= b$, so that 
$\la(V)=\{a,b\}=\mu(V)$.  Moreover,
for any pair $v_1,\,v_2$
 of adjacent vertices in $T$, $\la(v_1)=\mu(v_2)$, $\la(v_2)=\mu(v_1)$,
 and the
 $\R$ directions of $Y_{v_i}\simeq \bar Y_{v_i}\times\R$ determine
 orthogonal directions  in $Y_e$, where $e\defeq\ol{v_1v_2}$. In particular, there
 is a $G$-equivariant $2$-coloring of $V$ (i.e. a two coloring of the finite graph ${\cal G}=T/G$) such that $\la$ and $\mu$ are functions
 of the vertex color. 
\end{lemma}
\proof
Pick $e=\ol{v_1v_2}\in E$, and consider $G_e\otimes\R\simeq
\Z^2\otimes\R\simeq \R^2$.  For $i=1,2$ we have subspaces
$Z_i\defeq (Z(G_{v_i})\cap G_e)\otimes\R\simeq \Z\otimes\R\simeq\R$
determined by the centers of the $G_{v_i}$'s.  The action $G_e\acts Y_e$ 
induces an inner product $\langle\cdot,\cdot\rangle$ on $G_e\otimes\R$ 
by letting, for $g\in G_e$,  $\langle g,g \rangle=\delta_g^2$.  
(For $p\in Y_e$ there is an embedded Euclidean plane $\R^2\subset Y_e$ 
invariant under the action of $G_e\acts Y_e$ on which $G_e$ acts 
by translations of $\delta_g$.  So our metric is naturally related to 
the metric on this $\R^2$).  We can extend the maps $\tau_i:G_e\to R$ 
to linear maps $\hat\tau_i:G_e\otimes\R\to R$.  For $g\in G_e$ 
we see that $\tau_i(g)$ is just the change in the $i$-vertical 
component from $y$ to $g(y)$.  So for $x\in G_e\otimes\R$, $\hat \tau_i(x)$ 
is just the length of the orthogonal projection of $x$ onto $Z_i$. 
In particular, for $x\in Z_i$ we have $\langle x,x\rangle=\hat \tau_i(x)^2$, 
and $Ker(\hat \tau_i)$ is $\langle\cdot,\cdot\rangle$ perpendicular to $Z_i$.  
Similarly   $\hat {MLS}_i(x)$ corresponds to the length of the 
projection of $x$ perpendicular to $Z_i$ (i.e. to $Ker(\tau_i)$).
In general $Z_1$ and $Z_2$ are not perpendicular but they are always linearly independent (see figure \ref{twoscales}.

Similarly, using the action $G_e\acts Y_e'$, we get an 
induced inner product $\langle\cdot,\cdot\rangle'$ and linear 
maps $\hat \tau_1'$ and $\hat \tau_2'$.  By assumption $\hat\tau_i'=
\mu(v_i)\hat\tau_i$, hence $Ker(\hat\tau_i)
=Ker(\hat\tau_i')$.  So the space perpendicular to $Z_1$ (resp. $Z_2$) 
with respect to $\langle\cdot,\cdot\rangle$ is the same as that with 
respect to $\langle\cdot,\cdot\rangle'$.  Hence, by the 
independence of $Z_1$ and $Z_2$, either $\langle\cdot,\cdot\rangle'=a\langle\cdot,\cdot\rangle$ 
for some $a>0$ or $Z_1$ and $Z_2$ are perpendicular in both metrics so $Ker(\hat\tau_1)=Ker(\hat\tau_1')=Z_2$ and 
$Ker(\hat\tau_2)=Ker(\hat\tau_2')=Z_1$.
In the latter case, choosing $x\in Z_2$ we have 
$\langle x, x\rangle'=(\hat {MLS}_1(x))^2=\mu(v_2)^2(\hat {MLS}_1(x))^2=\mu(v_2)^2\langle x,x\rangle$ and
$\langle x,x\rangle'=
(\tau_2'(x))^2=\la(v_1)^2(\tau_2(x))^2=
\la(v_1)^2\langle x,x\rangle$ so
$\la(v_1)=\mu(v_2)$; similarly $\la(v_2)=\mu(v_1)$.
\qed

\medskip
We now fix a $G$-equivariant $(L,A)$-quasi-isometry $\Phi:X\ra X'$
for the rest of this section.  The constants defined will depend on the geometry of $G\stackrel{\rho}{\acts} X$ and 
$G\stackrel{\rho'}{\acts} X'$ as well as any other explicitly stated quantities.

\begin{lemma}
\label{sigmaspacespreserved}
There is a constant $D_1=D_1(L,A)$ so that for every 
$\si\in V\cup E$, the Hausdorff distance
$d_H(\Phi(X_\si),X_\si')$
is at most $D_1$.
\end{lemma}
\proof For $g\in G$ let $d_g$ and $d_g'$ denote the displacement
functions for $g$ in $X$ and $X'$ respectively.  In particular $d_g'(x)\leq L d_g(x) +A$.

For any $v\in V$,
if $g\in Z(G_v)$ is a generator of the center $Z(G_v)$ then 
(see section \ref{groupsonx}) $d_g$ is proper on $X/Z(Z(G_v),G)=X/G_v$ (by Lemma \ref{properdisp}) and hence 
$d_g$ (resp $d_g'$) grows with the distance from 
$X_v$ (resp $X_v'$).   This (along with the fact that there are only finitely many $\sigma$ modulo $G$) implies the lemma when
$\si\in V$.   If $e\in E$, and
$g_1,\,g_2\in G_e$ are a basis for $G_e$ then 
$\max\{d_{g_1},d_{g_2}\}$ (resp $\max\{d_{g_1}',d_{g_2}'\}$)
grows with the distance from $X_e$ (resp $X_e'$). 
This implies the lemma when $\si\in E$.
\qed

\medskip
We now assume for the remainder of this section that
$\rho$ and $\rho'$ have equivalent geometric data,
and we let  $\la:V\ra\R_+$
 and $\mu:V\ra\R_+$ be as in Definition \ref{geomdataequiv}.
For each $v\in V$, we have  nearest point
projections $p_v:X\ra Y_v$ and $p_v':X'\ra Y_v'$.
Modulo renormalization of
the metrics on the spaces $Y_v$, the equivariant quasi-isometry
$\Phi$ restricts to a  Hausdorff approximation:

\begin{lemma}
\label{closetoprod}
For every $v\in V$, the $G_v$-equivariant
quasi-isometry $\Phi_v\defeq p_v'\circ\Phi\restr_{Y_v}:Y_v\ra Y_v'$
has the following properties:

1. It is at distance $<D_2=D_2(L,A)$ from a $G_v$-equivariant map 
$\Psi_v:Y_v\ra Y_v'$ which respects 
the product structures $Y_v\simeq \bar Y_v\times\R$ and 
$Y_v'\simeq \bar Y_v'\times\R$.

2. If we stretch the metric on the $\bar Y_v$ factor of $Y_v$ by
$\la(v)$, and the metric on the $\R$ factor by $\mu(v)$, then
$\Phi_v$ becomes a $D_3=D_3(L,A)$-Hausdorff
approximation, and maps unit speed geodesic segments to within
$D_3$ of unit speed geodesic segments.
\end{lemma}
\proof
We first prove part 1.  
The $\R$ fibers of $Y_v$ and $Y_v'$ are within uniform Hausdorff
distance of $Z(G_v)$-orbits, so  $G_v$-equivariance
implies that $\Phi_v$ takes $\R$ fibers of $Y_v$ to within uniform
Hausdorff distance (say $C_1$) of $\R$ fibers of $Y_v'$.  The $\bar Y_v$
fibers of $Y_v$ are within uniform Hausdorff distance of
sets of the form $\{ g(p)\mid \mbox{$p\in Y_v$,$g\in G_v$,  $|\tau_v(g)|<C$}\}$
for sufficiently large $C$, and a similar characterization
of the $\bar Y_v'$ fibers of $Y_v'$ holds.  We now define the product map 
$\Psi_v=\bar\Psi_v\times\Psi_v^\R:\bar Y_v\times\R\ra
\bar Y_v'\times\R$.  Fix a basepoint $p\in Y_v$.  We may assume that the 
$\R$ factor of $p$ and $\Phi_v(p)$ are $0$ and  take $\Psi_v^\R(t)=\mu(v)t$.  
We let $S\subset \bar Y_v$ be a (set theoretical) cross section for the $H_v$ 
action.  For $s\in S$ choose a $g\in G_v$ such that $|\tau_v(g)|<C$ and 
$d(g(p),(s,0))\leq C_1$.  Now, since $|\tau'_v(g)|<\mu(v)C$, we can 
choose a point $\bar\Psi_v(s)$ such that $d(g(\Phi(p)),(\bar\Psi_v(s),0))\leq \mu(v)C$.  
We note that $$d(\Psi((s,0)),\Phi((s,0)))\leq  d((\bar\Psi_v(s),0),\Phi(g(p))) + LC_1+A\leq \mu(v)C + LC_1+A$$
Extend this to an $H_v$ equivariant map $\bar\Psi_v:\bar Y_v\to \bar Y_v'$.  
Thus, along with the fact that $\tau'(g)=\mu(v)\tau(g)$, we see that 
$\Psi_v$ is a $G_v$-equivariant map.  Now for every $q\in Y_v$ 
there is a $g\in G_v$ such that $g(q)=(s,t)$ for some $s\in S$ and some $-C<t<C$.  Now $d(\Psi(q),\Phi(q))=$
$$d(\Psi((s,t)),\Phi((s,t)))\leq 2C\mu(v)+d(\Psi((s,0)),\Phi((s,0)))\leq 3\mu(v)C + LC_1+A$$
This proves 1.

Part 2 follows from part 1 if we can show that $\bar\Psi_v$
and $\Psi_v^\R$ carry unit speed geodesics to within uniform distance
of unit speed geodesics.  The map $\Psi_v^\R$
clearly does, because the translation distance in the $\R$-direction
is measured by $\tau_v$ (resp. $\tau_v'$) and these have
ratio $\mu(v)$.   Recall that $H_v\defeq G_v/Z(G_v)$ acts discretely
and cocompactly on the hyperbolic metric spaces $\bar Y_v$
and $\bar Y_v'$. 
Since $MLS_v'=\la(v)MLS_v$, if we renormalize the metric
on $\bar Y_v$ by $\la(v)$  we can apply Lemma \ref{preservinggeodesics}
to see that $\bar\Psi_v:\bar Y_v\ra \bar Y_v'$ 
preserves unit speed geodesics up to uniform error.  This proves
2. 
\qed

\medskip
We may now use our quasi-isometry $\Phi:X\ra X'$ to transport standard $K$-templates
for $X$ to templates for $X'$ (see section \ref{templatesection}).  
Start with a standard $K$-template $(\Tm,f,\phi)$ for 
some segment or ray $\ga\subset T$.  Recall that the walls
of $\Tm$  come from flats $F_e\subset Y_e\subset X_e$.
To produce the new template $\Tm'$  distort the metric
on $\Tm$ by an affine change as follows.
For each $v\in V$, we think of scaling the
metric on $\bar Y_v$ by $\la(v)$ and the $\R$-factor
of $Y_v$ by $\mu(v)$; and then we distort the flats 
$F_e=F_{\ol{vv'}}\subset F_v\cap F_{v'}$ and strips
$\St_{ee'}\subset Y_{e\cap e'}$ used to build $\Tm$
accordingly.   Lemmas  \ref{sigmaspacespreserved} and
\ref{closetoprod}
imply that $(\Tm',f,\Phi\circ\phi)$ is
a $K'$-template for $\ga$ where $K'=K'(L,A)$.  Notice (using Lemma \ref{constants}) that the
identity map $\Tm\ra\Tm'$ is an affine map (i.e. maps
constant speed geodesics to constant speed geodesics),
and is a homothety when $\la=\mu=a\in\R$.  

\begin{theorem}
\label{equivdatasublinear}  Let $G\stackrel{\rho}{\acts} X$ and $G\stackrel{\rho'}{\acts} X'$ be admissible actions of an admissible group G on Hadamard spaces $X$ and $X'$ such that
$\rho$ and $\rho'$ have equivalent geometric data, and let $\Phi:X\to X'$ be a $G$-equivariant $(L,A)$-quasi-isometry. 
Then there is a function $\th:\R_+\ra\R_+$ (depending on K,L, A and  the geometry of X and X') with $\lim_{r\ra\infty}\th(r)=0$
so that for every $x,y\in X$, $z\in\ol{xy}$, we have 
\begin{equation}
\label{shadowtwice}
d_{X'}(\Phi(z),\ol{\Phi(x)\Phi(y)})\leq
(1+d_X(z,x))\th(d_X(z,x)).
\end{equation}  Consequently, by Lemma \ref{sublinearbending},
$\Phi$ extends to a unique
map $\bar\Phi:\bar X\ra\bar X'$ which is continuous at
$\geo X\subset\bar X$.  Setting $\geo\Phi\defeq\bar\Phi\restr_{\geo X}:
\geo X\ra \geo X'$,
 we obtain a $G$-equivariant homeomorphism.  If $\la=\mu=a\in\R$, then
 $\geo\Phi$ is  an isometry with respect to Tits metrics. 
\end{theorem}
\proof
Pick $x,y\in X$, and
then find a segment $\ga\subset T$ so that $x\in X_{e_1}$, 
$y\in X_{e_n}$,
and the $i^{th}$ edge of $\ga$ is $e_i$.  
Now let $(\Tm,f,\phi)$ be the standard $K$-template for $\ga$,
and define the $K'$-template $(\Tm',f,\phi')$ as in the 
paragraph preceding the statement of the Theorem.  We may
assume, after moving $x$ and $y$ a uniformly bounded distance if 
necessary (see part 3 of definition \ref{ktemplate}),  that $x=\phi(x_1)$, $y=\phi(y_1)$ for
some $x_1,\,y_1\in\Tm$.  We get (\ref{shadowtwice}) by applying 
Theorem \ref{shadthm} twice -- once to $(\Tm,f,\phi)$ and once to 
$(\Tm',f,\phi')$.  Specifically, Since $z\in \ol{\phi(x_1)\phi(y_1)}$ an application of Theorem \ref{shadthm} part 2 gives
$$d_X(z,\phi(\ol{x_1y_1}))\leq (1+d_X(z,x))\theta_{(X,K)}(d_X(z,x)).$$
So there is a $z_1\in \ol{x_1y_1}$ with 
$d_X(z,\phi (z_1))\leq (1+d_X(z,x))\theta_{(X,K)}(d_X(z,x)).$  
Hence we see
$$d_{X'}(\Phi(z),\phi'(z_1))\leq L(1+d_X(z,x))\theta_{(X,K)}(d_X(z,x))+A.$$
Now an application of Theorem \ref{shadthm} part 1 to $\phi'$ gives
$$d_{X'}(\phi'(z_1),\ol{\Phi(x)\Phi(y)})=d_{X'}(\phi'(z_1),
\ol{\phi'(x_1)\phi'(y_1)})\leq (1+d_{\Tm'}(x_1,z_1))\theta_{(X',K')}(d_{\Tm'}(x_1,z_1)).$$
We thus need to bound $d_{\Tm'}(x_1,z_1)$ linearly from above by $d_X(z,x)$.  
But this follows since  $d_{\Tm'}(x_1,z_1)\leq L'd_{X'}(\phi'(x_1),\phi'(z_1))+A'$ 
(where $L'$ and $A'$ come from Proposition \ref{ambienttotemplatepath} and 
depend only on $K'$ and the geometry of $X'$),  
$d_{X'}(\phi'(x_1),\phi'(z_1))\leq Ld_X(\phi(x_1),\phi(z_1))+A$, and 
$d_X(\phi(x_1),\phi(z_1))\leq d_X(x,z)+d_X(z,\phi(z_1))\leq d_X(x,z)+(1+d_X(z,x))\theta_{(X,K)}(d_X(z,x))$. 
Equation (\ref{shadowtwice}) will thus follow for an appropriate 
choice $\theta_{(K,L,A,X,X')}$ that will depend only on $L$, $A$, $K$, and the geometry of $X$ and $X'$.  
Lemma \ref{sublinearbending} 
then applies to $\Phi$, so we get an induced embedding $\geo\Phi:\geo X
\ra \geo X'$ which is automatically $G$-equivariant.  

Now suppose $\la=\mu=a$.  Pick $\eta\in\geo T$, and  a geodesic ray
$\ga\subset T$ with $\geo \ga=\eta$.   Then the identity map $\Tm\ra \Tm'$
between the associated templates $(\Tm,f,\phi)$ and $(\Tm',f,\phi')$
constructed as above is a homothety, and so part 5 of 
Theorem \ref{shadthm}  applied to $\geo\phi\restr_{\geo^\infty\Tm}$
and $\geo\phi'\restr_{\geo^\infty\Tm'}$ then shows that $\geo\Phi\restr_{\geo^\eta X}$
induces an isometry $\tits^\eta X\ra \tits^\eta X'$.   By Corollary
\ref{structuretits} it remains only to show that $\geo\Phi$ induces an isometry
$\tits X_v\ra\tits X_v'$.  But Lemma \ref{closetoprod} part 2 implies
that $\Phi_v:Y_v\ra Y_v'$ is at finite distance from a product of
Hausdorff approximations (up to rescaling of $X'$ by $\frac{1}{a}$)
and so $\Phi_v$ induces an isometry $\tits X_v\ra \tits X_v'$.
\qed

\subsection{Recovering the geometric data from the action on the ideal
boundary}
\label{recovering}

In this section we will prove the remaining implication of Theorem
\ref{main}: for any admissible action $G\acts X$, the topological conjugacy class of the
action $G\acts\geo X$ 
determines the functions 
$MLS_v:G_v\ra \R_+$ 
and $\tau_v:G_v\ra\R$ up to a multiplicative factor, for every
vertex  $v\in{\cal G}$.   Our strategy for proving this is as follows.
Using Lemma \ref{fixedinbdy} and
Corollary \ref{geoetadetection}, for any ideal boundary point $\eta$
of the Bass-Serre tree $T$, we may detect the subset
$\geo^\eta X\subset\geo X$; specifically, the action $G\acts \geo X$
determines the set of boundary points $\eta\in\geo T$ for which
$|\geo^\eta X|=1$.  To extract useful
information from this, we consider a special class of geodesic
rays $\ga\subset T$ (see Definition \ref{specialrays}) 
which admit templates $(\Tm,f,\phi)$ where
$\Tm$ is asymptotically self-similar: there is a (non-surjective)
map $\Tm\ra\Tm$ which stretches distances by a factor of $2$.
These templates have two key properties: their geometry
relates directly to the geometric data $MLS_v$ and $\tau_v$,
and at the same time we can tell explicitly when (in terms 
of the geometry) we have $|\geo^{\geo\ga}\Tm|=1$
(see section \ref{selfsimsec}).  Putting 
all this together we able to recover the geometric data 
from the action $G\acts \geo X$.

\subsubsection{$\tits^\infty\Tm$ is either a point or an interval}
\label{pointorinterval}

\begin{lemma}
\label{embedsin0pi}
Let $X$ be a Hadamard space, $p\in X$, and $\ga_i\subset X$ a
sequence of geodesics with $\lim_{i\ra\infty}d(p,\ga_i)=\infty$.
Set
$$S\defeq\{ \xi\in\tits X\mid \mbox{$\ol{p\xi}\cap \ga_i\neq\emptyset$
for all $i$}\}.$$
Then $S$ embeds isometrically in the interval $[0,\pi]$.
\end{lemma}

\proof
Pick $x,\,y,\,z\in S$, and for each $i$ choose $x_i\in
\ol{px}\cap\ga_i$,
$y_i\in \ol{py}\cap\ga_i$, $z_i\in \ol{pz}\cap\ga_i$.  After
passing to a subsequence and reordering $x,\,y,\,z$ we may assume
that $y_i$ lies between $x_i$ and $z_i$ on $\ga_i$, or that it
coincides with $x_i$ or $z_i$.  Then it follows that
$$\cangle_p(x_i,y_i)+\cangle_p(y_i,z_i)\leq \cangle_p(x_i,z_i).$$
Taking the limit as $i\ra\infty$, we get
$$\tangle(x,y)+\tangle(y,z)\leq \tangle(x,z).$$
Hence $(\{x,y,z\},\tangle)\subset\tits X$ embeds isometrically in
$[0,\pi]$.

If $|S|=1$ the lemma is immediate, so assume $|S|\geq 2$, and
construct a map $f:S\ra\R$ as follows.  Pick distinct points
$x_0,\,x_1\in S$ and for $i=0,\,1$ choose $f(x_i)\in\R$
so that $d(x_0,x_1)=d(f(x_0),f(x_1))$.  Now define $f$ uniquely
by the condition that $d(f(s),f(x_i))=d(s,x_i)$ for all $s\in S$
and $i=0,\,1$.  Clearly $f$ is an isometric embedding, and its
image lies in an interval of length at most $\pi$.
\qed

\begin{lemma}
\label{leqbeta}
Let $\Tm$ be a half template with walls $W_0,W_1,\ldots$.  For $i\geq 1$
set $\al_i\defeq \al(W_i)$ where $\al:Wall_\Tm^0\ra (0,\pi)$
is the angle function, and assume
$\min\{\alpha_i,\pi-\alpha_i\}\geq\beta$ for all $i\geq 1$.
Then

1. The diameter of $\tits^\infty\Tm$ with respect to the
Tits metric is at most $\pi-\beta$.

2. For $p\in W_0$ and every $i>0$ there is an $R_i$, depending only on $\beta$ and $d(o_i,p)$, so that if 
$q\in L_i^+$ then $\ol{pq}\cap L_i^-\subset B(o_i,R_i)$.
\end{lemma}

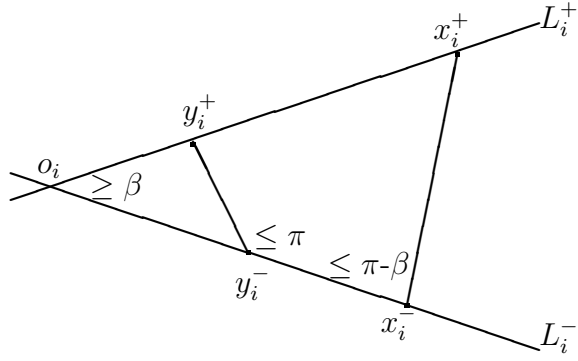
\begin{figure}
\begin{center}

\begin{picture}(100,100)
\thicklines
\put(-50,55){\line(3,-1){200}}
\put(-50,45){\line(3,1){200}}
\put(-40,55){\makebox{$o_i$}}
\put(40,25){\line(-1,2){21}}
\put(40,25){\framebox(0,0)}
\put(35,10){\makebox{$y_i^-$}}
\put(42,28){\makebox{$\leq\pi$}}
\put(100,5){\framebox(0,0)}
\put(90,-5){\makebox{$x_i^-$}}
\put(70,17){\makebox{$\leq\pi$-$\beta$}}
\put(100,5){\line(1,5){19}}
\put(19,66){\framebox(0,0)}
\put(15,75){\makebox{$y_i^+$}}
\put(119,100){\framebox(0,0)}
\put(110,106){\makebox{$x_i^+$}}
\put(150,-10){\makebox{$L_i^-$}}
\put(150,110){\makebox{$L_i^+$}}
\put(-20,47){\makebox{$\geq\beta$}}

\end{picture}

\bigskip

\caption{On the wall $W_i$ the angle between $L_i^-$ and $L_i^+$ is $\geq \beta$ hence the angle at $x_i^-$ is $\leq \pi-\beta$.  The angle at $y_i^-$ is always $\leq \pi$ (no mater where $y_i^+$ lies on $L_i^+$).}
\label{anglei}

\end{center}
\end{figure}

\proof
Pick $p\in W_0\subset\Tm$, and distinct points
$x,\,y\in\tits^\infty\Tm$.
Let $x_i^-\in L_i^-$ (resp. $x_i^+\in L_i^+$) be the point
where the ray $\ol{px}$ enters (resp. exits) the  wall
$W_i$; define $y_i^-,\,y_i^+$ similarly.   Since $x\neq y$,
there is an $i_0>0$ so that $x_i^-\neq y_i^-$
when $i\geq i_0$.  Pick $i\geq i_0$, and assume that
$d(x_i^-,o_i)\geq d(y_i^-,o_i)$ (the case when
$d(y_i^-,o_i)\geq d(x_i^-,o_i)$ is similar).  Clearly (see figure \ref{anglei})
we have
$$\angle_{x_i^-}(x_i^+,o_i)=\angle_{x_i^-}(x_i^+,y_i^-)\leq
\max\{\alpha_i,\pi-\alpha_i\}\leq \pi-\beta.$$
Now this implies 1 since by standard properties of Tits angles
$$\tangle(x,y)=\left
[\lim_{i\ra\infty}(\angle_{x_i^-}(y_i^-,x)+\angle_{y_i^-}(x_i^-,y))\right
]-\pi
\leq \max\{\alpha_i,\pi-\alpha_i\}\leq \pi-\beta.$$
It also implies 2 by applying triangle comparison to the angle at
$x_i$ of the triangle with vertices $p$, $o_i$, and $x_i$.
\qed

\begin{proposition}
\label{interval}
Let $\Tm$ be a uniform half template with
$\min\{\alpha_i,\pi-\alpha_i\}\geq
\beta$ for
all $i$.  Then
$\tits^\infty\Tm$ is isometric to an interval $[0,\theta]$
where $\theta\in [0,\beta]$.

\end{proposition}

\proof
Applying Lemma \ref{embedsin0pi} with $\ga_i=L_i^-$ and Lemma
\ref{leqbeta} we get that $\tits^\infty\Tm$ is isometric to
a subset of $[0,\beta]$.  We need only show that $\tits^\infty\Tm$
is connected.  Suppose $x,\,y\in \tits^\infty\Tm$, and
pick a point $z$ lying on the Tits interval
$\ol{xy}\subset\tits\Tm$.  Let $p, \,x_i^-,\,y_i^-$
be as in the proof of Lemma \ref{leqbeta}.  Then
$\lim_{i\ra\infty}\cangle_p(x_i^-,y_i^-)=\tangle(x,y)$,
so we can choose a sequence $z_i\in\ol{x_i^-y_i^-}\subset L_i^-$
so that $\lim_{i\ra\infty}\cangle_p(x_i^-,z_i)=\tangle(x,z)$
and $\lim_{i\ra\infty}\cangle_p(z_i^-,y_i)=\tangle(z,y)$.
The segments $\ol{pz_i}$ converge to the ray $\ol{pz}$.
Since for any $j>0$ the segment
$\ol{pz_i}$ crosses $L_j^+$ for sufficiently large $i$,
by Lemma \ref{leqbeta} part 2 we get
$$\ol{pz_i}\cap L_j^-\subset B(o_j,R_j)$$
and we conclude that $\ol{pz}\cap L_j^-\neq \emptyset$.
Hence $z\in \tits^\infty\Tm$.
\qed

The $\theta$ in the above proposition is referred to as the Tits angle
of $\Tm$ and is denoted $\theta(\Tm)$.

\subsubsection{Self-similar Templates}
\label{selfsimsec}

In this section we study a special class of full templates called {\em
self-similar templates}.  
\begin{definition}
\label{selfsimdef}
Let $\Tm$ be a full template with $Wall_\Tm=\{W_i\}_{i\in\Z}$ and 
$Strip_\Tm=\{\St_i\}_{i\in\Z}$, and set $\al_i\defeq\al(W_i)$,
$l_i\defeq l(\St_i)$, and $\eps_i\defeq \eps(\St_i)$ (we define
$\eps$ using the strip orientation compatible with the strip directions and the usual ordering
on $\Z$).  Then $\Tm$ 
is a {\em self-similar template} if for all $i,\,j\in\Z$ we have
$\al_i=\al_j$, $l_{i+2j}=2^jl_i$, and $\eps_{i+2j}=2^j\eps_i$.
In this case we say that $\Tm$ has data $(\beta;l_0,\eps_0,l_1,\eps_1)$
where $\beta=\al_i$ for all $i\in\Z$.
\end{definition}
Note that by the definition, if we rescale the metric on a self-similar
template $\Tm$ by a factor of $2$, then we get a template equivalent
to $\Tm$, i.e. there is a homothety  $\Phi: \Tm \to \Tm$ which preserves
strip directions, stretches distances by a factor of $2$, and which
shifts wall and strip indices by $2$.   A self-similar template
is determined up to equivalence by $\beta$ and the data
$\{l_0,\eps_0,l_1,\eps_1\}$.

$\Tm$ contains one point $v$ which is not on any wall or strip (since
the union of the walls and strips is not complete).  The point $v$ is the limit of
the Cauchy sequence $\{o_{-i}|i=0,1,2,...\}$.  On the other hand, for
each $j\in \Z$ the half template, $\Tm_j$, given by the union of planes
$W_i$ and strips $\St_i$ for $i\geq j$ is uniform (and complete).
The images of the embeddings $\geo^\infty\Tm_i\ra\geo\Tm$
and $\tits^\infty\Tm_i\ra\tits \Tm$ are independent of $i$,
and we use $\geo^\infty\Tm$ (respectively $\tits^\infty\Tm$)
to denote this common subspace.   We will
say that $\Tm$ is {\em trivial} if $|\geo^\infty\Tm|=1$ and 
{\em nontrivial} if $|\geo^\infty\Tm|>1$.

We note that if $r$ is a geodesic ray parameterized by arclength then
$\Phi \circ r$ is a geodesic ray parameterized by twice arclength.  So
$\Phi$ and $\Phi^{-1}$ take rays to rays.  Since $\Phi$ is a homothety
it preserves the Tits angles between rays.

If ${\cal R}_i$ represents the space of geodesic rays starting at $o_i$ and
intersecting $W_j$ for all $j\geq i$, then the above shows $\Phi^{-1}
({\cal R}_i)={\cal R}_{i-2}$ and that $\Phi^{-1}$ acts as a Tits isometry on
$\tits^\infty \Tm$.  In particular since $\tits^\infty \Tm$ is isometric
to an interval, $\Phi$ leaves the midpoint fixed, and $\Phi^2$ acts as
the identity.  Thus by repeated applications of $\Phi^{-1}$ we see that
we can represent $\tits^\infty \Tm$ as the set rays, ${\cal R}_{-\infty}$, that
start at $v$ and are invariant under $\Phi^2$.  Further, the middle ray
is invariant under $\Phi$.

We will show in the lemma below that all rays in ${\cal R}_{-\infty}$ are
$\Phi$ invariant.  In particular any such ray that intersects an $o_i$
must intersect all $o_{i+2n}$ for $n\in \Z$, and hence there are at most
two such rays, $r_{even}$ and $r_{odd}$, in ${\cal R}_{-\infty}$.

\begin{figure}
\begin{center}

\begin{picture}(400,200)
\put(0,0){\framebox(400,200)}
\thicklines
\put(0,100){\framebox(0,0)}

\put(24,106){\line(1,-4){26}}
\put(24,106){\line(-1,-4){24}}
\put(23,110){\makebox(0,0){$o_0$}}
\put(24,55){\makebox(0,0){\large I}}
\thinlines
\put(23,102){\line(1,0){2}}
\put(22,98){\line(1,0){4}}
\put(21,94){\line(1,0){6}}
\put(20,90){\line(1,0){8}}
\put(19,86){\line(1,0){10}}
\put(18,82){\line(1,0){12}}
\put(17,78){\line(1,0){14}}
\put(16,74){\line(1,0){16}}
\put(15,70){\line(1,0){18}}
\put(14,66){\line(1,0){20}}
\put(13,62){\line(1,0){22}}
\put(12,58){\line(1,0){24}}
\put(11,54){\line(1,0){26}}
\put(10,50){\line(1,0){28}}
\put(9,46){\line(1,0){30}}
\put(8,42){\line(1,0){32}}
\put(7,38){\line(1,0){34}}
\put(6,34){\line(1,0){36}}
\put(5,30){\line(1,0){38}}
\put(4,26){\line(1,0){40}}
\put(3,22){\line(1,0){42}}
\put(2,18){\line(1,0){44}}
\put(1,14){\line(1,0){46}}
\put(0,10){\line(1,0){48}}
\put(0,6){\line(1,0){49}}
\put(0,2){\line(1,0){50}}
\thicklines

\put(48,112){\line(1,-4){28}}
\put(48,112){\line(-1,-4){28}}
\put(48,116){\makebox(0,0){$o_2$}}
\put(48,58){\makebox(0,0){\large I}}
\thinlines
\put(46,104){\line(1,0){4}}
\put(44,96){\line(1,0){8}}
\put(42,88){\line(1,0){12}}
\put(40,80){\line(1,0){16}}
\put(38,72){\line(1,0){20}}
\put(36,64){\line(1,0){24}}
\put(34,56){\line(1,0){28}}
\put(32,48){\line(1,0){32}}
\put(30,40){\line(1,0){36}}
\put(28,32){\line(1,0){40}}
\put(26,24){\line(1,0){44}}
\put(24,16){\line(1,0){48}}
\put(22,8){\line(1,0){52}}
\thicklines

\put(96,124){\line(1,-4){31}}
\put(96,124){\line(-1,-4){31}}
\put(96,128){\makebox(0,0){$o_4$}}
\put(96,64){\makebox(0,0){\large I}}
\thinlines
\put(92,108){\line(1,0){8}}
\put(88,92){\line(1,0){16}}
\put(84,76){\line(1,0){24}}
\put(80,60){\line(1,0){32}}
\put(76,44){\line(1,0){40}}
\put(72,26){\line(1,0){48}}
\put(68,10){\line(1,0){56}}
\thicklines

\put(192,148){\line(1,-4){37}}
\put(192,148){\line(-1,-4){37}}
\put(192,152){\makebox(0,0){$o_6$}}
\put(192,76){\makebox(0,0){\large I}}
\thinlines
\put(184,116){\line(1,0){16}}
\put(176,84){\line(1,0){32}}
\put(168,52){\line(1,0){48}}
\put(160,20){\line(1,0){64}}
\thicklines

\put(384,196){\line(1,-4){15}}
\put(384,196){\line(-1,-4){49}}
\put(384,98){\makebox(0,0){\large I}}
\thinlines
\put(368,132){\line(1,0){32}}
\put(352,68){\line(1,0){48}}
\put(336,4){\line(1,0){64}}
\thicklines


\put(32,92){\line(1,4){27}}
\put(32,92){\line(-1,4){27}}
\put(34,88){\makebox(0,0){$o_1$}}
\put(34,144){\makebox(0,0){\large III}}
\thinlines
\put(31,96){\line(1,0){2}}
\put(30,100){\line(1,0){4}}
\put(29,104){\line(1,0){6}}
\put(28,108){\line(1,0){8}}
\put(27,112){\line(1,0){10}}
\put(26,116){\line(1,0){12}}
\put(25,120){\line(1,0){14}}
\put(24,124){\line(1,0){16}}
\put(23,128){\line(1,0){18}}
\put(22,132){\line(1,0){20}}
\put(21,136){\line(1,0){22}}
\put(20,140){\line(1,0){24}}
\put(19,144){\line(1,0){26}}
\put(18,148){\line(1,0){28}}
\put(17,152){\line(1,0){30}}
\put(16,156){\line(1,0){32}}
\put(15,160){\line(1,0){34}}
\put(14,164){\line(1,0){36}}
\put(13,168){\line(1,0){38}}
\put(12,172){\line(1,0){40}}
\put(11,176){\line(1,0){42}}
\put(10,180){\line(1,0){44}}
\put(9,184){\line(1,0){46}}
\put(8,188){\line(1,0){48}}
\put(7,192){\line(1,0){50}}
\put(6,196){\line(1,0){52}}
\thicklines

\put(64,84){\line(1,4){29}}
\put(64,84){\line(-1,4){29}}
\put(64,80){\makebox(0,0){$o_3$}}
\put(64,140){\makebox(0,0){\large III}}
\thinlines
\put(62,92){\line(1,0){4}}
\put(60,100){\line(1,0){8}}
\put(58,108){\line(1,0){12}}
\put(56,116){\line(1,0){16}}
\put(54,124){\line(1,0){20}}
\put(52,132){\line(1,0){24}}
\put(50,140){\line(1,0){28}}
\put(48,148){\line(1,0){32}}
\put(46,156){\line(1,0){36}}
\put(44,164){\line(1,0){40}}
\put(42,172){\line(1,0){44}}
\put(40,180){\line(1,0){48}}
\put(38,188){\line(1,0){52}}
\put(36,196){\line(1,0){56}}
\thicklines

\put(128,68){\line(1,4){33}}
\put(128,68){\line(-1,4){33}}
\put(128,64){\makebox(0,0){$o_5$}}
\put(128,132){\makebox(0,0){\large III}}
\thinlines
\put(124,84){\line(1,0){8}}
\put(120,100){\line(1,0){16}}
\put(116,116){\line(1,0){24}}
\put(112,132){\line(1,0){32}}
\put(108,148){\line(1,0){40}}
\put(104,164){\line(1,0){48}}
\put(100,180){\line(1,0){56}}
\put(96,196){\line(1,0){64}}
\thicklines

\put(256,36){\line(1,4){41}}
\put(256,36){\line(-1,4){41}}
\put(256,32){\makebox(0,0){$o_7$}}
\put(256,116){\makebox(0,0){\large III}}
\thinlines
\put(248,68){\line(1,0){16}}
\put(240,100){\line(1,0){32}}
\put(232,132){\line(1,0){48}}
\put(224,164){\line(1,0){64}}
\put(216,196){\line(1,0){80}}

\put(197,200){\vector(1,-4){5}}
\put(212,140){\vector(1,-4){5}}
\put(227,80){\vector(1,-4){5}}
\put(242,20){\vector(1,-4){5}}
\put(218,100){\makebox{$\St_6$}}

\put(183,200){\vector(-1,-4){5}}
\put(168,140){\vector(-1,-4){5}}
\put(153,80){\vector(-1,-4){5}}
\put(138,20){\vector(-1,-4){5}}
\put(158,100){\makebox{$\St_5$}}

\put(0,100){\vector(4,1){400}}
\put(0,100){\vector(4,-1){400}}
\put(320,170){\makebox(0,0){$r_{even}$}}
\put(320,30){\makebox(0,0){$r_{odd}$}}
\end{picture}

\bigskip

\caption{Development of a self-similar template from quarter planes
$Q^I_i$. Developements of quadrants for nonnegative integers are shaded
and labeled with their type, while those with negative indices are not
shown.  The direction of $\St_5$ and $\St_6$ are shown via arrows.
The interval of rays between $r_{even}$ and $r_{odd}$  correspond to rays in the
template.}
\label{scale}

\end{center}
\end{figure}
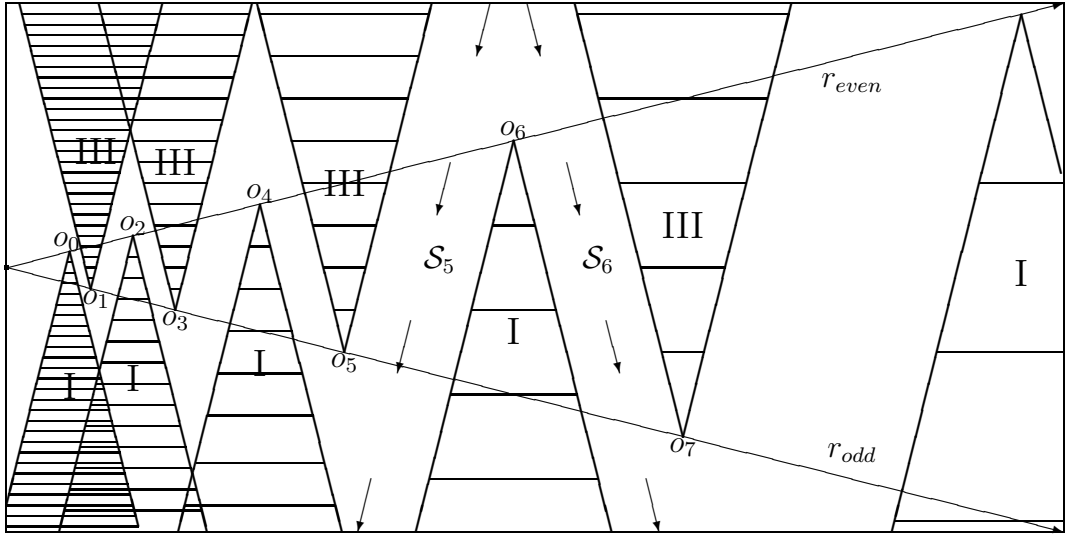

Each choice $N\in \{I,II,III,IV\}$ determines a choice of
quarter planes $Q^N_i\subset W_i$ as follows: for all
$k\in\Z$, $Q^N_{2k}=Q_N\subset W_{2k}$
while $Q^N_{2k+1}=-Q_N\subset W_{2k+1}$.  Straightforward Euclidean
geometry shows that the corresponding development map, ${\cal D}_N$,
with ${\cal D}_N(v)=0$ has the property that ${\cal D}_N\circ\Phi \circ
{\cal D}_N^{-1}$ is multiplication by 2 wherever (and however) it is
defined (see Figure \ref{scale}).  (The easiest way to see this is to first
develop $Q^N_0$, $\St_0$, $Q^N_1$, $\St_1$, and $Q^N_2$ into the plane,
then shift the origin $(0,0)$ so that it lies on the line through ${\cal
D}(o_0)$ and ${\cal D}(o_2)$ and such that ${\cal D}(o_0)$ lies between
$(0,0)$ and ${\cal D}(o_2)$ and such that the distance from $(0,0)$ to
${\cal D}(o_2)$ is twice the distance to ${\cal D}(o_0)$.  We note that
${\cal D}(Q^N_2)=2{\cal D}(Q^N_0)$.  We can now define a map ${\cal D}$
uniquely so that ${\cal D}\circ\Phi \circ {\cal D}^{-1}$ is
multiplication by 2.  It is easy to check that this map is a
development and hence by uniqueness is ${\cal D}_N$ - up to an element
of $O(2)$.)  There are rays $r_{even}$ and $r_{odd}$ from the origin
such that ${\cal D}_N(o_{2i})\in r_{even}$ and ${\cal D}_N(o_{2i+2})\in
r_{odd}$.  (We can fix ${\cal D}_N$ completely if desired by taking
$r_{even}$ to be the positive $x$-axis and to make $r_{odd}$ point in
the upper half plane.)

\begin{lemma}
\label{sitemplates}
Let $\Tm$ be a nontrivial self-similar template.  Then there is a
choice $N\in \{I,II,III,IV\}$ such that for $r\in {\cal R}_{-\infty}$, $r\cap
W_i\subset Q^N_i$.  Each $r\in {\cal R}_{-\infty}$ is $\Phi$ invariant, and
$r_{even}$ and $r_{odd}$ are the ${\cal D}_N$ images of the boundary
rays of ${\cal R}_{-\infty}$ as a closed interval.  In particular the set of
$r\in {\cal R}_{-\infty}$ span a space isometric, via ${\cal D}_N$ to a convex
cone in $R^2$.
\end{lemma}

\proof The fact that ${\cal R}_{-\infty}$ is a closed interval follows from
Proposition \ref{interval}.  If $r_1,r_2\in {\cal R}_{-\infty}$ then their
invariance under the homothety $\Phi^2$ implies that $\angle
_T(r_1,r_2)=\angle_v(r_1,r_2)$, and hence
the rays between $r_1$ and $r_2$ span a space isometric to a convex cone
in $R^2$.  In particular, no ray between $r_1$ and $r_2$ can intersect
an origin $o_i$.  Thus the middle ($\Phi$ invariant) ray $r_m$
determines $N,M\in\{I, II, III, IV\}$ such that $r_m\cap W_{2n}\subset
Q_N$, and $r_m\cap W_{2n+1} \subset Q_{M}$.  The above says that for any
interior ray, $r$, $r\cap W_{2n}\subset Q_N$, and $r\cap W_{2n+1}
\subset Q_{M}$.  Consider the development map ${\cal D}$ determined by
this choice of quarter planes.  We will assume that ${\cal D}(v)$ is the
origin, so in particular ${\cal D}(Q_i)$ and ${\cal D}(\St_i)$ miss the
origin.  Now ${\cal D}\circ\Phi^2 \circ {\cal D}^{-1}$ is defined when
restricted to the image of each strip or quarter plane.  It is clearly
just multiplication by 4, since that is what it does to each ray in the
$\cal D$ image of ${\cal R}_{-\infty}$.

We assume that $N=I$ (the other cases are similar) and must prove that
$M=III$.  If $M=II$ (resp. $M=IV$) then for $i=1$ (resp. $i=0$), ${\cal
D}(Q_i)\cup{\cal D}(\St_i)\cup {\cal D}(Q_{i+1})$ will contain a half
plane which in turn contains one of ${\cal D}(Q_i)$ or ${\cal
D}(Q_{i+1})$.  But this cannot happen since any ray entering that half
plane will never leave it.  If $M=I$ then for $k\geq max\{\frac {\pi}
{\beta},\frac {\pi} {\pi-\beta}\}$ we have $\cup_{i=1}^k ({\cal
D}(Q_i)\cup {\cal D}\St_i)$ contains a half plane which in turn contains
one of the quarter planes and the same argument works.  Thus we conclude
that $M=III$, and $Q_i=Q^N_i$.

We know that ${\cal D}_N\circ\Phi \circ {\cal D}_N^{-1}$ is just
multiplication by 2.  This means in particular that all rays in
${\cal R}_{-\infty}$ are invariant under $\Phi$.  Now since at least one ray
between $r_{even}$ and $r_{odd}$ hits every ${\cal D}_N(Q_i)$ this is
true of all such rays, hence all are ${\cal D}_N$ of a ray in $\Tm$.
This gives a 1-1 correspondence between rays in $\Tm$ and rays in the
plane between $r_{even}$ and $r_{odd}$ completing the lemma.

\qed

The argument in the last proof shows that a self-similar template is
nontrivial if and only if there is an $N \in \{I,II,III,IV\}$ such that
some (and hence any) ray between the corresponding $r_{even}$ and
$r_{odd}$ intersects every ${\cal D}_N(Q^N_i)$.  But 
by self-similarity this will be true
 if and only if the line segment from ${\cal
D}_N(o_0)$ to ${\cal D}_N(o_2)$ intersects ${\cal D}_N(Q^N_1)$.  This
leads to the following lemma:

We will use the notation $\R^4_0=\{(x_1,x_2,x_3,x_4)\subset
\R^4\mid x_1,\, x_3>0\}$.

Let $A_\beta$ be:
$$\{(x,y)\in (-\frac \pi 2,\frac \pi 2)\times (-\frac \pi 2,\frac \pi
2)\mid\mbox{$x+\beta \geq y\geq x-\beta$  and $\ -x+(\pi-\beta)\geq y \geq
-x-(\pi-\beta)$}\}.$$

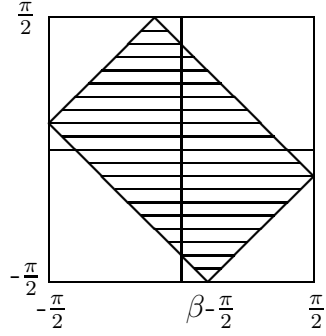
\begin{figure}
\begin{center}
\begin{picture}(110,110)

\put(5,5){
\begin{picture}(100,100)
\put(0,0){\framebox(100,100)}
\thicklines
\put(0,60){\line(1,1){40}}
\put(60,0){\line(1,1){40}}
\put(0,60){\line(1,-1){60}}
\put(40,100){\line(1,-1){60}}
\thinlines
\put(0,50){\line(1,0){100}}
\put(50,0){\line(0,1){100}}

\put(50,10){\line(1,0){20}}
\put(40,20){\line(1,0){40}}
\put(30,30){\line(1,0){60}}
\put(20,40){\line(1,0){80}}
\put(10,50){\line(1,0){80}}
\put(0,60){\line(1,0){80}}
\put(10,70){\line(1,0){60}}
\put(20,80){\line(1,0){40}}
\put(30,90){\line(1,0){20}}

\put(55,5){\line(1,0){10}}
\put(45,15){\line(1,0){30}}
\put(35,25){\line(1,0){50}}
\put(25,35){\line(1,0){70}}
\put(15,45){\line(1,0){80}}
\put(5,55){\line(1,0){80}}
\put(5,65){\line(1,0){70}}
\put(15,75){\line(1,0){50}}
\put(25,85){\line(1,0){30}}
\put(35,95){\line(1,0){10}}

\end{picture}
}
\put(10,0){\makebox(0,0)[t]{-$\frac \pi 2$}}
\put(110,0){\makebox(0,0)[t]{$\frac \pi 2$}}
\put(70,0){\makebox(0,0)[t]{$\beta$-$\frac \pi 2$}}
\put(0,105){\makebox(0,0){$\frac \pi 2$}}
\put(0,0){\makebox(0,0)[b]{-$\frac \pi 2$}}

\end{picture}

\bigskip

\caption{The graph of $A_\beta$ for $\beta=\frac {3\pi} 5$.}
\label{figureab}

\end{center}
\end{figure}

See Figure \ref{figureab}.

\begin{lemma}
\label{goodset}
Fix $\pi>\beta>0$ and let ${\cal A}_\beta$ be the set of
$(l_0,\eps_0,l_1,\eps_1)\subset \R^4_0$ such that the self-similar
template with  data $(\beta;l_0, \eps_0, l_1, \eps_1)$ is
trivial.  Then
$${\cal A}_\beta = \{(l_0,\eps_0,l_1,\eps_1)\in \R^4_0\mid (arctan(\frac {\eps_0}
{l_0}),arctan(\frac {\eps_1}{l_1}))\in A_\beta\}$$

\end{lemma}
\proof
Let $\psi_i=arctan(\frac {\eps_i}
{l_i})$.
The proof follows once we show that $R^4-{\cal A}_\beta$
consists of 4 components given in order for $N$=$I$, $II$, $III$, and
$IV$ by:
$$\psi_0>\psi_1+\beta,\ \ \  -\psi_1-(\pi-\beta)>\psi_0,\ \ \
\psi_1-\beta>\psi_0,\ \ \ {\rm and}\ \  \psi_0>-\psi_1+(\pi-\beta).$$

Here we do the case $N=I$, i.e. the line segment from ${\cal D}_I(o_0)$
to
${\cal D}_I(o_2)$ intersects ${\cal D}_I(Q^I_1)$ .  The cases $N=II,III,
{\rm and}\  IV$ are similar.

\begin{figure}
\begin{center}

\begin{picture}(100,100)
\put(0,0){\framebox(100,100)}
\thicklines
\put(40,0){\line(0,1){60}}
\put(15,30){\vector(0,1){70}}
\put(15,30){\vector(-1,2){15}}
\put(40,60){\line(1,-2){30}}
\put(80,50){\vector(-1,2){25}}
\put(80,50){\vector(0,1){50}}
\put(40,10){\vector(0,1){0}}
\put(65,9){\vector(-1,2){0}}

\put(8,80){\makebox(0,0){$I$}}
\put(72,80){\makebox(0,0){$I$}}
\put(50,20){\makebox(0,0){$III$}}

\thinlines

\put(15,90){\line(-1,0){15}}
\put(15,80){\line(-1,0){15}}
\put(15,70){\line(-1,0){15}}
\put(15,60){\line(-1,0){15}}
\put(15,50){\line(-1,0){10}}
\put(15,40){\line(-1,0){5}}

\put(40,10){\line(1,0){25}}
\put(40,20){\line(1,0){20}}
\put(40,30){\line(1,0){15}}
\put(40,40){\line(1,0){10}}
\put(40,50){\line(1,0){5}}

\put(80,60){\line(-1,0){5}}
\put(80,70){\line(-1,0){10}}
\put(80,80){\line(-1,0){15}}
\put(80,90){\line(-1,0){20}}

\put(40,60){\line(-5,-6){15}}
\put(40,60){\line(4,-1){15}}

\put(25,20){\vector(-1,0){9}}
\put(30,20){\vector(1,0){9}}
\put(70,45){\vector(2,1){10}}
\put(62,41){\vector(-2,-1){10}}
\put(28,20){\makebox(0,0){$l_0$}}
\put(66,43){\makebox(0,0){$l_1$}}

\put(20,47){\vector(0,1){13}}
\put(20,42){\vector(0,-1){13}}
\put(71,58){\vector(1,-2){5}}
\put(69,62){\vector(-1,2){5}}
\put(20,45){\makebox(0,0){$\epsilon_0$}}
\put(70,60){\makebox(0,0){$\epsilon_1$}}

\put(34,41){\makebox(0,0){$\theta_0$}}
\put(45,41){\makebox(0,0){$\theta_1$}}
\put(52,50){\makebox(0,0){$\theta_2$}}

\put(10,50){\makebox(0,0){$\beta$}}
\put(75,70){\makebox(0,0){$\beta$}}

\thinlines

\end{picture}

\caption{Developement in case I; $\epsilon_0>0$, $\epsilon_1<0$ .}
\label{figure8.6}

\end{center}
\end{figure}
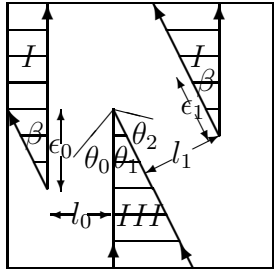

Let $\pi>\theta_0>0$ be the angle between the line segment from ${\cal
D}_I(o_1)$ to ${\cal D}_I(o_0)$ and the ``incoming'' edge of ${\cal
D}_I(Q^I_1)$ (see Figure \ref{figure8.6}).  Here the ``incoming'' edge is the ${\cal D}_I$ image of the
negative half line of $L_1^-$ (since $Q^I_1$ is of type $III$).
Similarly let $\theta_2$ be the angle between the line segment from
${\cal D}_I(o_1)$ to ${\cal D}_I(o_2)$ and the ``outgoing'' edge of
${\cal D}_I(Q^I_1)$.  $\theta_1$ will be the angle at ${\cal D}_I(o_1)$
of the sector ${\cal D}_I(Q^I_1)$ (which in our case is just $\beta$).
It is easy to see that our condition is equivalent to
$$\pi> \theta_0 + \theta_1 + \theta_2.$$
The rest of the argument is just that the definitions (being careful
about how the sign of $\eps_i$ and the orientation of the strips
interact)
give us, when $N=I$,
$$\theta_0=\frac \pi 2 - arctan(\frac{\eps_0} {l_0})\ \ \ {\rm and} \ \
\theta_2=\frac \pi 2 - arctan(\frac{-\eps_1} {l_1})$$
which completes the argument.
\qed

We extract the information we need with the following elementary (but
somewhat non-trivial):

\begin{corollary}
\label{geodata}
Pick $a_i>0$ for $i=1,\ldots,4$, $b_1>0$, $b_2\in\R$, and 
a subset ${\cal B}\subset \R^4_0$.  Suppose there is a $\be\in(0,\pi)$ so that
${\cal B}$ is precisely the set of 
$(x_1,x_2,x_3,x_4)\in \R^4_0$  for which
$$(a_1x_1+b_1,a_2x_2+b_2,a_3x_3,a_4x_4)\in {\cal A}_\beta.$$
Then  $\beta$, $\frac {b_1} {a_1}$, and $\frac {b_2} {a_2}$ are uniquely
determined by ${\cal B}$.   If this unique $\be$ is not $\frac \pi 2$ then
$\frac {b_1} {a_2}$ and $\frac {b_2} {a_1}$ are also determined by
${\cal B}$.
\end{corollary}

\proof We consider the map $$\Psi(x_1,x_2,x_3,x_4)=(arctan(\frac
{a_2x_2+b_2}{a_1x_1+b_1}),arctan(\frac {a_4x_4}{a_3x_3}))\eqdef
(\psi_0(x_1,x_2),\psi_1(x_3,x_4)).$$
the previous lemma along with our assumption says that for
$a_1x_1+b_1>0$ and $x_3>0$
$(x_1,x_2,x_3,x_4)\in {\cal B}\Longleftrightarrow
\Psi(x_1,x_2,x_3,x_4)\in A_\beta$.  $\Psi$ maps
$\{(x_1,x_2,x_3,x_4)|a_1x_1+b_1>0\ {\rm and}\ x_3>0\}$ onto $(-\frac \pi
2,\frac \pi 2)\times (-\frac \pi 2,\frac \pi 2)$.

{\em Step 1:  $\beta$ is determined by ${\cal B}$.}
We notice that $(x,y^0)\in A_\beta$ for all $\frac \pi 2 -\epsilon
<x<\frac \pi 2$ if and only if $y^0=\frac \pi 2-\beta$.  Similarly
$(x^1,y)\in A$ for all large $\frac \pi 2 -\epsilon <y<\frac \pi 2$ if
and only if $x^1=\frac \pi 2 -\beta$. Thus $(1,x_2,x^0_3,x^0_4)\in {\cal
B}$ for all sufficiently large $x_2$ if and only if
$\psi_1((x^0_3,x^0_4)) = \frac \pi 2-\beta$, while
$(x^1_1,x^1_2,1,x_4)\in {\cal B}$ for all sufficiently large $x_4$ if
and only if $\psi_0((x^1_1,x^1_2)) = \frac \pi 2-\beta$.  Fix such an
$(x^1_1,x^1_2)$, then, if $\beta \not= \frac \pi 2$, there is a unique
$x^1_4$ such that for all $x_4\geq x^1_4$, $(x^1_1,x^1_2,1,x_4)\in
{\cal B}$.  If $\beta <\frac \pi 2$ we have $\psi_1(1,x^1_4)=\frac \pi 2
-2\beta$, while if $\beta>\frac \pi 2$ we have $\psi_1(1,x^1_4)=
2\beta-\frac {3\pi} 2$.

Let $\bar {a_i}$, $\bar {b_1},\bar {b_2}$, and $\bar \beta$ be another
choice of parameters that work and $\bar \Psi$ and $\bar \psi_i$ the
corresponding functions.  Then $\tan(\bar \psi_1) = c\tan(\psi_1)$ where
$c=\frac {\bar a_4a_3}{\bar a_3 a_4}>0$.  Plugging in $(x^0_3,x^0_4)$
from the previous paragraph allows us to conclude
$\tan(\frac \pi 2 - \bar \beta) = c \tan(\frac \pi 2 - \beta)$,
(i.e. $\tan(\bar \beta)=\frac 1 c \tan(\beta)$ if $\beta\not= \frac \pi
2$).  Note in particular that the sign of $\frac \pi 2 -\bar \beta$ is
the same as the sign of $\frac \pi 2 -\beta$.  Plugging in $(1,x^1_4)$
from the above paragraph we get $\tan(\frac \pi 2-2\bar
\beta)=c\tan(\frac \pi 2-2\beta)$ if $\beta <\frac \pi 2$ and
$\tan(2\bar \beta-\frac {3\pi} 2)=c\tan(2\beta-\frac {3\pi} 2)$ if
$\beta>\frac \pi 2$.  However this can only happen if c=1 and $\beta =
\bar \beta$.  Thus we conclude that $\beta$ is determined by ${\cal B}$
and that either $\beta=\frac \pi 2$ or else $c=1$ and $\psi_1=\bar\psi_1$.

{\em Step 2: If $\beta \not= \frac \pi 2$ then $\psi_0=\bar \psi_0$, while if
$\beta = \frac \pi 2$ then $\tan (\bar \psi_0)=\frac 1 c \tan(\psi_0)$.}
We will first show that given ${\cal B}$, $\beta$ and $\psi_1$ then
there is at most one choice for $\psi_0$.  To see this fix $(x_1,x_2)$
and consider $S\subset R^2$ such that ${\cal B}\cap ((x_1,x_2)\times
R^2)=(x_1,x_2)\times S$.  Then $\psi_1(S)$ is precisely the interval
such that $A_\beta\cap (\psi_0(x_1,x_2)\times R)=(\psi_0(x_1,x_2)\times
\psi_1(S))$.  However the shape of $A_\beta$ is such that $\psi_1(S)$
thus determines $\psi_0(x_1,x_2)$, unless $\beta=\frac \pi 2$ in which
case it determines $\psi_0(x_1,x_2)$ up to sign.  The sign is determined
by continuity and the fact that for $x_2$ large and positive (resp.
negative) then $\psi_0(x_1,x_2)$ is positive (resp. negative).

Thus if $\beta \not= \frac \pi 2$ then $\psi_1=\bar\psi_1$ and we see that $\bar \psi_0 =\psi_0$ is
the unique solution.  If $\beta = \frac \pi 2$ then $\tan(\bar \psi_1) = c\tan(\psi_1)$ and $A_\beta =\{(x,y)\in
R^2|-x+\frac \pi 2\geq y \geq x-\frac \pi 2$ when $x>0$ and $x+\frac \pi
2\geq y \geq -x-\frac \pi 2$ when $x<0\}$.  Thus for $c>0$ the map
$(x,y)\rightarrow (arctan(\frac 1 c\tan(x)),arctan(c\tan(y))$ preserves
$A_\beta$ and step 2 follows.

The rest of the proof follows from the following equation that holds for
all positive $x_1$ and all $x_2$:
$$\frac {\bar a_2x_2+\bar b_2}{\bar a_1x_1+\bar b_1}=\frac 1 c \frac
{a_2x_2+b_2}{a_1x_1+b_1}$$
and the fact that $c=1$ when $\beta \not= \frac \pi 2$.

\qed

\subsubsection{Recovering the data}

In this section ${\cal G}$ will denote a fixed admissible graph of groups,
$G\defeq \pi_1({\cal G})$  the fundamental group of ${\cal G}$, and 
$G\acts T$ the Bass-Serre action for ${\cal G}$.    For every vertex $v$
of ${\cal G}$ we choose a generator $\zeta_v\in Z(G_v)$ for the center of
$Z(G_v)$ as in Definition \ref{geometricdatadef}.     Also, we will fix
an admissible action $G\acts X$.  We will use the template notation 
from section \ref{subsectemplates}.  Recall that
when $\Tm$ is a half-template, then $\geo^\infty\Tm\subset \geo\Tm$
denotes the set of boundary points corresponding to rays which
intersect all but finitely many walls.

The goal of this section is to prove the remaining half of Theorem \ref{main}:
the topological conjugacy class of the
action $G\acts\geo X$ determines the functions 
$MLS_v:G_v\ra \R_+$ 
and $\tau_v:G_v\ra\R$ up to a multiplicative factor, for every
vertex  $v\in{\cal G}$.

\begin{definition}
An element $g$ of a vertex group $G_v$ is {\em restricted}
if $g$ acts on $\bar Y_v$ (see section \ref{sectionvertexandedge})
as an axial isometry and its fixed points  in $\geo \bar Y_v$
are distinct from 
the fixed points of  $G_e$ where $e$ is any edge incident to $v$.
\end{definition}
It will be convenient to choose, for each vertex $v$ of ${\cal G}$,
a restricted  element $\de_v\in G_v$ which lies in the commutator subgroup
$[G_v,G_v]\subset G_v$:
\begin{lemma}
\label{lotsofrestricted}
For every vertex $v$, the commutator subgroup $[G_v,G_v]$ contains
restricted elements.
\end{lemma}
\proof We first recall that $H_v\defeq G_v/Z(G_v)$ is a nonelementary hyperbolic
group, and the induced action $H_v\acts \bar Y_v$ is discrete and
cocompact.       

Choose  a free nonabelian subgroup $S\subset H_v$ \cite[p. 212]{hypgps},
and elements $\bar g_1,\,\bar g_2\in S$ which belong
to a free basis for the commutator subgroup $[S,S]\subset [H_v,H_v]$,
and let $g_i\in G_v$ be  lift of $\bar g_i$ under the projection $G_v\ra H_v$.
Then $g_i$ acts axially on $\bar Y_v$ since $\bar g_i$ has infinite order
in $H_v$  and $H_v$ acts discretely on $\bar Y_v$.  Note that
$Fix(g_1,\geo\bar Y_v)\cap Fix(g_2,\geo\bar Y_v)=\emptyset$, since
otherwise by Lemma \ref{discreteaxis} we would have $Fix(g_1,\geo\bar Y_v)=Fix(g_2,\geo\bar Y_v)$,
forcing $\<g_1,g_2\>$ to be virtually cyclic, which is absurd.
Set $h_n\defeq g_1^ng_2^n$.  Lemma \ref{almostfree} tells us that
$h_n$ is axial for large $n$ and $Fix(h_n,\geo\bar Y_v)$ converges
to $\{\xi_1,\xi_2\}\subset \geo\bar Y_v$ where $\xi_i\in Fix(g_i,\geo\bar Y_v)$.  The induced action of $G_e$ on $\bar Y_v$ translates a geodesic $\ga_e\subset \bar Y_v$.  By the finiteness of $\cal G$ we can choose elements $g_e\in G_e$ such that the induced translation of $g_e$ is nonzero but uniformly bounded.
Thus Lemma \ref{discreteaxis} says that subsets $Fix(G_e,\geo\bar Y_v)=Fix(g_e,\geo\bar Y_v)$ define a discrete subset
of $(\geo\bar Y_v\times\geo \bar Y_v)/\Z_2$, so either

a) $Fix(h_n,\geo\bar Y_v)\cap Fix(G_e,\geo\bar Y_v)=\emptyset$ for all edges
$e$ incident to $v$ when $n$ is large,

or

b) There is a subsequence $h_{n_i}$ and an edge $e$ incident to $v$ 
so that $Fix(h_{n_i},\geo\bar Y_v)= Fix(G_e,\geo\bar Y_v)$.  

But if b) held then
we would have $ Fix(G_e,\geo\bar Y_v)=\{\xi_1,\xi_2\}$, which, by Lemma \ref{discreteaxis},
would force the absurd conclusion that $Fix(g_1,\geo\bar Y_v)=Fix(g_e,\geo\bar Y_v)=Fix(g_2,\geo \bar Y_v)$.
Hence case a) holds and the lemma is proved.
\qed

Notice that for every $v$,  $MLS_v(\de_v)\neq 0$ (since $\de_v$
acts on $\bar Y_v$ as an axial isometry), $\tau_v(\de_v)=0$ since
$\de_v\in[G_v,G_v]$,  and
$\tau_v(\zeta_v)\neq 0$.

\begin{lemma}
\label{needthis}
In order to determine $MLS_v$ and $\tau_v$ up to a multiplicative
factor, it suffices to determine the ratios
$$\frac{MLS_{v}(\si)}{MLS_{v}(\de_{v})}\quad\mbox{and}\quad \frac{\tau_{v}(\si)}
{\tau_{v}(\zeta_{v})}$$
for every restricted element $\si\in G_v$ whose fixed point set in 
$\geo\bar Y_v$ is disjoint from $Fix(\de_v,\geo\bar Y_v)$.
\end{lemma}
\proof Choose an arbitrary $\si\in G_v$.  

We first discuss $MLS_v$.
Note that $MLS_v(\si)=0$
iff $\si$ projects to an element of finite order in $H_v$; hence
we may assume that $\si$ acts on $\bar Y_v$ as an axial isometry.
First assume that $Fix(\si,\geo \bar Y_v)=Fix(\de_v,\geo \bar Y_v)$.
Let $\bar\si,\bar \de_v\in H_v$ be the projections to $H_v$.
Then $\bar\si$ and $\bar \de_v$ generate a virtually cyclic
subgroup $S$ because they have a common axis.   Hence there
is a finite subset $\{s_1,\ldots,s_k\}\subset G_v$ so that for
any $n$ we have $\si^n=s_{i_n}\de_v^{j_n}\zeta_v^{k_n}$ for suitable $i_n,\,j_n,\,k_n$.
Then 
$$MLS_v(\si)=\lim_{n\ra\infty}\frac{1}{n}MLS_v(\si^n)=\lim_{n\ra \infty}
MLS_v(\de_v^{j_n})$$
so we can recover the ratio above in this case.

Now assume $Fix(\si,\geo \bar Y_v)\cap Fix(\de_v,\geo \bar Y_v)=\emptyset$.
Setting $h_k\defeq \de_v^k\si^{k^2}$, we argue as in the proof of
Lemma \ref{lotsofrestricted} to see that for large $k$, $h_k$ is restricted,
$Fix(h_k,\geo \bar Y_v)\cap Fix(\de_v,\geo \bar Y_v)=\emptyset$,
and $|MLS_v(h_k)-kMLS_v(\de_v)-k^2MLS_v(\si)|$ is uniformly bounded (by Lemma \ref{almostfree}).
We may then recover the desired ratios from the formula
$$MLS_v(\si)=\lim_{k\ra\infty}\frac{1}{k^2}MLS_v(h_k).$$

We now consider the behavior of $\tau_v$.  Suppose $\si$ projects to an 
element of finite order in $H_v$, or $Fix(\si,\geo \bar Y_v)\cap Fix(\de_v,\geo\bar Y_v)\neq
\emptyset$.   In either case we have a finite set $\{s_1,\ldots,s_k\}\subset G_v$
so that $\si^n=s_{i_n}\de_v^{j_n}\zeta_v^{k_n}$ for suitable $i_n,\,j_n,\,k_n$.
Then 
$$
\tau_v(\si)=\lim_{n\ra\infty}\frac{1}{n}\tau_v(\si^n)=\lim_{n\ra\infty}
[\frac{j_n}{n}\tau_v(\de_v)+\frac{k_n}{n}\tau_v(\zeta_v)].$$
For the case $Fix(\si,\geo \bar Y_v)\cap Fix(\de_v,\geo \bar Y_v)=\emptyset$, use the same $h_k$ as above for large $k$ since $\tau_v(h_k)=k^2\tau(\sigma)$. 

\qed

We now focus our attention on a vertex $\bar v_1$ of ${\cal G}$.
Choose an edge $\bar e$ of ${\cal G}$ incident to $\bar v_1$, and
lift $\bar e$ to an edge  $e$ of the Bass-Serre tree.   We adopt the
notation from section \ref{basicbassserre}   for the associated
graph of groups ${\cal G}'$, $G'\defeq \pi_1({\cal G}')$, $T'\subset T$,
etc.    We fix some restricted element $\si\in G_{\bar v_1}$ with
$Fix(\si,\geo\bar Y_{v_1})\cap Fix(\de_{v_1},\geo\bar Y_{v_1})= \emptyset$.

\begin{definition} (Special rays)
\label{specialrays}
For every $(p,q,r,s)\in\R^4_0$  we define a geodesic ray $\ga\subset T' $
and sequences $\hat l_i=l_i(p,q,r,s)$, $\hat \eps_i=\hat\eps_i(p,q,r,s)$ 
as follows.

{\em Case 1: $\bar e\subset{\cal G}$ is embedded.}  Let $e=\ol{v_1v_2}$. For each $i\in\N$
we set $s_{2i-1}\defeq \si^{2^{i-1}}\de_{v_1}^{[p2^{i-1}]}\zeta_{v_1}^{[q2^{i-1}]}$
and $s_{2i}\defeq\de_{v_2}^{[r2^{i-1}]}\zeta_{v_2}^{[s2^{i-1}]}$, where
$[x]$ denote the integer part of $x\in\R$.  Then we
let $\ga\subset T$ be the geodesic ray with successive edges
\begin{equation}
\label{amalgray}
e,\,s_1e,\,s_1s_2e,\ldots,s_1\ldots s_ke,\ldots,
\end{equation}
and define, for all $i\in\N$,
$\hat l_{2i-1}=2^{i-1}(MLS_{v_1}(\si)+|p|MLS_{v_1}(\de_{v_1}))$,
$\hat l_{2i}=2^{i-1}|r|MLS_{v_2}(\de_{v_2})$; and 
$\hat \eps_{2i+1}\defeq 2^i(\tau_{v_1}(\si)+q\tau_{v_1}(\zeta_{v_1}))$,
and $\hat\eps_{2i}\defeq 2^{i-1}s\tau_{v_2}(\zeta_{v_2})$. 

{\em Case 2: $\bar e\subset {\cal G}$ is a loop.}  Let $t$ be as in section \ref{basicbassserre}.  For each $i\in \N$ 
we set $s_{2i-1}\defeq \si^{2^{i-1}}\de_{v_1}^{[p2^{i-1}]}\zeta_{v_1}^{[q2^{i-1}]}$
and $s_{2i}\defeq\de_{v_1}^{[r2^{i-1}]}\zeta_{v_1}^{[s2^{i-1}]}$;
then we
let $\ga\subset T$ be the geodesic with successive vertices
\begin{equation}
\label{hnnray}v_1,\,tv_1,\,ts_1t^{-1}v_1,\,ts_1t^{-1}s_2tv_1,\ldots 
ts_1t^{-1}s_2t\ldots s_{2k}tv_1,\ldots,
\end{equation}
and define, for all $i\in\N$,
$\hat l_{2i-1}=2^{i-1}(MLS_{v_1}(\si)+|p|MLS_{v_1}(\de_{v_1}))$,
$\hat l_{2i}=2^{i-1}|r|MLS_{v_1}(\de_{v_1})$; and 
$\hat \eps_{2i+1}\defeq 2^i(\tau_{v_1}(\si)+q\tau_{v_1}(\zeta_{v_1}))$,
and $\hat\eps_{2i}\defeq 2^{i-1}s\tau_{v_1}(\zeta_{v_1})$. 
\end{definition}
 These rays are useful because they admit templates that
are asymptotically self-similar:

\begin{lemma}
\label{itinscaleinvariant}
There is a half template  $(\Tm,f,\phi)$ for $\ga$ with walls
$Wall_\Tm=\{ W_i\}_{i=1}^\infty$ and strips
$Strip_\Tm=\{\St_i\}_{i=1}^\infty$, so that for $i$ sufficiently
large, 
$l(\St_i)=\hat l_i$, $\eps(\St_i)=\hat\eps_i$, and $\al(W_i)=\beta$,
where $\beta$ is the Tits angle between the (positively oriented)
$\R$-factors of $Y_{v_1}$ and $Y_{v_2}$.
\end{lemma}
\proof
We will treat the case when $\bar e$ is embedded; the other case
is similar.  To prove the lemma we will show that the desired template
may be obtained from a standard template for $\ga$ by changing
the strip widths and strip gluings by a bounded amount.  Thus by "uniformly" we will just mean independent of $i$, but possibly dependent on all other choices.

For $i\geq 1$ we set $e_i\defeq s_1\ldots s_{i-1}e$ and 
$v_i\defeq e_i\cap e_{i+1}$.   Recall (section \ref{templatesection})
that the standard template $(\Tm,f,\phi)$ for $\ga$ is constructed using flats
$F_{e_i}\subset Y_{e_i}$ and flat strips $\St_{e_i,e_{i+1}}\subset
Y_{v_i}$, where 
$\St_{e_i,e_{i+1}}=\ga_{e_i,e_{i+1}}\times\R\subset\bar Y_{v_i}\times\R
=Y_{v_i}$.   We choose $x_1\in F_{e_1}$ and set
$x_i\defeq s_1\ldots s_{i-1}x_1\in F_{e_i}$.

\medskip
{\em Step 1: There is a constant $c_1$ so that  $d(x_i,\St_{e_i,e_{i+1}})<c_1$, 
$d(x_{i+1},\St_{e_i,e_{i+1}})<c_1$ and  $|Width(\St_{e_i,e_{i+1}})-\hat l_i|<c_1$.}
We will do the case when $i=2j-1$ is odd;  the even case is similar.
Let $\pi:Y_{v_1}=\bar Y_{v_1}\times\R\ra\bar Y_{v_1}$ be the projection map, 
and set $\bar x_1\defeq\pi(x_1)$, $\bar F_e\defeq\pi(F_e)$, $\bar F_{s_ie}\defeq
\pi(\bar s_iF_e)=\pi(F_{s_ie})$, and $\ga_{e,s_ie}\defeq\pi(\St_{e,s_ie})$.
We apply $(s_1\ldots s_{i-1})^{-1}$ and then $\pi$ to everything, and 
 are thereby reduced to showing that there is a $c_1$ so that
$d(\bar x_1,\ga_{e,s_ie})<c_1$,  $d(s_i\bar x_1,\ga_{e,s_ie})<c_1$, and 
$|d(s_i\bar x_1,\bar x_1)-\hat l_i|<c_1$. 

>From the definition of $s_i$ we have $s_i\bar x_1=\si^{2^{j-1}}\de_{v_1}^{[p2^{j-1}]}\bar x_1$
since $\zeta_{v_1}$ acts trivially on $\bar Y_{v_1}$.    Since $\si,\,\de_{v_1}
\in G_{v_1}$ are restricted elements, the sets $\geo\bar F_e\subset\geo \bar Y_{v_1}$
and $\geo \bar F_{s_ie}\subset\geo\bar Y_{v_1}$ are disjoint from the 
fixed point sets of $\si$ and $\de_{v_1}$; the latter two sets are disjoint
by assumption.  Therefore we may apply Lemma \ref{almostfree} to 
conclude that when $j$ is sufficiently large,

a) $s_i:\bar Y_{v_1}\ra \bar Y_{v_1}$ is an axial isometry with 
an axis  $\ga_i\subset \bar Y_{v_1}$ at uniformly bounded distance from $\bar x_1$
and $s_i\bar x_1$.

b) $|d(s_i\bar x_1,\bar x_1)-\hat l_i|$
is uniformly bounded.

c) The attracting (resp. repelling) fixed point of $s_i$ in $\geo\bar Y_{v_1}$
is close to the attracting fixed point, $\xi^+$, of $\si$
(resp. repelling fixed point, $\xi^-$, of $\de_{v_1}^{sign(p)}$).

\smallskip
\no
Let $\ga_{e,s_ie}$ be the shortest path from 
$\bar F_e$ to $\bar F_{s_ie}$ with endpoints $\bar z_i\in \bar F_e$ and $\bar w_i \in \bar F_{s_ie}$
Let $\ga$ be a geodesic with endpoints $\xi^-$ and $\xi^+$. The Gromov hyperbolicity of $\bar Y_{v_1}$ and Lemma \ref{geodesicunion} part 3 imply that for $t$ fixed and large and $i$ large the axis of $s_i$ comes within a uniform distance of $\ga_{e,s_ie}(t)$ and hence $\ga_{e,s_ie}(t)$ stay a uniform distance from $\ga$.  Thus $\bar z_i$, which is the point on $\bar F_e$ closest to $\ga_{e,s_ie}(t)$, must stay a uniform distance from the set of points on $\bar F_e$ closest (in a Buseman function sense) to $\xi^+$ (which we see by taking $t$ large) and hence stay uniformly close to $\bar x_1$ .
Similarly $s_i^{-1}(\bar w_i)$ approaches the set of points on $\bar F_e$ closest to $\xi^-$ and hence $\bar w_i$ stays uniformly close to $s_i(\bar x_1)$.

{\em Step 2: There is a $c_2$ so that the standard template $(\Tm,f,\phi)$ satisfies
$|\eps(\St_i)-\hat\eps_i|<c_2$ and $\al(W_i)=\beta$ for all $i>1$.}
The assertion that $\al(W_i)=\beta$ is clear from the 
definition of $\beta(W_i)$ and the construction of standard
templates.   Step 1 then implies that there is a $c_3$ so that 
the origin $o_i\in W_i$ maps under $f:\Tm\ra X$ to within
distance $c_3$ of $x_i$, for $i>1$.  Hence from the definition
of $\eps:Strip^o_\Tm\ra\R$ we see that $\eps(\St_i)$ agrees with
$\tau_{v_i}((s_1\ldots s_{i-1})s_i(s_1\ldots s_{i-1})^{-1})$ to
within $2c_3$.  The conjugacy invariance of $\tau_v$
then gives $|\eps(\St_i)-\hat\eps_i|<c_2$ for a suitable $c_2$.

{\em Step 3: Adjusting $(\Tm,f,\phi)$.} In steps 1 and 2 we have
shown that the standard template satisfies conditions of Lemma
\ref{itinscaleinvariant} to within bounded error.  So we now
modify the construction of $\Tm$ by changing the metric\footnote{We
do this in the simplest way: we start with the metric product decomposition
$\hat\St_{e_i,e_{i+1}}\simeq I\times\R$ and then scale the metric 
on the $I$ factor.} on $\hat \St_{e_i,e_{i+1}}$so that $Width(\hat\St_{e_i,e_{i+1}})=\hat l_i$
if $l_i\geq 1$, and leaving $\hat\St_{e_i,e_{i+1}}$ untouched otherwise,
and by modifying the gluings $\partial\hat\St_{e_i,e_{i+1}}\ra W_{e_i}\,\amalg\,
W_{e_{i+1}}$ by a bounded amount so that $\eps(\St_i)=\hat\eps_i$.
Finally, if we redefine $f:\Tm\ra X$ to agree with the original 
$f$ on $\amalg_e \,W_e$ and on $\Tm-(\amalg_e \,W_e)$, then we get
the desired template for $\ga$.\qed

\medskip
\no
{\em Proof the Theorem \ref{main} concluded.}
Consider the subset ${\cal B}\subset\R^4_0$ of $4$-tuples $(p,q,r,s)$
for which the geodesic ray $\ga\subset T$ defined above gives a trivial subset
$\geo^{\geo\ga}X$ (i.e. a single point); by Lemma \ref{fixedinbdy} and
Corollary \ref{geoetadetection} the subset $\geo^{\geo\ga}X$
can be detected just using the action $G\acts\geo X$, and so
${\cal B}$ is also determined by the action $G\acts\geo X$.  On the
other hand, by Theorem \ref{shadthm}, $\geo^{\geo\ga}X$ is
trivial iff $\geo^\infty\Tm$ is trivial (i.e. a single point) where
$(\Tm,f,\phi)$ is any template for $\ga$.   Using Lemma \ref{itinscaleinvariant}
we arrive at the following:
the subset ${\cal B}$ of $(p,q,r,s)\in\R^4_0$
so that any template $\Tm$ with $l(\St_i)=\hat l_i(p,q,r,s)$, $\eps(\St_i)=\hat\eps_i(p,q,r,s)$,
and $\al(W_i)=\beta$ (for $i$ sufficiently large) is trivial,
is determined by the action $G\acts\geo X$.
But since a template $\Tm$ with $l(\St_i)=\hat l_i(p,q,r,s)$, $\eps(\St_i)=\hat\eps_i(p,q,r,s)$,
and $\al(W_i)=\beta$ (for $i$ sufficiently large) is trivial
iff the self-similar template with data $\{\beta;\hat l_3,\hat\eps_3,\hat l_4,\hat\eps_4\}$
is trivial, we may apply Corollary \ref{geodata} to conclude that 
the ratios 
$$\frac{MLS_{v_1}(\si)}{MLS_{v_1}(\de_{v_1})}\quad\mbox{and}\quad
\frac{\tau_{v_1}(\si)}{\tau_{v_1}(\zeta_{v_1})}$$
as well as $\beta$ are determined by the action $G\acts\geo X$.
Moreover, unless $\beta=\frac{\pi}{2}$ then 
$$\frac{MLS_{v_1}(\si)}{\tau_{v_1}(\zeta_{v_1})}\quad\mbox{and}\quad
\frac{\tau_{v_1}(\si)}{MLS_{v_1}(\de_{v_1})}$$
are also determined.

\subsection{Examples}
\label{examplesection}
In this section we construct the example mentioned in the introduction: we
will find two locally compact Hadamard spaces $X_0$ and $X_r$ on which an
admissible group $G$ acts discretely and cocompactly with the same geometric
data (i.e. the induced actions of $G$ on $\geo X_0$ and $\geo X_r$ are
topologically conjugate), an equivariant quasi-isometry $\tilde F:X_0\to X_r$,
and a geodesic ray $\gamma$ of $X_0$ such that $\tilde F(\gamma)$ does not lie
within a bounded distance of a geodesic ray in $X_r$.

To do that we consider for each (small) real $r$ a complex $M_r$, built out of
four flat square tori $T_i$.  On each $T_i$ we use standard angle coordinates
$(s,t)_i$ with $(s+2n\pi,t+2m\pi)_i=(s,t)_i$ for any integers $n$ and $m$.  We
let $M_r=T_1\cup T_2\cup T_3\cup T_4/\sim_r$, where $(0,t)_1\sim_r(t,0)_2$,
$(\pi,t)_2\sim_r(t,\pi)_3$, $(\pi,t)_3\sim(t,\pi)_4$ and
$(0,t)_4\sim_r(t,r)_1$.  We will let $X_r$ represent the universal cover.

The $X_r$ are Hadamard spaces with admissible fundamental groups which all have
the same geometric data.  This is easiest to see by considering $M_r$ as the
union of four spaces $T_1\cup T_2$, $T_2\cup T_3$, $T_3\cup T_4$, and $T_4\cup
T_1$ each of which is isometric to a figure eight cross a circle, where both
circles in the figure eight and the product circle have length $2\pi$.  Each of
the three spaces is glued to adjacent spaces along (product) boundary tori
(reversing the factors).  We note that the underlying finite graph of {\cal G}
is a square with four edges $e_i$ corresponding naturally to $T_i$.

It is easy to see that the $M_r$ are homeomorphic.  In fact for small $r$ the
fundamental groups are identified in a natural way.  However, we will find it
more useful to consider the map $F_r:M_0\to M_r$ defined by $F(s,t)_i=(s,t)_i$,
except that $F((0,t)_4)=(t,0)_1$.  This is not continuous along the closed
geodesic $(0,t)_4$, but the induced equivariant quasi-isometry $\tilde F:X_0\to
X_r$ is relatively easy to study.

We now want to choose a geodesic ray $\gamma$ in $X_0$.  For these spaces
geodesics are just geodesics in templates.  It is easier to first choose the
degenerate half template $\Tm$ that it will lie in.  By degenerate template we
mean a template where we alow the strip widths to be 0.  In fact in our case all
the strip widths will be 0, all the angles will be $\frac \pi 2$, and the
displacements will all be odd multiples of $\pi$.  We will choose $\Tm$ so that
the edges of the ray in the Bass-Serre tree project to the finite graph
periodically in the order $e_1$, $e_2$, $e_3$, $e_4$, $e_1\ldots$.  Subject to this
constraint we can still choose the itinerary such that the displacements in
$\Tm$ are any odd multiples of $\pi$ we please.  Hence we may choose the
itinerary (e.g. make it nontrivial) so that there is a geodesic $\gamma$ in
$\Tm$ which misses all the vertices such that $\gamma\cap W_i$ gets arbitrarily
long.  (The easiest way to see this is to consider the developement $\cal D$.  
First choose a ray $r$ that you want to be ${\cal D}(\ga)$.  
The choice of itinerary at each step amounts to a choice of 
quarter planes among those shifted by  $2n\pi$.  We can thus make 
the choice so that the quarter plane intersects $r$ in increasingly long intervals).

Of course $\gamma$ is also a geodesic ray in $X_0$.  $\tilde F$ will take $\Tm$
to a corresponding template $\Tm_r$, and the geodesic ray corresponding to
$F(\gamma)$ must be a geodesic in this template.  $\tilde F$ will be
discontinuous exactly where walls corresponding to $e_4$ and $e_1$ are glued.
The discontinuity is a translation by $r$ perpendicular to the gluing line.

We develop $\Tm$ (and $\Tm_r$) to the plane in such a way that $\gamma$ goes to
a ray with angle $\frac \pi 2>\theta>0$, the quarter planes map to planes of
type II and IV where walls of type $e_1$ and $e_3$ yield quarter planes of type
II and walls of type $e_2$ and $e_4$ yield quarter planes of type IV.  $\tilde
F$ will induce a map of developments which is discontinuous precisely on the
horizontal lines where the quarter planes coming from $e_4$ meet those coming
from $e_1$.  The discontinuity will be precisely a vertical shift by $r$.  Thus
the image of $\gamma$ consists of arbitrarily long line segments of slope
$\theta$ with infinitely many vertical jumps of size $r$.  Since the segments
get arbitrarily long, the only rays that can stay a bounded distance from
$\tilde F(\gamma)$ must also have slope $\theta$.  But then, because of the
infinitely many vertical jumps of size $r$ in $\tilde F(\gamma)$, no such ray
can stay a bounded distance from $\tilde F(\gamma)$.

\bigskip

One can make similar examples on singular piecewise Euclidean graph manifolds.
We construct such examples by gluing together two pieces.  Each piece is
topologically a twice punctured torus cross $S^1$.  The boundary of each piece
will consist of two totally geodesic square flat two-tori.  The first space is
constructed by gluing both of the corresponding boundary tori together flipping
the coordinates.  The second space is similar except that for one of the
boundary tori the gluing map is coordinate flipping composed with a small
translation.

The metric of each piece is a product metric where the circle has length 1 and
where the metric on the torus is the completion of the flat square torus minus
two line segments (slits) of length $\frac 1 2$.  The torus with the slits is a
compact flat singular space with boundary being two closed geodesics of length 1
(i.e. going ``around a slit'').

The argument that this gives an example is very similar to the above since
geodesics in the space are in fact geodesics in the corresponding templates and
since the induced map on the development of appropriate templates has properties
similar to the above example.

We also suspect there are such examples on smooth graph manifolds.  In fact one
may be able to construct such an example by smooth approximations to the above
example.  To do this carefully it would be necessary to be careful with how
closely template geodesics shadow actual geodesics in this case.

\bibliography{refs}
\bibliographystyle{alpha}
\addcontentsline{toc}{subsection}{References}

\no
Chris Croke:\\
Department of Mathematics\\
University of Pennsylvania\\
209 S. 33rd St.\\
Philadelphia, PA 19104-6395\\
ccroke@math.upenn.edu\\

\no
Bruce Kleiner:\\
Department of Mathematics\\
University of Utah\\
Salt Lake City, UT 84112-0090\\
bkleiner@math.utah.edu\\
Current address:\\
Department of Mathematics\\
University of Michigan\\
Ann Arbor, MI 48109-1109\\
bkleiner@math.lsa.umich.edu\\
\end{document}